\documentclass[a4paper,11pt,reqno]{amsart}
\pdfoutput=1
\usepackage[dvipsnames]{xcolor}
\usepackage{epsfig,url,paralist}
\usepackage{amssymb}
\usepackage[normalem]{ulem}
\usepackage{fullpage,dynkin-diagrams}
\usepackage{amsmath}
\usepackage{tikz}
\usetikzlibrary{patterns}
\usetikzlibrary{backgrounds,scopes}
\usepackage{fullpage}
\usepackage{hyperref,verbatim}
\usepackage[capitalize]{cleveref}
\usepackage{graphicx}
\usepackage{calrsfs}            
\usepackage{setspace,multirow}
\usepackage{enumitem,lineno}
\setlist{nolistsep}
\usepackage{colortbl}
\usepackage{lscape}
\usepackage{wrapfig}
\newtheorem{theo}{Theorem}

\numberwithin{equation}{section}
\newtheorem{theorem}{Theorem}[section]

\newtheorem{theorem*}{Main Result}

\newtheorem{lem}[theorem]{Lemma}
\newtheorem{coro}[theorem]{Corollary}
\newtheorem{cor*}[theorem*]{Corollary}

\newtheorem{prop}[theorem]{Proposition}
\theoremstyle{definition}
\newtheorem{remark}[theorem]{Remark}
\newtheorem{defn}[theorem]{Definition}

\newtheorem{exm}[theorem]{Example}
\numberwithin{theorem}{section}
\usepackage{fancyhdr}
\def\<{\langle}
\def\>{\rangle}
\newcommand{\PG}{\mathsf{PG}}
\DeclareMathOperator{\Res}{\mathsf{Res}}
\newcommand{\proj}{\mathsf{proj}}
\newcommand{\id}{\mathsf{id}}

\DeclareMathOperator{\ch}{\mathsf{char}}
\newcommand{\cP}{\mathcal{P}}
\newcommand{\cL}{\mathcal{L}}
\newcommand{\cM}{\mathcal{M}}
\newcommand{\K}{\mathbb{K}}
\newcommand{\FF}{\mathbb{F}}
\renewcommand{\L}{\mathbb{L}}
\renewcommand{\H}{\mathbb{H}}
\renewcommand{\O}{\mathbb{O}}
\newcommand{\A}{\mathbb{A}}
\newcommand{\B}{\mathbb{B}}
\newcommand{\PGL}{\mathsf{PGL}}

\newcommand{\SL}{\mathsf{SL}}
\newcommand{\PSL}{\mathsf{PSL}}

\newcommand{\pperp}{\perp\!\!\!\perp} 
\DeclareMathOperator{\Aut}{\mathsf{Aut}}
\DeclareMathOperator{\Norm}{\mathsf{norm}}

\DeclareMathOperator{\cotyp}{\mathsf{cotyp}}
\DeclareMathOperator{\kar}{\mathsf{char}}
\def\per{\barwedge}

\keywords{polar spaces, metasymplectic spaces, projectivities, Levi factor}
\makeatletter
\@namedef{subjclassname@2020}{\textup{2020} Mathematics Subject Classification}
\makeatother

\begin{document}

\setcounter{tocdepth}{3}

\title{Groups of Projectivities in Spherical Buildings of Non-Simply Laced Type}
\author{Sira Busch \and Hendrik Van Maldeghem}
\address{Sira Busch\\ Department of Mathematics, M\"unster University, Germany}
\email{s\_busc16@uni-muenster.de}
\address{Hendrik Van Maldeghem \\ Department of Mathematics, statistics and Computer Science, Ghent University, Belgium} \email{Hendrik.VanMaldeghem@UGent.be}
\thanks{The first author is funded by the Claussen-Simon-Stiftung and by the Deutsche Forschungsgemeinschaft (DFG, German Research Foundation) under Germany's Excellence Strategy EXC 2044 --390685587, Mathematics M\"unster: Dynamics--Geometry--Structure. The second author is partially supported by Ghent University Special Research Fund, grant BOF.24Y.2021.0029.01. This work is part of the PhD project of the first author.} 

\begin{abstract}
In this article we describe the general and special projectivity groups for all irreducible residues of all thick, irreducible, spherical buildings of type $ \mathsf{B_{n}}$, $ \mathsf{C_{n}}$ and $\mathsf{F_4}$, and rank at least $3$. This determines the exact structure and action of Levi subgroups of parabolic subgroups of groups of Lie type related to those buildings.  

\end{abstract}
\maketitle
\section{Introduction}

In \cite{Bus-Sch-Mal:24}, all general and special projectivity groups for all irreducible residues of all thick, irreducible, spherical  buildings of simply-laced type and rank at least $3$ are determined. Previously, Norbert Knarr and the second author determined these groups for many spherical buildings of rank $2$ in \cite{Kna:88} and \cite[Chapter~8]{Mal:98}, respectively. In this article, we determine the general and special projectivity groups for all irreducible residues of all thick irreducible buildings of type $ \mathsf{B_{n}}$, $ \mathsf{C_{n}}$ and $ \mathsf{F_{4}}$ and rank at least $3$. As described in \cite{Bus-Sch-Mal:24}, this will determine the exact structure and action of Levi subgroups of parabolic subgroups of simple groups of Lie type on the corresponding residues. 

Buildings of type $ \mathsf{B_{n}}$ and $ \mathsf{C_{n}}$  correspond to polar spaces, and there is a large variety of such structures. The precise special and general projectivity groups sometimes heavily depend on the field and the associated pseudo-quadratic form (see below), but we provide as much general information as possible, in particular, we provide information about generating sets that should suffice to determine the exact groups for any given situation. 

Referring for the notation and terminology to further sections, we summarise our results for type $ \mathsf{B_{n}}$ as follows. 

\begin{theo}
Let $\Delta$ be a polar space of rank $r\geq 3$, and let $U$ be a singular subspace of $\Delta$. Then the following hold.
\begin{compactenum}[$(i)$]
\item If $U$ is a projective space of dimension $d\leq r-2$, then $\Pi^+_\leq(U)$ is the full linear type-preserving group of $U$ and $\Pi_\leq(U)$ is the full linear group including (linear) dualities.
\item If $U$ is a maximal singular subspace, that is, a projective space of dimension $r-1$, then, if $\Delta$ is embeddable, $\Pi^+(U)$ is the linear group generated by homologies with factors in a certain norm set (see \emph{\cref{SPG}}), whereas $\Pi(U)$ extends $\Pi^+(U)$ with a duality having companion field automorphism the involution of the pseudo-quadratic form defining $\Delta$. If $\Delta$ is non-embeddable, then $\Pi^+(U)$ is the full  projective group of the corresponding Cayley plane and $\Pi(U)$ extends this group with a standard polarity.
\item If $U$ has dimension at most $r-3$, then $\Pi^+_\geq(U)$ is generated by products of two reflections, that is, collineations that pointwise fix a geometric hyperplane. In the cases that $\Pi^+_\geq(R)\neq\Pi_\geq(R)$, we have that $\Res_\Delta(U)$ arises from a non-degenerate quadratic form with non-degenerate associated polarity and such that $r$ minus the rank of $\Res_\Delta(U)$  is an odd integer. In that case, $\Pi_\geq(U)$ is generated by all reflections. If $\Delta$ is non-embeddable, then $U$ is a point and $\Pi^+(U)=\Pi(U)$ is the full group of direct similitudes of the corresponding quadratic form. 
\item If $U$ has dimension $r-2$, then, if $\Res_\Delta(U)$ can be identified with an embeddable polar space of rank $1$ in the natural way, the same conclusion holds as in the previous case. If $\Delta$ is non-embeddable, then $\Pi_\geq^+(U)=\Pi_\geq(U)$ is the full linear group of the corresponding projective line. 
\end{compactenum}
\end{theo}

Buildings of type $\mathsf{F_4}$ come in exactly five well specified flavours, and hence, it is possible to write down the exact special and general projectivity groups for all irreducible residues in each of these cases. We refer to \cref{tableF4} at the end of the paper for a complete enumeration.

In \cite{Bus-Sch-Mal:24}, the authors showed that, in the simply laced case, there is a general diagrammatic rule to predict when the special and general projectivity groups coincide \emph{generically} (meaning, over all fields). The present paper shows that this rule does not hold anymore when the diagram is not simply laced. In particular, for buildings of type $\mathsf{F_4}$, the rule would imply that only for the residues of type $\{2,3,4\}$, $\{2.3\}$ and all rank 1 residues,  the special and general projectivity groups generically coincide. In the split case, however, also residues of type $\{1,2,3\}$ show this behaviour, whereas we have rank $1$ residues for which the special and general projectivity groups do not coincide (but the latter was to be expected as this phenomenon also exists in the non-simply laced rank $2$ case, see \cite{Kna:88} and \cite[Chapter~8]{Mal:98} that we mentioned before). 

As a direct application of our results, we come to know the exact action of the (irreducible) Levi subgroups of the little projective group of a spherical building on the corresponding residues. (The little projective group of a spherical building is the group generated by all root elations and is usually simple.)  Indeed, for that, it suffices to generalise \cite[Proposition 3.2]{Bus-Sch-Mal:24}, which proves a connection between the little projective groups of Moufang spherical buildings of simply-laced type  and the special projectivity groups of those buildings, to the non-simply laced case. In fact the proof given in \cite{Bus-Sch-Mal:24} is valid for all Chevalley groups. However, the little projective group of a building of type $ \mathsf{B_{n}}$, $ \mathsf{C_{n}}$ or $\mathsf{F}_4$ is not always a Chevalley group (not even an algebraic group). We use a result proved in \cite{Muh-Ped-Wei:15} to prove the same in the general case, see  \cref{MR0} below. This generalises the results of \cite{Bus-Sch-Mal:24}. 

In \cite[Section~8.1]{Bus-Sch-Mal:24}, a purely algebraic method, using the root lattice and the weight lattice of a Chevalley group, is given to determine the action of the (irreducible) Levi subgroups of the little projective group of the associated spherical building on the corresponding residues, and it is applied to the simply laced case. We think that this application could be extended to the non-simply laced case when the corresponding building is associated to an adjoint Chevalley group. Perhaps this can also be done, if the building is associated to an algebraic group (but we did not try to do so); in the remaining cases (when the building is associated to a classical group over a skew field, possibly infinite-dimensional over its centre, or to a group of mixed type $\mathsf{F_4}$), it is not obvious to us, how to extend that method. In any case, one of the requirements of the method explained in \cite[Section~8.1]{Bus-Sch-Mal:24} is a certain computation in rank $2$ residues. Such a computation would be similar to the one performed in \cref{SPG}, for instance. Although certain shortcuts may arise to determine the special projectivity groups in some situations, if we extended the method of \cite[Section~8.1]{Bus-Sch-Mal:24} to the non-simply laced case, like in the split cases, we decided to treat all cases uniformly using geometric arguments. We refer the reader to \cite{Bus-Sch-Mal:24} for more details on the algebraic approach. 

The paper is structured as follows. In \cref{pre}, we introduce terminology and notation, in particular, we define the geometries (polar spaces and metasymplectic spaces) that correspond to spherical buildings of types $ \mathsf{B_{n}}$, $ \mathsf{C_{n}}$ and $\mathsf{F_4}$ and review some basic properties. In \cref{buildings} we recall some old general observations for spherical buildings that we will use in subsequent sections and make a new one. These observations are independent of the type. Then, in \cref{reductionth}, we specialise to non-simply laced types and translate some reduction theorems, proved in \cite{Bus-Sch-Mal:24}, to these types. These reduction theorems come with some conditions, and we check in \cref{prep} that the conditions are satisfied in the cases we will need. In \cref{upperPS} we determine the projectivity groups for the upper residues in polar spaces, and in \cref{projres} we do the same for the lower residues. In \cref{pointsMSS} we determine the projectivity groups of the residues of points in a metasymplectic space (which are the residues of vertices of type $1$ and type $4$ in buildings of type $\mathsf{F_4}$), and we handle the other residues of metasymplectic spaces in \cref{nonmax}. In \cref{conclusionF4} we put our results for metasymplectic spaces in a table. 

\section{Preliminaries}\label{pre}

\subsection{Point-line geometries and partial linear spaces}
Throughout, we will work with incidence structures called polar spaces and metasymplectic spaces, which are instances of partial linear spaces. In this subsection, we introduce the general definitions we will need. Note that thick polar spaces correspond to buildings of type $ \mathsf{B_{n}}$ and $ \mathsf{C_{n}}$, and (thick) metasymplectic spaces will correspond to buildings of type $\mathsf{F_4}$.

	\begin{defn}
		A \textit{point-line geometry} is a pair $\Delta=(\cP,\cL)$ with $\cP$ a set and $\cL$ a set of subsets of $\cP$. The elements of $\cP$ are called \textit{points}, the members of $\cL$ are called \textit{lines}. If $p \in \cP$ and $L \in \cL$ with $p\in L$, we say that the point $p$ \emph{lies on} the line $L$, and the line $L$ \emph{contains} the point $p$, or \emph{goes through} $p$.  If two points $p$ and $q$ are contained in a common line, they are called \textit{collinear}, denoted $p \perp q$.  If they are not contained in a common line, we say that they are \textit{non-collinear}. For any point $p$ and any subset $P \subset \cP$, we denote \[p^\perp := \{q \in \cP\mid q \perp p\} \text{ and } P^\perp := \bigcap_{p \in P} p^\perp.\]
		A\textit{(thick) partial linear space} is a point-line geometry in which every line contains at least three points, and where there is a unique line through every pair of distinct collinear points $p$ and $q$, which is then denoted with $pq$. We will usually omit the adjective ``thick''. 
	\end{defn}
	
\begin{exm}\label{example}
Let $V$ be a vector space of dimension at least $3$. Let $\cP$ be the set of $1$-spaces of $V$, and let $\cL$ be the set of $2$-spaces of $V$, each of them regarded as the set of $1$-spaces it contains. Then $(\cP,\cL)$ is called a \emph{projective space (of dimension $\dim V-1$)}, and denoted as $\PG(V)$.   
\end{exm}

	\begin{defn}
		Let $\Delta = (\cP,\cL)$ be a partial linear space.  
			\begin{enumerate}
				\item A \textit{path of length $n$} in $\Delta$ from point $x$ to point $y$ is a sequence $(x= p_0,p_1,\dots, p_{n-1},p_n=y)$ of points of $\Delta$ such that $p_{i-1} \perp p_{i}$ for all $i \in \{1,\dots,n-1\}$. It is called a \textit{geodesic} when there exist no paths of $\Delta$ from $x$ to $y$ of length strictly smaller than $n$, in which case \emph{the distance} between $x$ and $y$ in $\Delta$ is defined to be $n$, notation $d_\Delta(x,y) = n$. 
				\item The partial linear space $\Delta$ is called \emph{connected} when for any two points $x$ and $y$, there is a path (of finite length) from $x$ to $y$. If moreover the set $\{d_\Delta(x,y)\mid x,y \in \cP\}$ has a supremum in $\mathbb{N}$, this supremum is called the \textit{diameter} of $\Delta$. 
				\item A subset $S$ of $\cP$ is called a \textit{subspace} of $\Delta$ when every line $L\in\cL$ that contains at least two points of $S$, is contained in $S$. A subspace that intersects every line in at least a point, is called a \textit{hyperplane}. A hyperplane is called \emph{proper} if it does not consist of the whole point set. A subspace is called \textit{convex} if it contains all points on every geodesic that connects any pair of points in $S$. 
We usually regard subspaces of $\Delta$ in the obvious way as subgeometries of $\Delta$.  
				\item A subspace $S$ in which all points are collinear, or equivalently, for which $S \subseteq S^\perp$, is called a \textit{singular subspace}. If $S$ is moreover not contained in any other singular subspace, it is called a \textit{maximal singular subspace}. A singular subspace is called \textit{projective} if, as a subgeometry,  it is a projective space (cf.~\cref{example}). Note that every singular subspace is convex. 
				
				\item For a subset $P$ of $\cP$, the \textit{subspace generated by $P$} is denoted $\<P\>_\Delta$ and is defined to be the intersection of all subspaces containing $P$. 
 A subspace generated by three mutually collinear points, not on a common line, is called a \textit{plane}. Note that, in general, this is not necessarily a singular subspace; however we will only deal with polar and metasymplectic spaces, which implies that subspaces generated by pairwise collinear points are singular; in particular planes will be singular subspaces isomorphic to projective planes, that is, projective subspaces of dimension~2. 
			\end{enumerate}
	\end{defn}
	
	\subsection{Polar spaces} We recall the definition of a polar space, and gather some basic properties.  
	We take the viewpoint of Buekenhout--Shult \cite{Bue-Shu:74}. All results in this section are well known. 

	\begin{defn}
		A \textit{polar space} is a point-line geometry in which every line contains at least three points and for every point $x$ the set $x^\perp$ is a proper geometric hyperplane. We will only consider polar spaces of finite frank, that is, such that maximal singular subspaces are generated by a finite number of points. The minimal such number is called the \emph{rank of the polar space}. A \emph{submaximal} singular subspace is a hyperplane of a maximal singular subspace.
	\end{defn}	
	
	One can show that all singular subspaces are either empty, points, lines or  projective spaces of finite dimension (see \cite[Theorem~7.3.6 and Lemma~7.3.8]{Shu:11} or \cite[Theorem~1.3.7]{Mal:24}). Consequently, a polar space is a partial linear space. A polar space is called \emph{top-thin} if every submaximal singular subspace is contained in exactly two maximal singular subspaces. It is well known that a polar space is either top-thin, or each submaximal singular subspace is contained in at least three maximal singular subspaces (see \cite[Theorem~1.7.1]{Mal:24}). In the latter case the polar space is called \emph{thick}. In the former case, the polar space corresponds to a building of type $\mathsf{D}_n$, which have a simply laced diagram. Since these were treated in \cite{Bus-Sch-Mal:24}, we only consider the thick polar spaces. 
	
	\begin{remark} 
		If a partial linear space contains no points, or contains at least two points but no lines, it is automatically a polar space, of rank $0$ or rank $1$, respectively.   
	\end{remark}
	
\begin{defn}
Let $\Delta=(\cP,\cL)$ be a polar space, and let $\cM$ be the set of maximal singular subspaces. For each submaximal singular subspace $U$, let $P(U)$ be the set of maximal singular subspaces  containing $U$, and let $ \mathfrak{M}$ be the collection of all such sets $P(U)$ for $U$ ranging over all submaximal singular subspaces. Then $\Delta^*=(\cM,\mathfrak{M})$ is a point-line geometry that we call a \emph{dual polar space}. The \emph{dual} of a dual polar space $\Delta^*$ is the original polar space $\Delta$.
\end{defn}	

For two non-collinear points in a polar space, we call the set $\{x,y\}^{\perp\perp}=(\{x,y\}^\perp)^\perp$ a \emph{hyperbolic line}. 

\subsection{Metasymplectic spaces}	\label{metaSS}
Let $\Delta=(\cP,\cL)$ be a point-line geometry. Then we call $\Delta$ a \emph{metasymplectic space}, if the following conditions are satisfied.
\begin{enumerate}
\item If a point $p \in \Delta$ has distance $2$ to some other point $q \in \Delta$, then either there exists a unique point $u \in \Delta$, such that $p \perp u \perp q$, or $p$ and $q$ are contained in a convex subspace of $\Delta$. These convex subspaces are --- viewed as point-line geometries on their own --- isomorphic to thick polar spaces of rank $3$ and we call them \emph{symplecta}, or \emph{symps} for short. If two points $p$ and $q$ span a symplecton, it is commonly denoted by $\xi(p,q)$. Every line is contained in a symplecton.
\item For every point $p \in \Delta$, there exists at least one point $q \in \Delta$, such that $p$ and $q$ have distance $3$.
\item For every point $p$ and every symp $\xi$ not containing $p$,  the set $p^\perp\cap\xi$ is never a point. 
\end{enumerate}
	
It follows from the Main Theorem of \cite{Shu:10} that metasymplectic spaces are in one-to-one correspondence to (thick) buildings of type $\mathsf{F_4}$, together with the choice of an extremal vertex of its diagram. We summarise some more properties of such spaces, while, at the same time, introducing some more terminology. In the following, let $\Delta=(\cP,\cL)$ be a metasymplectric space. 

\begin{compactenum}[$(i)$]
\item The maximal distance between two points is $3$. Points at distance $3$ are called \emph{opposite} as they correspond to opposite vertices of the associated spherical building. If two points $p$ and $q$ are opposite, we write $p\equiv q$. Points $p,q$ at distance $2$ not contained in a symp are called \emph{special}, and we write $p\Join q$. For points $p,q$ at distance $2$ that are contained in a symp, we say that $p$ and $q$ are \emph{symplectic}, and we write $p\pperp q$. For a subset $S\subseteq\cP$, we denote by $S^{\pperp}$ the set of all points symplectic to all points of $S$.
\item Let $\Xi$ be the set of symps of $\Delta$. Singular subspaces of $\Delta$ of dimension 2 are planes. For each plane $\pi$ of a given symp, let $\Xi(\pi)$ be the set of symps of $\Xi$ containing $\pi$. Let $\Pi$ be the set of all such sets $\Xi(\pi)$, for $\pi$ ranging over all planes of all symps of $\Delta$. Then $\Delta^*=(\Xi,\Pi)$  is also a metasymplectic space. This is the \emph{principle of duality}. 
\item Under the duality mentioned in $(ii)$, collinear points correspond to symps intersecting in a plane, symplectic points correspond to symps intersecting in just a point, special points correspond to disjoint symps, for which there exists a unique symp intersecting both in respective planes and opposite symps correspond to disjoint symps with the property that being symplectic defines a bijection between their point sets inducing an isomorphism. 
\item Given a point $p\in\cP$, we define $\cL_p$ as the set of lines containing $p$ and $\Pi_p$ as the set of planes containing $p$. If we view each member $\pi$ of $\Pi_p$ as the set of all lines contained in $\pi$ that go through $p$, then it follows from the principle of duality that the geometry $\Res_\Delta(p)=(\cL_p,\Pi_p)$ is a thick dual polar space of rank $3$. We call this geometry the \emph{residue of $\Delta$ at $p$}. The set of lines through a certain point in a certain plane is called a \emph{planar line pencil}. 
\item Given a point $p$ and a symp $\xi$, there are the following three possibilities.
\begin{compactenum}[$(a)$]
\item The point $p$ belongs to $\xi$.
\item The point $p$ is collinear to (all points of) a unique line $L=p^\perp\cap\xi$ of $\xi$ and each point of $\xi$ collinear to $L$ is symplectic to $p$, while each point of $\xi$ not collinear to $L$ is special to $p$. WE say that $p$ and $\xi$ are \emph{close}.
\item The point $p$ is symplectic to a unique point $q\in\xi$ and each point of $\xi$ collinear to $q$ is special to $p$, while each other point of $\xi$ is opposite $p$.
\end{compactenum} 
\item Given two opposite points $p,q\in\cP$, let $E(p,q)$ be the set of points symplectic to both $p$ and $q$. If we denote by $\Xi(p,q)$ the set of intersections of a symp containing at least two points of $E(p,q)$ with $E(p,q)$ itself, then $(E(p,q),\Xi(p,q))$ is a polar space  isomorphic to the dual of $\Res_\Delta(p)\cong\Res_\Delta(q)$.  It is called the \emph{equator geometry (with poles $p$ and $q$)}. The poles are not necessarily unique, and the set of points which can be poles for a given equator geometry is called an \emph{imaginary line}. It follows from this  that there is a unique imaginary line through two given opposite points $p,q$ and that it coincides with $(p^{\pperp}\cap q^{\pperp})^{\pperp}$. Below, in \cref{classMSS}, we will mention precisely when an imaginary line contains more than two points. 
\end{compactenum}

\subsection{Opposition, residues, projections and groups of projectivities}
Projections and groups of projectivities in spherical buildings are defined in general, see \cref{buildings} below. In this paragraph we specialise to buildings of type $\mathsf{B}_n$ and $\mathsf{F_4}$, viewed as polar spaces and metasymplectic spaces. 

Let $\Delta$ be a polar space of rank $r$. We call two singular subspaces $U$ and $V$ \emph{opposite} if $\dim U = \dim V$ and no point of $U$ is collinear to all points of $V$. This coincides with the notion of ``opposition'' in the corresponding building.  If $d=\dim U\leq r-1$, then we consider the set of all singular subspaces of dimension $d+1$ that contain $U$ as the point set of a new geometry denoted $\Res_\Delta(U)$, or briefly $\Res(U)$, where each singular subspaces of dimension $d+2$ containing $U$ defines in the natural way a line. That geometry is a polar space of rank $r-d-1$ (empty, if $d=r-1$, of rank $1$, if $d=r-2$), called the \emph{upper residue of $U$ in $\Delta$}. The geometry with point set $U$ and lines those of $U$ itself is called the \emph{lower residue of $U$ in $\Delta$}; it is just the singular subspace $U$ viewed as a projective space, and we simply denote it by $U$, for obvious reasons. Lower residues are always projective spaces; upper residues are polar spaces.  

Let $U$ and $V$ be two opposite singular subspaces. By \cite[Theorem~1.4.11]{Mal:24}, the map $U\per V$ mapping any subspace $S\subseteq U$ to $S^\perp\cap V$ defines an isomorphism from $U$ to the dual of $V$ and is called a\emph{ (lower)} \emph{perspectivity}.  In general, we denote $S^\perp \cap V$ more systematically as $\proj_V(S)$, or, when we want to emphasize the origin, as $\proj^U_V(S)$. The composition of (a finite number of lower) perspectivities is a\emph{ lower} \emph{projectivity}. When a lower projectivity has the same image as origin, then we call it a\emph{ lower} \emph{self-projectivity}. The set of all lower self-projectivities of a subspace $U$ forms a group denoted by $\Pi_\le(U)$ and called the \emph{general lower projectivity group of $U$}. The set of lower self-projectivities, which are the composition of an even number of lower perspectivities forms a normal subgroup of $\Pi_\le (U)$ of index at most $2$, called the \emph{special lower projectivity group of $U$} and denoted as $\Pi_\le^+(U)$. If $U$ is a maximal singular subspace, then the upper residue is empty and so we denote the respective lower projectivity groups of $U$ as $\Pi^+(U)$ and $\Pi(U)$. 

We can now do the same for upper residues. Let $U$ and $V$ still be opposite singular subspaces. The map $U\per V$ (it will always be clear from te context whether lower or upper residues is meant) mapping a subspace $U'\supseteq U$ of dimension $1+\dim U$ to the unique subspace $V'\supseteq V$ of dimension $1+\dim V$ intersecting $U'$ in a point, and also denoted by $\proj_V(U')$, induces an isomorphism from $\Res(U)$ to $\Res(V)$, called an \emph{upper perspectivity}. The composition of (a finite number of upper) perspectivities is an\emph{ upper} \emph{projectivity}. When an upper projectivity has the same image as origin, then we call it an\emph{ upper} \emph{self-projectivity}. The set of all upper self-projectivities of a subspace $U$ forms a group denoted by $\Pi_\ge(U)$ and called the \emph{general upper projectivity group of $U$}. The set of upper self-projectivities, which are the composition of an even number of upper perspectivities, forms a normal subgroup of $\Pi_\ge (U)$ of index at most $2$, called the \emph{special upper projectivity group of $U$} and denoted by $\Pi_\ge^+(U)$. If $U$ is a point, then the lower residue is trivial and so, since there is no confusion, we denote the respective upper projectivity groups of $U$ as $\Pi^+(U)$ and $\Pi(U)$. 

Now let $\Delta$ be a metasymplectic space. Recall from above that  the residue of $\Delta$ at a point is always a dual polar space of rank $3$. Viewed in the dual $\Delta^*$, the residue at a point $p$ is just the rank 3 polar space given by the symplecton corresponding to $p$. For a line $L$, the \emph{upper residue}  $\Res_\Delta(L)$ is the point-line geometry with point set the set of planes of $\Delta$ containing $L$, and line set the set of sets of such planes  lying in a common given symp through $L$. This is always a projective plane, which corresponds to the \emph{lower residue} of $L$ in the dual $\Delta^*$. The \emph{lower residue} of $L$ in $\Delta$ is the projective line $L$ itself, viewed as just a set of points (hence viewed as a polar space of rank $1$). Similarly, the upper residue $\Res_\Delta(\pi)$ of a plane $\pi$ is the polar space of rank $1$ with point set the set of symplecta through $\pi$. The lower residue of $\pi$ is just the projective plane $\pi$. 

 Let $p,p'$ be two opposite points. Then for each symplecton $\xi$ containing $p$, there exists a unique symplecton $\xi'$ through $p'$ intersecting $\xi$ in a unique point. The map $p\per p'$ mapping $\xi\mapsto\xi'$ induces an isomorphism of $\Res_\Delta(p)$ to $\Res_\Delta(p')$, called a \emph{perspectivity}. We define in exactly the same way as above \emph{projectivities}, \emph{self-projectivities} and the \emph{special} and \emph{general} projectivity group of $p$, denoted $\Pi^+(p)$ and $\Pi(p)$, respectively. 
 
Now let $L,L'$ be two opposite lines, that is, each point $x$ of $L$ is not opposite a unique point $x'$ of $L'$ and vice versa. The map $L\per L'$ given by $x\mapsto x'$ is again a \emph{lower perspectivity}, and we similarly as before define \emph{lower projectivities}, \emph{lower self-projectivities}, and the \emph{lower special} and  \emph{general projectivity groups} of $L$, denoted as $\Pi^+_\leq(L)$ and $\Pi_\leq(L)$, respectively. For each symplecton $\xi$ containing $L$, there exists a unique plane $\pi'$ through $L'$ all points of which are close to $\xi$. The map $L\per L'$ (again, it will always be clear from the context whether this is between lower or upper residues) given by $\xi\mapsto \pi'$ induces an isomorphism from $\Res_\Delta(L)$ to $\Res_\Delta(L')$ (not preserving types), called an \emph{upper perspectivity}. Again this gives rise to \emph{upper projectivities}, \emph{upper self-projectivities}, and the \emph{upper special} and  \emph{general projectivity groups} of $L$, denoted as $\Pi^+_\geq(L)$ and $\Pi_\geq(L)$, respectively. For a plane $\pi$ of $\Delta$, the \emph{upper (lower) special} and \emph{general projectivity groups} of $\pi$ in $\Delta$, denoted $\Pi_\geq^+(\pi)$ ($\Pi_\leq^+(\pi)$) and $\Pi_\geq^+(\pi)$ ($\Pi_\leq^+(\pi)$), respectively, are the lower (upper) special and  general projectivity groups of $\pi$ as a line in the dual $\Delta^*$.
 
 Now let $p,p'$ be opposite points and $\xi,\xi'$ opposite symps, with $p\in\xi$ and $p'\in\xi'$. Let $L$ be a line with $p\in L\subseteq\xi$. There is a unique line $L''\in\xi'$, which is the projection of $L$ onto $\xi'$, and there is a unique line $L'\ni p'$ intersecting $L''$ in some point. It follows from  \cite[3.19.5]{Tits:74} that the map $\{p,\xi\}\per\{p',\xi'\}$ given by $L\mapsto L'$ induces an isomorphism from $\Res_\Delta(\{p,\xi\})=:\Res_\xi(p)$ to $\Res(\{p',\xi'\})=:\Res_{\xi'}(p')$, called a \emph{perspectivity}. This defines, as before, \emph{(self-)projectivities} and the \emph{special} and \emph{general projectivity group} of $\{p,\xi\}$, denoted as $\Pi^+(\{p,\xi\})$ and $\Pi(\{p,\xi\})$, respectively. 
 
 Finally, let $P$ be a planar line pencil in a polar space or a metasymplectic space. Such a pencil is determined by its vertex $x$ and a plane $\pi$ and consists of all lines through $x$ in $\pi$.  Similarly as in the previous paragraph one defines \emph{perspectivities} $P\per P'$ and {(self-)projectivities}, and the \emph{special} and \emph{general projectivity groups} of $P$, which we denote by $\Pi^+(P)$ (or $\Pi^+(\{x,\pi\})$) and  $\Pi(P)$ (or $\Pi(\{x,\pi\})$), respectively, and which are permutation groups of $P$. 

\subsection{Classification of polar spaces}
Polar spaces of rank at least $3$ have been classified and we will give an overview in the following. For that purpose, recall from \cite[Section~8.3]{Tits:74} (see also \cite[Section~3.2]{Mal:24}) that a \emph{polarity} of a projective space is a symmetric relation $\Box$ between the points, such that, (with obvious notation,) $p^{\Box}$ is a hyperplane for each point $p$. A polarity is non-degenerate, if $p^{\Box}$ is always a proper hyperplane. For a subspace $S$, we define $S^{\Box}$ to be the intersection of all $x^{\Box}$, for $x \in S$. An \emph{absolute subspace} for $\Box$ is a subspace $U$ with the property that $U\subseteq U^{\Box}$. The \emph{radical} of a polarity $\Box$ is the image under $\Box$ of the full point set. The radical is then the intersection of all $x^{\Box}$, for $x$ varying over all points. A polarity is called a \emph{null-polarity}, if each point is absolute. Most polar spaces we will review are embeddable into some projective space. We also include the rank $1$ and $2$ polar spaces that arise as residues in higher rank polar spaces. In the following, we assume that the reader is familiar with the basic algebraic notions of bilinear and quadratic forms. The null set of a quadratic form is a \emph{quadric}.  If a subspace $U$ is disjoint from that quadric, then the quadratic form is said to be \emph{anisotropic over $U$}. For a $\K$-vector space and bilinear form $b \colon V \times V\to\K$, the \emph{associated polarity $\Box$} is defined as $p\Box q$ if $b(v,w)=0$, where $p$ and $q$ are the $1$-spaces of $V$ generated by $v$ and $w$, respectively. 
 
\subsubsection{Symplectic polar spaces} These are polar spaces arising from non-degenerate null-polarities, that is, non-degenerate polarities in finite-dimensional projective spaces, such that all points of the projective space are absolute points. The dimension of the projective space is always odd, say $2n-1$, and the rank of the polar space is $n$.  The coordinatising field is always commutative. These polar spaces are completely and uniquely determined by giving the collinearity relation between arbitrary points of the projective space. In standard form, this is given by the alternating bilinear form \[f((x_{-n},\ldots,x_{-1},x_1,\ldots,x_n),(y_{-n},\ldots,y_{-1},y_1,\ldots,y_n))=
x_{-n}y_n-y_{-n}x_n+\cdots +x_{-1}y_1-y_{-1}x_1.\]
  
\subsubsection{Orthogonal polar spaces} These are polar spaces arising from non-degenerate quadratic forms of finite Witt index at least $1$. Recall that a quadratic form is non-degenerate, if it is anisotropic over the radical of the associated bilinear form. In characteristic different from $2$, this just means that the associated bilinear form is non-degenerate (has trivial radical). The polarity associated to the bilinear form is also non-degenerate.  In characteristic $2$, however, we distinguish between non-degenerate quadratic forms with degenerate associated bilinear form (and call them \emph{inseparable}), and those with non-degenerate associated bilinear form (and call them \emph{separable}). We extend this terminology to the other characteristics by calling every non-degenerate quadratic form in characteristic not 2 \emph{separable}. We further extend this terminology to the orthogonal polar spaces themselves in the obvious way. In characteristic $2$, the associated polarity of a separable orthogonal polar space is a non-degenerate null-polarity. The one of an inseparable orthogonal polar space is a degenerate null-polarity.  Separable orthogonal polar spaces admit, up to isomorphism, a unique embedding in some projective space. The difference of the dimension of the ambient projective space with twice the rank, plus one, will be called the \emph{anisotropic corank} of the polar space. It is the (maximal) vector dimension of an anisotropic form that is needed to describe the corresponding quadric.    Quadrics with anisotropic corank equal to $0$ are \emph{hyperbolic quadrics} and correspond to buildings of type $\mathsf{D}_n$. If the anisotropic corank is $1$, we speak about \emph{parabolic quadrics} and \emph{parabolic polar spaces}.

Since hyperbolic quadrics are related to buildings of type $\mathsf{D}_n$, and we determined the projectivity groups of buildings of type $\mathsf{D}_n$ in \cite{Bus-Sch-Mal:24}, we will not be much concerned with them in this article, except that they will show up in some proofs, when we extend the ground field (which we explain now in some more detail). According to Theorems~3.4.3 and~4.4.4 of \cite{Mal:24}, the standard equation of a quadric in $\PG(V)$, with $V$ a $\K$-vector space, that corresponds to a polar space of rank $r$, is
\[X_{-r}X_r+\cdots+X_{-2}X_2+X_{-1}X_1=f_0(v_0,v_0),\]
where $V=V'\oplus V_0$, with $\dim V'=2r$, $\{e_{-r},\ldots,e_{-1},e_1,e_2,\ldots,,e_r\}$ is a basis of $V'$, and a vector $v\in V$ is written as $v=x_{-r}e_{-r}+\cdots+x_{-1}e_{-1}+x_1e_1+x_2e_2+\cdots x_re_r+v_0$, with $v_0\in V_0$, and $f_0$ is an anisotropic quadratic form in $V_0$. In the finite dimensional case, if $\dim V_0$ is even, then the above equation becomes the equation of a hyperbolic quadric over a \emph{splitting field}, that is, an overfield of $\K$ over which $f_0$ becomes completely reducible (one can take the algebraic or quadratic closure of $\K$). We will use this and note that hyperbolic quadrics have two natural systems of maximal singular subspaces. Two maximal singular subspaces belong to the same system if and only if the intersection has even codimension in each of them (the \emph{codimension of $U\cap V$ in $U$} is $\dim U-\dim(U\cap V)$). 

\def\L{\mathbb{L}}
\def\q{\mathfrak{q}}

\subsubsection{Pseudo-quadratic polar spaces} These are polar spaces arsing from $\sigma$-quadratic forms, where $\sigma$ is an involution of a skew field $\L$ (and $\sigma\neq\id$). We briefly describe these. Let $V$ be a right vector space over $\L$. Let $g:V\times V\to\L$ be a $(\sigma,\id)$-linear form (meaning $g$ is additive in both variables and $g(vk,w\ell)=k^\sigma g(v,w)\ell$, for all $k,\ell\in\L$ and $v,w\in V$). We define the following:
\begin{align*}
\L_\sigma &:=  \{t-t^\sigma\mid t\in\L\} \\
f\colon & V\times V\to\L, ~ (v,w)\mapsto g(v,w)+g(w,v)^\sigma \\
\mathfrak{q}\colon & V\to\L/\L_{\sigma}, ~ v\mapsto g(v,v)+\L_\sigma
\end{align*}

Then $f$ is a Hermitian sesquilinear form (meaning it is $(\sigma,\id)$-linear and $f(v,w)^\sigma=f(w,v)$) and we denote its radical by $R_f$, that is, $R_f=\{v\in V\mid f(v,w)=0,\forall w\in V\}$. We say that $f$ \emph{is associated to} $\q$. Suppose $\mathfrak{q}$ is anisotropic over $R_f$, that is, $\mathfrak{q}(v)=0$ if, and only if, $v=0$, for all $v\in R_f$. Let $X$ be the set of vectors $v\in V$ for which $\mathfrak{q}(v)=0$. Suppose the subspaces of $V$ of maximal dimension, that are contained in $X$, have finite dimension $r$. Then the point set $\cP=\{\<v\>\mid v\in X\}$ of $\PG(V)$, together with the lines induced from $\PG(V)$, is a  polar space of rank $r$, called a \emph{pseudo-quadratic polar space}. This description also makes sense for $\sigma=1$ (and all embeddable polar spaces, except for the symplectic ones in characteristic distinct from $2$, can be described like this). We call it the \emph{pseudo-quadratic description}. 

In the commutative case, the Hermitian sesquilinear form $f$ determines the polar space (without using the associated pseudo-quadratic form; its points correspond to vectors $v\in V$ for which $f(v,v)=0$) and we call the polar spaces \emph{Hermitian}. Sometimes the form $V\to V$ with $v\mapsto f(v,v)$ is called a \emph{Hermitian form}. 

\subsubsection{Non-embeddable polar spaces} \label{NEPS}

There are two different kinds of non-embeddable polar spaces. Non-embeddable polar spaces of one kind are top-thin and arise as the line-Grassmannian of a $3$-dimensional projective space over a non-commutative skew field. Since the projectivity groups of projective spaces are all well known (see for instance \cite{Bus-Sch-Mal:24}), we will not need to consider these in the sequel. 

The other non-embeddable polar spaces are thick and each of them arises as the fixed point structure of an involution in a building of type $\mathsf{E_7}$. In fact, by the main result of \cite[Chapter~9]{Tits:74}, for each non-Desarguesian Moufang plane $\pi$, there exists a unique polar space of rank $3$, whose planes are isomorphic to $\pi$. We will not need an explicit construction of those polar spaces and refer to \cite{Bru-Mal:14} for an elementary one. We will use some properties of these polar space derived in \cite{Pas-Mal:23}, where thick non-embeddable polar spaces are called \emph{Freudenthal-Tits polar spaces}. However, we will need the standard description of the Moufang plane associated to a Cayley division algebra $\O$ over a field $\K$. We denote that plane as $\PG(2,\O)$. Recall that $\O$ is an $8$-dimensional, alternative, quadratic, non-associative composition division algebra with a standard involution $x\mapsto\overline{x}$, such that, for all $x\in \O$, we have $x\overline{x}\in\O$ and $x+\overline{x}\in\O$. An affine plane can be described as the set of ordered pairs $(x,y)\in\O\times\O$, with lines of two types: (1) For $m,k\in\O$, the line $[m,k]$ contains the points $(x,mx+k)$, for all $x\in\O$; (2) for $x\in\O$, the line $[x]$ contains the points $(x,y)$, for all $y\in\O$.  We will need this description in \cref{planeNE}. 

\subsection{Classification of metasymplectic spaces}\label{classMSS}

According to \cite[Chapter~10]{Tits:74}, buildings of type $\mathsf{F_4}$ are in one-to-one correspondence with pairs $(\K,\A)$, where $\K$ is a field and $\A$ is a quadratic, alternative (composition) division algebra over $\K$. Usually, one labels the $\mathsf{F_4}$-diagram linearly in such a way that the type set $\{1,2\}$ corresponds to the residues isomorphic to $\PG(2,\K)$, and $\{3,4\}$ to those isomorphic to $\PG(2,\A)$ (this is also called the \emph{Bourbaki labelling}). We denote such building as $\mathsf{F_{4}}(\K,\A)$. The metasymplectic spaces $\mathsf{F_{4,1}}(\K,\A)$ take the vertices of type $1$ as points, and hence, its planes are defined over $\K$, whereas $\mathsf{F_{4,4}}(\K,\A)$ takes the vertices of type $4$ as points and has planes defined over $\A$. Quadratic alternative division algebras come in exactly five flavours. Note in advance that each such algebra is equipped with a \emph{standard involution}, that is, and involutive anti-automorphism $\sigma$, such that $xx^\sigma\in\K$ and $x+x^\sigma\in\K$. The form $x\mapsto xx^\sigma$ is the \emph{norm form}. 
\begin{compactenum}
\item  $\A=\K$. Then $\mathsf{F_4}(\K,\K)$ is the so-called \emph{split building of type $\mathsf{F_4}$}. The symps of $\mathsf{F_{4,1}}(\K,\K)$ are parabolic polar spaces over $\K$; those of $\mathsf{F_{4,4}}(\K,\K)$ are symplectic polar spaces over $\K$.
\item $\A$ is a separable quadratic extension of $\K$.  The symps of $\mathsf{F_{4,1}}(\K,\A)$ are orthogonal polar spaces with anisotropic corank $2$. The anisotropic quadratic form needed to describe those, is exactly the norm form of the Galois extension. The symps of $\mathsf{F_{4,4}}(\K,\A)$ are Hermitian polar spaces, naturally embedded into $\PG(5,\A)$, and associated to the standard involution of $\A$ (which is the Galois involution).
\item $\A$ is a quaternion division algebra over $\K$.  The symps of $\mathsf{F_{4,1}}(\K,\A)$ are orthogonal polar spaces with anisotropic corank $4$. The anisotropic quadratic form needed to describe these is exactly the norm form of the quaternion division algebra. The symps of $\mathsf{F_{4,4}}(\K,\A)$ are pseudo-quadratic polar spaces, naturally embedded into $\PG(5,\A)$, and associated to the standard involution $\sigma$ in $\A$. 
\item $\A$ is a Cayley division algebra (or octonion division algebra) over $\K$. The symps of $\mathsf{F_{4,1}}(\K,\A)$ are orthogonal polar spaces with anisotropic corank $8$. The anisotropic quadratic form needed to describe these is exactly the norm form of the Cayley division algebra. The symps of $\mathsf{F_{4,4}}(\K,\A)$ are the thick non-embeddable polar spaces associated to $\A$. 
\item $\A$ is an inseparable extension of $\K$ in characteristic $2$. The symps of both $\mathsf{F_{4,1}}(\K,\A)$ and $\mathsf{F_{4,4}}(\K,\A)$ are inseparable orthogonal polar spaces.  
\end{compactenum}
In general, we denote a symp of $\mathsf{F_{4,1}}(\K,\A)$ as $\mathsf{B_{3,1}}(\K,\A)$, and one of $\mathsf{F_{4,4}}(\K,\A)$ as $\mathsf{C_{3,1}}(\A,\K)$. The corresponding dual polar spaces, which are point residues in $\mathsf{F_{4,4}}(\K,\A)$ and $\mathsf{F_{4,1}}(\K,\A)$, respectively, are denoted as $\mathsf{B_{3,3}}(\K,\A)$, and  $\mathsf{C_{3,3}}(\A,\K)$, respectively.

In \cref{metaSS} above, we defined the equator geometry of two opposite points of a metasymplectic space. We can now define the extended equator geometry for the metasymplectic space $\mathsf{F_{4,4}}(\K,\A)$.  Let $p,q$ be two opposite points of $\mathsf{F_{4,4}}(\K,\A)$. Then the union of all equator geometries $E(x,y)$, for $x,y$ opposite points varying over $E(p,q)$, together with all lines in these geometries, and in $E(p,q)$, is called the \emph{extended equator geometry}, denoted as $\widehat{E}(p,q)$. Note that is contains $p$ and $q$. It follows from \cite{Pet-Mal:23} that $\widehat{E}(p,q)$ is a polar space whose residues are isomorphic to $E(p,q)$ (which are isomorphic themselves to $\mathsf{B_{3.1}}(\K,\A)$, and hence, we denote such polar space as $\mathsf{B_{4,1}}(\K,\A)$).  

We now take a closer look at imaginary lines. Let $p,q$ be two opposite points. By \cite[Proposition~2.10.5]{Lam-Mal:24}, the imaginary line through $p$ and $q$ coincides with the hyperbolic line through $p$ and $q$ in the polar space $E(x,y)$, for any pair of opposite points $x,y$ in $E(p,q)$. Hence, imaginary lines have size $2$ in $\mathsf{F_{4,4}}(\K,\A)$ when $\A$ is not an inseparable extension of $\K$, including the case $\A=\K$ when $\kar\K=2$. In alll other cases, imaginary lines have size at least $3$.


\section{Some general observations for buildings of non-simply laced type}\label{buildings}
The present section is the only section, where we use some pure building theory. We briefly introduce the notions that we will need to use. For more background we refer to the book by Abramenko \& Brown \cite{Abr-Bro:08}. 

A (thick) spherical building is a (thick) simplicial chamber complex, such that every pair of simplices is contained in a finite thin chamber subcomplex, called an \emph{apartment}, and every two apartments are isomorphic through an isomorphism fixing every vertex in their intersection. A \emph{panel} is a simplex, which can be completed to a chamber by properly adding precisely one vertex. \emph{Adjacent chambers} are chambers sharing a panel. The adjacency graph on the set of chambers induces a numerical distance functiion on the set of chambers. For each pair of adjacent chambers $C,C'$ of an apartment there exists a unique \emph{folding}, that is, an idempotent morphism with the property that every chamber is the image of zero or precisely two chambers, mapping $C'$ to $C$. The image of a folding is called a \emph{root}. If $\alpha$ is the root determined by the folding $C'\mapsto C$, then we denote by $-\alpha$ the \emph{opposite root}; namely the one determined by the opposite folding $C\mapsto C'$. The intersection $\alpha\cap(-\alpha)$ is called the \emph{boundary} of both $\alpha $ and $-\alpha$ and denoted $\partial\alpha$. The \emph{interior} of $\alpha$ is $\alpha\setminus\partial\alpha$. 

For a simplex $F$, the residue $\Res(F)$ is the simplicial complex induced on the set of vertices $v\notin F$ such that $F\cup\{v\}$ is a simplex. If we work in a building $\Delta$, then, for clarity, $\Res(F)$ is sometimes also denoted by $\Res_\Delta(F)$. The \emph{type} of $F$ is the set of types of its members; the \emph{cotype} is the the complementary set of types (with respect to the whole type set). The \emph{type} of the residue $\Res(F)$ is the cotype of $F$. 
If $F$ and $F'$ are opposite simplices, then for each chamber $C$ containing $F$ there exists a chamber $C'$ containing $F'$ at minimal distance. The mapping $C\setminus F\mapsto C'\setminus F'$ induces an isomorphism from $\Res(F)$ to $\Res(F')$ (see \cite[Theorem~3.28]{Tits:74}), which we call a \emph{perspectivity}. As in the previous section, this gives rise to the notions \emph{(even) projectivity, self-projectivity}, the \emph{special projectivity group $\Pi^+(F)$}, and the \emph{general projectivity group $\Pi(F)$}.  

Chambers in an apartment are called \emph{opposite}, if they are never contained in the same root. This is independent of the chosen apartment and therefore, we say that chambers are opposite in a building, if they are opposite in an apartment. If the building is thick, then this is equivalent to saying that there exists a unique apartment containing the two chambers. All spherical buildings of rank at least $3$ are \emph{Moufang}, that is, for each root $\alpha$, the group $U_\alpha$ fixing every chamber having a panel in the interior of $\alpha$, acts transitively on the set of apartments containing~$\alpha$. 

The aim of this section is to extend \cite[Theorem A]{Bus-Sch-Mal:24} to spherical Moufang buildings, which are not necessarily of simply laced type. We will prove the following theorem.

\begin{theorem}\label{MR0}
Let $F$ be a simplex of a Moufang spherical building $\Delta$. Let $\Aut^+(\Delta)$ be the automorphism group of $\Delta$ generated by the root groups. Then $\Pi^+(F)$ is permutation equivalent to the action of the stabiliser $\Aut^+(\Delta)_F$ of $F$ in $\Aut^+(\Delta)$ on the residue $\Res_\Delta(F)$ of $F$ in $\Delta$. 
\end{theorem}

The proof of this theorem, given in \cite{Bus-Sch-Mal:24} for the simply laced case,  requires that the unipotent radical of a parabolic subgroup in a Moufang spherical building pointwise stabilises the corresponding residue, and acts transitively on the simplices opposite the given residue. For the simply laced case, this follows from the Levi decomposition of parabolic subgroups in Chevalley groups. In general, we can use \cite[Proposition~24.21]{Muh-Ped-Wei:15}: 
\begin{lem}\label{L1}
Let $\Delta$ be a spherical Moufang building and let $F$ be a simplex of $\Delta$. Let $G_{F}$ be the stabiliser of $F$ in $\Aut^+(\Delta)$. Then there exists a subgroup $U_F\leq G_{F}$ which acts sharply transitively on the set $F^\equiv$ of simplices opposite $F$, and which pointwise fixes $\Res_\Delta(F)$.  
\end{lem}

For completeness we describe how the subgroup $U_F$ is constructed. 

Let $F'$ be a simplex in $\Delta$ opposite $F$.
Choose an apartment $\Sigma$ containing $F$ and $F'$ and a chamber $C$ in $\Sigma$ containing $F$.
For a root $\alpha$ of $\Sigma$, let $U_{\alpha}$ be the corresponding root group; that is the group of automorphisms $g$ of $\Delta$ that fixes all chambers that have a panel in the interior of $\alpha$.

Then $U_F$ is the group generated by all root groups $U_{\alpha}$ corresponding to roots $\alpha$, such that $\alpha$ contains $C$, and $F$ is in $\alpha$, but not in the boundary $\partial \alpha$ of $\alpha$. 

Now exactly the same arguments as in \cite{Bus-Sch-Mal:24} lead to \cref{MR0}.
We also recall the following general rule, see \cite[Observation 3.1]{Bus-Sch-Mal:24}.

\begin{prop}\label{O1}
Let $\Delta$ be a spherical building over the type set $I$ and let $J\subseteq I$ be self-opposite. Let $F$ be a simplex of type $J$. Then $\Pi^+(F)=\Pi(F)$ if, and only if,  the identity in $\Pi(F)$ can be written as the product of an odd number of perspectivities. 
\end{prop} 

We will also use \cite[Lemma 7.1]{Bus-Sch-Mal:24}, which we state now.

\begin{lem}\label{gate2}
Let $\Delta$ be a spherical building over the type set $I$ and let $F_K$ be  a simplex of type $K\subseteq I$.  Let $K\subseteq J\subset I$ and let $F_J$ be a simplex of type $J$ containing $F_K$. Let $\Pi^+_K(F_J)$ be the special projectivity group of $F_J\setminus F_K$ in $\Res_\Delta(F_K)$. Then $\Pi^+_K(F_J)\leq \Pi^+(F_J)$.
\end{lem}

The following result follows directly from the fact that two chambers are always contained in an apartment, and each chamber of an apartment has a unique opposite chamber in that apartment. 

\begin{lem}\label{Opp-notopp}
In a thick spherical building $\Delta$, given a pair $(C,C')$ of distinct chambers, there exists a chamber $D$ opposite $C$ and not opposite $C'$.
\end{lem}

Finally we recall \cite[Proposition~8.2]{Bus-Sch-Mal:24}. A building is said to have \emph{thickness at least $t+1$}, if every panel is contained in at least $t+1$ chambers. 

\begin{prop}\label{Tits3.30}
If a spherical building has thickness at least $t+1$, then there exists a chamber
opposite $t$ arbitrarily given chambers. In particular, there exists a vertex opposite $t$ arbitrarily given vertices of the same self-opposite type.
\end{prop}


\section{General reduction theorems}\label{reductionth}
In this section, we prove some results that reduce the computation of the projectivity groups to rather special cases. We also establish when the special and general projectivity groups generically coincide. We begin with the latter.

\subsection{Special versus general projectivity groups in polar spaces}
It is clear that the general projectivity group of the lower residue of a singular subspace contains dualities and the special group does not. Hence, these are always different. 

Concerning the projectivity groups of the upper residues, the almost completely opposite situation holds. Indeed, we have the following result. 

\begin{prop}\label{special=general}
Let $\Delta$ be a thick polar space of rank at least $3$. Let $U$ be a non-maximal singular subspace of odd dimension at least $1$. Then $\Pi^+_\geq(U)\equiv\Pi_\geq(U)$. The same conclusion holds, if $\Delta$ is not a separable orthogonal polar space and $U$ has arbitrary dimension (but is non-maximal). 
\end{prop}

\begin{proof}
Pick a singular subspace $U'$ opposite $U$. If $\Delta$ is non-embeddable, then $U$ and $U'$ are points or lines and $(U^\perp\cap U'^\perp)^\perp$ is a thick polar space $\Gamma$ of rank $1$ or $2$, respectively (this follows from Proposition 5.11 of \cite{Pas-Mal:23}). Hence, there exists a point or line~$U''$ in~$\Gamma$ opposite both~$U$ and~$U'$. The projectivity of upper residues $U\per U'\per U''\per U$ is the identity, showing the assertion in this case.  

Now suppose that $\Delta$ is embedded in $\PG(V)$, with $V$ minimal; that is, the associated polarity is non-degenerate. In this case, the subspace $\Gamma$ of $\Delta$, induced by the subspace of $\PG(V)$ spanned by~$U$ and~$U'$, is a polar space distinct from a hyperbolic space of odd rank. Therefore, $\Gamma$ contains a singular subspace $U''$ opposite both $U$ and $U'$. 
Since the associated polarity in $\PG(V)$ is non-degenerate, we have $U''\subseteq(U^\perp\cap U'^\perp)^\perp$ and so $U\per U'\per U''\per U$ is the identity. Now the assertion follows from \cref{O1}. 
\end{proof}

\subsection{Special versus general projectivity groups in metasymplectic spaces}\label{specialgeneralMSS}
\begin{prop}\label{special=generalMSS}
Let $p$ be a point of a metasymplectic space $\Gamma$. Then $\Pi(p)=\Pi^+(p)$ as soon as there exists an imaginary line in $\Gamma$ of size at least $3$ and containing $p$.
\end{prop}

\begin{proof}
Let $p,q,r$ be three points of the same imaginary line. Then all points symplectic to both $p$ and $q$ are also symplectic to $r$, which implies that $p\per q\per r\per p$ fixes each symplecton through $p$, and is hence the identity. Now apply \cref{O1}
\end{proof}

From \cite[Lemma~5.2]{Bus-Sch-Mal:24} it now follows, for a simplex $S$ of type $\{1\},\{1,4\}$ or $\{1,3,4\}$ of $\mathsf{F_4}(\K,\A)$, with $\K$ a field and $\A$ an alternative quadratic division algebra over $\K$, that $\Pi(F)=\Pi^+(F)$. For simplices of type $\{1,4\}$ we will come back to this in an explicit way in the proof of  \cref{14F4}. 

Also, note that the condition stated in \cref{special=generalMSS} is not necessary, as for split buildings the conclusion will hold without the condition being satisfied, see the first line of row (B3) in \cref{tableF4}.

\subsection{Reduction to the product of three perspectivities}
In principle, to determine the projectivity groups, one has to consider arbitrarily long sequences of perspectivities. However, the following results will lead to the fact that projectivity groups are generated by self-projectivities which are products of at most four perspectivities in a particular sequence.

We say that a set $\Pi$ of automorphisms of a polar space $\Delta$ is \emph{geometric}, if its members are characterised by their fix structure. Formally, this means that an automorphism belongs to $\Pi$ if, and only if, its fix set is a member of a certain given set of subsets of the point set of $\Delta$, closed under the action of the full automorphism group of $\Delta$. We will mainly apply this notion for fix structures being hyperplanes or subhyperplanes.  We now phrase \cite[Lemma~8.1]{Bus-Sch-Mal:24} for our situation of polar spaces.

 \begin{lem}\label{triangles}
 Let $\Delta$ be a polar space of rank $r$ and let $j$, $0\leq j<r$, be an arbitrary natural number.  If $j>0$, then suppose that for each quadruple of singular subspaces of dimension $j$ containing at least one opposite pair, there exists a singular subspace of dimension $j$ opposite all the members of the given quadruple. If $j=0$, suppose the same conclusion holds for each quadruple of points with the property that the pairwise intersections of the perps are not all the same. Let $F, F', F''$ be three pairwise opposite singular subspaces of dimension  $j$ and denote by $\theta_0$ the projectivity $F\per F'\per F''\per F$ of upper residues, if $j<r-1$, and of lower residues, if $j=r-1$. Denote by $\Pi_3(F)$ the set of all corresponding self-projectivities of $F$ of length $3$ and suppose that $\Pi_3(F)$ is geometric. 
 Then $\Pi(F)=\<\Pi_3(F)\>$ and $\Pi^+(F)=\<\theta_0^{-1}\theta\mid \theta \in\Pi_3(F)\>$. 
 \end{lem}

The conditions on the quadruple of singular subspaces do not appear in \cite[Lemma~8.1]{Bus-Sch-Mal:24}. However, in the proof, each considered quadruple $(F_0,F_1,F_2,F_3)$ is part of a chain of perspectivities $F_0\barwedge F_1\barwedge F_2\barwedge F_3$. Hence, $F_0$ is opposite $F_1$, which is opposite $F_2$, and $F_2$ is opposite $F_3$. Also, if $j=0$, and if $p_0,p_1,p_2$ are points, such that $p_0^\perp\cap p_1^\perp=p_1\cap p_2^\perp$, then $p_0\barwedge p_1\barwedge p_2=p_0\barwedge p_2$. So we may shorten the chain without altering the projectivity defined by the chain.  

\subsection{Reduction to the product of four perspectivities}

Now we phrase \cite[Lemma~8.17]{Bus-Sch-Mal:24} in terms of polar spaces and our situation. In the subsequent remark, we improve on the conditions. But first a definition.

 \begin{defn}\label{defsim}
Let $\Delta$ be a polar space of rank $n \geq 2$. Let $S$ be a set of singular subspaces of dimension $s$ in $\Delta$. We define the \textit{$s$-space-graph} $\Gamma(S)=(S,\sim)$, with $\sim$ denoting the adjacency, as follows:
\begin{enumerate}
\item[(V)] The vertices are the elements of $S$.
\item[(E)]
\begin{itemize}
\item For $s=0$, we draw an edge between two vertices in $\Gamma(S)$, if the corresponding points are collinear in $\Delta$.
\item For $s \in [1, n-2]$, we draw an edge between two vertices in $\Gamma(S)$, if the corresponding subspaces in $\Delta$ intersect in a subspace of dimension $s-1$ and are both contained in a common subspace of dimension $s+1$.
\item For $s = n-1$, we draw an edge between two vertices in $\Gamma(S)$, if the corresponding maximal subspaces in $\Delta$ intersect in a subspace of dimension $s-1$.
\end{itemize}
\end{enumerate} 
\end{defn}

\begin{lem}\label{upanddown}
 Let $\Delta$ be a polar space of rank $r$ and let $j$, $0\leq j<r$, be an arbitrary natural number.  Suppose that for each pair of singular subspaces $F$, $F'$ of dimension $j$, the graph $\Gamma(S)$, where $S$ is the set of singular subspaces opposite both $F$ and $F'$, is connected. Suppose also that there exists a singular subspace opposite any given set of three singular subspaces of dimension $j$. Let $F$ be a given singular subspace of dimension $j$. Denote by $\Pi_4(F)$ the set of all self-projectivities $F\per F_2\per F_3\per F_4\per F$ of $F$ of length $4$ with $F\sim F_3$, $F_2\sim F_4$. Suppose that $\Pi_4(F)$ is geometric. Then $\Pi^+(F)=\<\Pi_4(F)\>$. 
 \end{lem}

\begin{remark}\label{openneer}
In \cref{upanddown} we may also assume that $F\cap F_3$ and $F_2\cap F_4$ are opposite. Indeed, if not, then we can select in $\<F,F_3\>$ a subspace $F_3'$ of dimension $j$ disjoint from $\proj_F(F_2\cap F_4)$. We then write \[F\barwedge F_2\per F_3\per F_4=(F\per F_2\per F_3'\per F_4\per F)\cdot(F\per F_4)\cdot (F_4\per F_3'\per F_2\per F_3\per F_4)\cdot F_4\per F.\]

Likewise, if $j\leq r-2$, we may also assume that $\<F,F_3\>$ and $\<F_2,F_4\>$ are opposite. Indeed, noting that a singular subspace $U_1$ through $F\cap F_3$ is opposite a singular subspace $U_2$ through $F_2\cap F_4$ if, and only if, $U_1\cap (F\cap F_3)^\perp\cap (F_2\cap F_4)^\perp$ is opposite $U_2\cap (F\cap F_3)^\perp\cap (F_2\cap F_4)^\perp$ in the polar space $ (F\cap F_3)^\perp\cap (F_2\cap F_4)^\perp$ (where we may assume, due to the previous paragraph, that $F\cap F_3$ and $F_2\cap F_4$ are opposite), it suffices to verify the claim for $i=0$. 

First suppose the line $L_1:=\<F,F_3\>$ intersects the line $L_2:=\<F_2,F_4\>$ in a point $p$. Select a plane $\pi$ through $L_1$ not in a $3$-space with $L_2$, and choose a point $F_3'$ in $\pi$ not on $L_1$ and not collinear to $L_2$. As above, we may substitute $F_3$ by $F_3'$ in our sequence of perspectivities. Now none of $\<F,F_3'\>$ or $\<F_3,F_3'\>$ intersect $L_2$, and so we may assume now that $L_1$ does not intersect $L_2$. Then there is a unique point $p_2$ on $L_2$ collinear to all points of $L_1$, defining the plane $\pi:=\<p_2,L_1\>$. We select a plane $\alpha$ through $L_1$ not in a $3$-space with $\pi$. It follows that $\proj_\alpha L_2\in L_1$. Thus we may consider any point $F_3'$ of $\alpha\setminus L_1$ and we find that $L_2$ is opposite both $\<F,F_3'\>$ and $\<F_3,F_3'\>$. This completes the argument.   
\end{remark}

\begin{lem}\label{upanddownF4}
 Let $\Delta$ be a metasymplectic space.  Suppose that for each pair of points  $p$, $p'$, the graph $\Gamma(S)$, where $S$ is the set of points opposite both $p$ and $p'$, is connected. Suppose also that there exists a point opposite any given triple of points. Let $p$ be a given point. Denote by $\Pi_4(p)$ the set of all self-projectivities $p\per p_2\per p_3\per p_4\per p$ of $p$ of length $4$ with $p\perp p_3$, $p_2\perp p_4$ and $pp_3$ opposite $p_2p_4$. Suppose that $\Pi_4(F)$ is geometric. Then $\Pi^+(F)=\<\Pi_4(F)\>$. 
 \end{lem}

\begin{proof}
The statement, without the condition that $pp_3$ is opposite $p_2p_4$, is \cite[Lemma~8.17]{Bus-Sch-Mal:24} specialised to buildings of type $\mathsf{F_4}$ and vertices of type $1$ or $4$. 

Now let  $p\per p_2\per p_3\per p_4\per p$ be a self-projectivity of $p$ with $p\perp p_3$ and $p_2\perp p_4$. Let $L_3$ be the projection of $pp_3$ onto the residue of $p_2$. Set $L_4:=p_2p_4$. In $\Res_\Delta(p_2)$, the elements $L_3$ and $L_4$ have distance $0,1,2$ or $3$ from each other. If they have distance $3$, then \cite[Proposition~3.29]{Tits:74} implies that $pp_3$ and $p_2p_4$ are opposite. Hence, we may assume that their distance $d$ is $0,1$ or $2$. Then there exists a line $L_2$ through $p_2$ at distance   $d+1$ from $L_3$ in $\Res_\Delta(p)$ and at distance $1$ from $L_4$ in that residue, that is, coplanar with $L_4$ in $\Delta$.   Let $\pi_0$ be the plane spanned by $L_2$ and $L_4$. Let $L_1$ be the projection of $L_2$ onto (the residue of) $p$, and let $q_2$ be a point on $L_2$ opposite $p$. The projection of the line $p_4q_2$ from $p_4$ onto $p$ is contained in the projection of $\pi_0$ onto $p$, is different from the projection of $p_4p_2$ onto $p$ and hence, has distance $d+1$ to $L_3$. If we call the \emph{width}  of   $p\per p_2\per p_3\per p_4\per p$ the distance $d$ in $\Res_\Delta(p)$ between the lines $pp_3$ and the projection of $p_2p_4$ onto $p$, then we can write
\[p\per p_2\per p_3\per p_4\per p = (p\per p_2\per p_3\per q_2\per p)\cdot (p\per q_2\per p_3\per p_4\per p),\]where the width of both self projectivities of $p$ on the right hand side is $d+1$ (look at the inverse of the second to see this). An induction argument on $d$ proves the assertion. 
\end{proof}

\begin{remark}
One might wonder, why one would bother to reduce the computation of the special projectivity groups to the computation of self-projectivities of length $4$, when we can reduce it to the computation of those of length $3$ with \cref{triangles}. The reason is, firstly, that in general, a generic subspace opposite two given subspaces has an algebraically complicated form with many parameters (especially when the dimension of $F$ is large). In the length $4$ case (see \cref{openneer}), there are essentially only two single parameters: one to fix $F_3$, and one to fix $F_4$, whereas $F,F_2,F\cap F_3$ and $F_2\cap F_4$ can be chosen freely. Secondly, the conditions are different, and sometimes those of \cref{upanddown} are easier to meet than those of \cref{triangles}. 
\end{remark}


\section{Preparations and auxiliary results}\label{prep}

\subsection{Conditions for reduction for polar spaces}
In this subsection, we check the conditions of \cref{triangles} and \cref{upanddown} in the cases where we shall apply them.  We start with the polar spaces.

Regarding \cref{triangles}, we can be brief: the following is proved in \cite{Bus-Mal:24}.

\begin{prop}\label{lines}
Let $\Delta$ be a polar space of rank $r\geq 2$. Suppose first that each line contains at least four points and let $0\leq j<r-1$. Then, if $j>0$, for each quadruple of singular subspaces of dimension $j$ containing at least one opposite pair, there exists a singular subspace of dimension $j$ opposite all the members of the given quadruple. If $j=0$, the same conclusion holds for each quadruple of points with the property that the pairwise intersections of the perps are not all the same. Secondly, suppose that each submaximal subspace is contained in at least four maximal singular subspaces, and $r\geq 3$, then each quadruple of maximal singular subspaces, containing at least one opposite pair, admits a common opposite singular subspace. 
\end{prop}

Now we consider the main condition in \cref{upanddown}. First we treat the rank $2$ case, although we only consider the projectivity groups for rank at least $3$ --- we need the rank $2$ case for induction purposes. Strictly speaking, we could assume that the rank $2$ polar spaces we deal with are Moufang, but we prove the result for all generalised quadrangles (with at least four points per line and four lines through a point, respectively). 

\begin{prop}\label{rank2case}
Let $\Delta$ be a thick polar space of rank $2$. Let $p_{1}$ and $p_{2}$ be two points. Let $S$ be the set of points opposite both $p_{1}$ and $p_{2}$. Then $\Gamma(S)$ is connected, if each line contains at least four points. 
\end{prop}

\begin{proof}Assume each line contains at least four points. 
Let  $q_{1}$ and $q_{2}$ be two points in $\Delta \setminus (p_{1}^{\perp} \cup p_{2}^{\perp})$.  If $q_{1}$ and $q_{2}$ are collinear in $\Delta$, then they are also collinear in $\Delta \setminus (p_{1}^{\perp} \cup p_{2}^{\perp})$. Suppose $q_{1}$ and $q_{2}$ are opposite in $\Delta$. Let $L$ be a line through $q_{1}$ in $\Delta$. 

If $\proj_{L} (q_{2})$ is not equal to $\proj_{L} (p_{1})$ or $\proj_{L} (p_{2})$, then $\proj_{L} (q_{2})$ is a vertex of $\Delta \setminus (p_{1}^{\perp} \cup p_{2}^{\perp})$ and the lines $q_{2}\proj_{L} (q_{2})$ and $q_{1}\proj_{L} (q_{2})$ give rise to edges in $\Gamma(\Delta \setminus (p_{1}^{\perp} \cup p_{2}^{\perp}))$ that form a path between the vertices corresponding to $q_{1}$ and $q_{2}$ via the vertex corresponding to $\proj_{L} (q_{2})$.

Only if for every line $L$ through $q_{1}$ it is the case that $\proj_{L} (q_{2})$ is equal to either $\proj_{L} (p_{1})$ or $\proj_{L} (p_{2})$, we can not immediately find such a path. So assume this is the case. Then every point in $q_{1}^{\perp} \cap q_{2}^{\perp}$ is either collinear to $p_{1}$ or to $p_{2}$. 

Suppose $p^\perp_1\cap p^\perp_2\cap q_1^\perp\cap q_2^\perp\neq\varnothing$ and let $x$ be collinear to all of $p_1,p_2,q_1,q_2$. Let $L$ be a line through $q_1$ not containing $x$. Without loss of generality, we may assume that the unique point $y$ on $L$ collinear to $q_2$ is collinear to $p_1$. Let $s_1$ be a point on $L$ not collinear to $p_2$, and distinct from $q_1$ and $y$. Let $s_2$ be the point on $q_2x$ collinear to $s_1$. None of $s_1$ or $s_2$ is collinear to either $p_1$ or $p_2$ and we have the path $q_1\sim s_1\sim s_2\sim q_2$. 

Hence, we may suppose that each point of $q_1^\perp\cap q_2^\perp$ is collinear to either $p_1$, or $p_2$, but never to both. Without loss of generality, we may assume that $p_1$ is collinear to at least $2$ points $o_1,o_2$ of $q_1^\perp\cap q_2^\perp$. If there are at least $5$ points on a line, then there is a point $r_1$ on $q_1o_1$ distinct from all of $q_1,o_1,\proj_{q_1o_1}p_2$ and $\proj_{q_1o_1}\proj_{q_2o_2}p_2$.  Then we have the path $q_1\sim r_1\sim\proj_{q_2o_2}r_1\sim q_2$ . Consequently, we may assume that there are exactly $4$ points per line, and hence, in view of the main result of \cite{Bro:91}, there are a finite number of lines through each point, say $n$. By assumption, $p_1$ is collinear to $m\geq 2$ points $o_1,o_2,\ldots,o_m$ of $q_1^\perp\cap q_2^\perp$. Hence, there are $n-m$ lines through $p_2$ containing a point of $q_1^\perp\cap q_2^\perp$. The other $m$ lines through $p_2$ each have to meet at least one of the lines $q_10_i$, $q_2o_i$, $1\leq i\leq m$, and each such line must meet one of those $m$ lines through $p_2$. It follows that there is a line  $L$ through $p_2$ meeting, without loss of generality, the line $q_1o_1$ and $q_2o_2$.  Then the unique points $r_1$ and $r_2$ on $q_1o_1$ and $q_2o_2$, respectively, distinct from $q_1,q_2,o_1,o_2, q_1o_1\cap L$ and $q_2o_2\cap L$, are collinear and opposite both $p_1$ and $p_2$. Hence, we have the path $q_1\sim r_1\sim r_2\sim q_2$. 
\end{proof}

Dually, we have:

\begin{coro}\label{rank2casedual}
Let $\Delta$ be a thick polar space of rank $2$. Let $L_{1}$ and $L_{2}$ be two lines. Let $S$ be the set of lines opposite both $L_{1}$ and $L_{2}$. Then $\Gamma(S)$ is connected, if each point is on at least four lines.
\end{coro}

Now we treat the higher rank cases.

\begin{prop}\label{higherrank}
Let $\Delta$ be a thick polar space of rank $n \geq 3$. Let $U_{1}$ and $U_{2}$ be two singular subspaces of dimension $s$. Let $S_{U_1,U_2}$ be the set of all singular subspaces of dimension $s$ opposite both $U_{1}$ and $U_{2}$. Then $\Gamma(S_{U_1,U_2})$ is connected for $s\leq n-2$, if each line has at least four points (this additional condition is not needed for $s=0$) and for $s=n-1$, if either each submaximal singular subspace is contained in at least four maximal singular subspaces or each line contains at least four points.  
\end{prop}

\begin{proof}The proof goes with induction on the rank $n$, the base case being \cref{rank2case} and \cref{rank2casedual}. However, for $s=0$, we provide an independent proof, neglecting the additional conditions on the sizes of the lines. 

\textbf{Case 1: $s=0$.} Let $p_{1}$ and $p_{2}$ be two points. Let $\Delta \setminus (p_{1}^{\perp} \cup p_{2}^{\perp})$ be the point-line geometry containing all points of $\Delta$ that are not collinear to either $p_{1}$ or $p_{2}$ and all lines of $\Delta$ between these points. These lines are exactly the lines $L$, which are not in a plane with either $p_{1}$ or $p_{2}$ and for which $L \setminus (\proj_{L} (p_{1}) \cup \proj_{L} (p_{2}))$ contains more than one point. 
We have to show that the point graph of $\Delta \setminus (p_{1}^{\perp} \cup p_{2}^{\perp})$ is connected.

Let $q_{1}$ and $q_{2}$ be two arbitrary points of $\Delta \setminus (p_{1}^{\perp} \cup p_{2}^{\perp})$. If $q_{1}$ and $q_{2}$ are collinear in $\Delta$, then they are also collinear in $\Delta \setminus (p_{1}^{\perp} \cup p_{2}^{\perp})$. Suppose $q_{1}$ and $q_{2}$ are opposite in $\Delta$. As in the proof of \cref{rank2case}, we may assume that every point in $q_{1}^{\perp} \cap q_{2}^{\perp}$ is either collinear to $p_{1}$ or to $p_{2}$. Since $q_{1}$ and $q_{2}$ are opposite, $q_{1}^{\perp} \cap q_{2}^{\perp} =: \Lambda$ defines a thick polar space of rank $n-1 \geq 2$. Both $p_{1}^{\perp} \cap \Lambda$ and $p_{2}^{\perp} \cap \Lambda$ define geometric hyperplanes that we will denote by $H_{1}$ and $H_{2}$ and we have $\Lambda = H_{1} \cup H_{2}$.

Suppose $H_{1}$ and $H_{2}$ are proper geometric hyperplanes. Then we obtain a contradiction, because a polar space can never be the union of two proper geometric hyperplanes (see Exercise 2.5 in \cite{Mal:24}): Let $x$ be a point of $H_{1} \setminus H_{2}$, such that $x^{\perp} \neq H_{1}$. Then $H_{1}$ induces a proper geometric hyperplane in $\Res_{\Lambda} (x)$ and $H_{2}$ has to contain the complement. Let $M$ be a line through $x$ that is not contained in $H_{1}$. Then $M \setminus \{x\}$ has to be contained in $H_{2}$ and with that, $x$ has to be contained in $H_{2}$, which is a contradiction.

So we may assume, without loss of generality, $p_{1}^\perp \cap \Lambda = \Lambda$. Since $H_{2}$ is a hyperplane, it contains two opposite points $o_{1}$ and $o_{2}$ and since $H_{2} \subseteq H_{1} = \Lambda$, the points $o_{1}$ and $o_{2}$ are both collinear to $q_{1}$, $q_{2}$, $p_{1}$ and $p_{2}$ and the lines $q_{1}o_{1}$ and $q_{2}o_{2}$ are opposite. Let $x_{1}$ be a point on $q_{1}o_{1}\setminus\{q_1\}$ that is not in $\Lambda$. The projection of $x_{1}$ onto $q_{2}o_{2}$ is a point $x_{2}$ that is not in $\Lambda$ and that does not coincide with $q_{2}$. Furthermore, $x_{1}$ and $x_{2}$ are not collinear to either $p_{1}$ or $p_{2}$ and thus contained in $\Delta \setminus (p_{1}^{\perp} \cup p_{2}^{\perp})$. The lines $q_{1}x_{1}$, $x_{1}x_{2}$ and $q_{2}x_{2}$ give rise to edges in $\Gamma(S_{p_1,p_2})$ that form a path $q_1\sim x_1\sim x_2\sim q_2$ between the vertices corresponding to $q_{1}$ and $q_{2}$.

This concludes the proof for the case $s=1$. For cases 2 and 3 below we assume that each line has at least four points. 
But before moving to these other cases, we prove a common claim under our general assumptions.

\textbf{Claim:} \emph{If $1\leq s\leq n-1$, then different components of $\Gamma(S_{U_1,U_2})$ only contain vertices corresponding to disjoint subspaces.}
Indeed, let $V_{1}$ and $V_{2}$ belong to $\Gamma(S_{U_1,U_2})$, with $ x \in V_{1} \cap V_{2}$. The projection of $x$ onto $U_{i}$ is a hyperplane $H_{i}$ of $U_i$, $i \in \{1,2\}$. With \cite[Proposition~3.29]{Tits:74}, it follows that every subspace $X$ through $x$, that is opposite $\langle x, H_{i} \rangle$ in $\Res(x)$, is opposite $U_{i}$ in $\Delta$. Since $\Res(x)$ is a polar space of rank $n-1$, it follows from the induction hypothesis that we can find a path between the vertices corresponding to $V_{1}$ and $V_{2}$ in $\Gamma(S_{U_1,U_2})$. The claim is proved.

\textbf{Case 2: $s = 1$.} Remember the rank of $\Delta$ is at least 3, so $\Delta$ contains planes.  Let $L_{1}$ and $L_{2}$ be two lines in $\Delta$ and let $S_{L_{1}, L_{2}}$ be the set of all lines opposite both $L_{1}$ and $L_{2}$. Let $J_{1}$ and $J_{2}$ be two lines of $S_{L_{1}, L_{2}}$. 

(I) If $J_{1}$ and $J_{2}$  intersect in a point, then by the above claim, we can find a path in $\Gamma(S_{L_1,L_2})$ connecting $J_1$ and $J_2$. 

(II) Now suppose $J_{1}$ and $J_{2}$ do not intersect in a point, but there is some point $x\in J_1$ collinear to all points of $J_2$. Set $\pi:=\<x,J_2\>$. Since $J_j$ is opposite $L_i$, $i,j=1,2$, the projection of $L_i$ onto $\pi$ is a point $u_i\in\pi\setminus(J_2\cup\{x\})$. Select a line $J_1'$ through $x$ avoiding both $u_1$ and $u_2$, then $J_1'$ is opposite both $L_1$ and $L_2$ and intersects both $J_1$ and $J_2$. By (I) we can connect $J_1$ with $J_1'$ and $J_1'$ with $J_2$.

(III) Finally, suppose that $J_1$ and $J_2$ are opposite. Let $\alpha_{1}$ and $\beta_{1}$ be two arbitrary planes containing  $J_{1}$ and let $\alpha_2$ and $\beta_2$ be the planes through $J_{2}$ intersecting $\alpha_1$ and $\beta_1$, respectively. Set $a:=\alpha_1\cap\alpha_2$ and $b:=\beta_1\cap\beta_2$. If either $a$ or $b$ is not collinear to either $L_1$ or $L_2$, say $a$, then we find a line $J_1'$ through $a$ in $\alpha_1$ opposite both $L_1$ and $L_2$, and we are reduced to Cases (I) and (II). In the other case, say $a\perp L_1$ and $b\perp L_i$, $i\in\{1,2\}$, there is a unique point $a'$ in $\alpha_1$ collinear to $L_2$, and there is a unique point $b'$ in $\beta_2$ collinear to $L_j$, $\{j\}=\{1,2\}\setminus\{i\}$.  Since lines have at least four points, we can now find a line $L_2'$ in $\beta_2$ distinct from both $J_2$ and $\proj_{\beta_2}a'$, and  avoiding both $b$ and $b'$. Then we can find a line $J_1'$ in $\alpha_1$ through $\proj_{\alpha_1}J_2'$ avoiding $a$ and $a'$. Hence, both $J_1'$ and $J_2'$ are opposite both $L_1$ and $L_2$ and by (I), $L_1$ and $L_1'$ are connected in $\Gamma(S_{L_1,L_2})$ and $L_2$ and $L_2'$ are connected in $\Gamma(S_{L_1,L_2})$; by (II) also $J_1'$ and $J_2'$ are connected in $\Gamma(S_{L_1,L_2})$. Hence, $L_1$ and $L_2$ are connected in $\Gamma(S_{L_1,L_2})$.

\textbf{Case 3: $2\leq s \leq n-1$.} 
Let $U_{1}$ and $U_{2}$ be two singular subspaces of dimension $s$ and $V_{1}$ and $V_{2}$ both opposite both $U_{1}$ and $U_{2}$. On top of the induction on the rank of $\Delta$, we also perform an induction on the dimension $s$ of $U_1$ and $U_2$. The base cases here are Case 1 and Case 2. 

(I) If $V_1$ and $V_2$ are not disjoint, then by the above claim, we can find a path in  $\Gamma(S_{U_1,U_2})$ connecting $V_1$ and $V_2$.

(II) Secondly, assume $V_1\cap V_2=\varnothing$.  Let $A_{1}$ be some $(s-1)$-dimensional subspace in $V_{1}$. Then $A_{1}$ is opposite some $(s-1)$-dimensional subspaces $B_{1}$ of $U_{1}$ and $B_{2}$ of $U_{2}$. Let $A_{2}$ be some $(s-1)$-dimensional subspace in $V_{2}$ that is opposite both $B_{1}$ and $B_{2}$. Then we can find a path between $A_{1}$ and $A_{2}$ in $\Gamma(S_{B_1,B_2})$. Let $A_{1}=X_{1},~ X_{2},~ \dots,~ X_{r}=A_{2}$ be the $(s-1)$-dimensional subspaces corresponding to all vertices of that path, such that $X_{j}$ is adjacent to $X_{j+1}$ for $j \in \{1,2,\ldots,r-1\}$. Then $X_{j}$ is an $(s-1)$-dimensional subspace that intersects both $X_{j-1}$ and $X_{j+1}$ in $(s-2)$-dimensional subspaces and that is opposite both $B_{1}$ and $B_{2}$. 

The projection of $U_{i}$ onto $X_{j}$ (for $j \in \{2,3,\ldots,r-1\}$ and $i \in \{1,2\}$) is a singular subspace $W_{ij}$ of dimension $s$ containing $X_j$.  Let $W_{j}$ be an $s$-dimensional subspace containing $X_{j}$ opposite in $\Res(X_j)$ both $W_{1j}$ and $W_{2j}$.  Then, again using \cite[Proposition~3.29]{Tits:74},  $W_{j}$ is opposite both $U_{1}$ and $U_{2}$. Set $W_1:=V_1$ and $W_r=V_2$. Since $X_{j}$ intersects $X_{j+1}$, $j\in\{1,2,\ldots,r-1\}$, and $X_{j} \subseteq W_{j}$, we find that $W_{j}$ intersects $W_{j+1}$, for all $j\in\{1,2\ldots,r-1\}$. Now we can use (I) again the find paths between $W_j$ and $W_{j+1}$ in $\Gamma(S_{U_1,U_2})$, for $j\in\{1,2\ldots,r-1\}$, which, taken together, form one big path between $W_1=V_1$ and $W_r=V_2$ in $\Gamma(S_{U_1,U_2})$. This concludes the proof of Case~3.

For the final case, we assume that each submaximal singular subspace is contained in at least four maximal singular subspaces.

\textbf{Case~4: $s=n-1$.} Let $M_1$ and $M_2$ be two maximal singular subspaces and let $N_1,N_2$ be two maximal singular subspaces opposite both $M_1$ and $M_2$. As in the previous case, by the Claim above, we may assume that $N_1\cap N_2=\varnothing$, hence $N_1$ and $N_2$ are themselves opposite.  We again proceed by induction on $n$, the case $n=2$ being \cref{rank2casedual}. So we assume $n>2$. Select $x\in N_1$, and let $N$ be the unique maximal singular subspace containing $x$ and intersecting $N_2$ in a hyperplane of $N_2$. Then $M_i\cap N$ is at most a point, $i=1,2$, and we find a hyperplane $H$ of $N$ containing $x$ and disjoint from both $M_1$ and $M_2$. Through $H$, there are at most two maximal singular subspaces not opposite either $M_1$ or $M_2$, and so there exists a maximal singular subspace $N'$ through $H$ opposite both $M_1$ and $M_2$. Applying the claim above to $N_1$ and $N'$, and to $N'$ and $N_2$ (which meet in $H\cap N_2\neq\varnothing$), Case~4 is proved. 
\end{proof}


\subsection{Conditions for reduction for metasymplectic polar spaces}
We now check the conditions of \cref{upanddownF4}.  We will only need the result for points. 

First, we present two lemmas for dual polar spaces. We only need them in rank 3, but the proofs are the same for general rank.

\begin{lem}\label{DPS1}
A dual polar space $\Gamma$ is not the union of two proper geometric hyperplanes.
\end{lem} 

\begin{proof}
Let $H$ and $H'$ be two geometric hyperplanes of $\Gamma$ with $H\cup H'=X$, where $X$ is the point set of $\Gamma$. Recall that a \emph{deep point} of $H$ is a point of $H$ such that $x^\perp\subseteq H$. Each point outside $H$ belongs to $H'$, as well as all points of $H$ that are not deep. We now claim that a deep point $x$ of $H$ belongs to $H'$ as soon as $x^\perp$ contains a point that is not deep. Indeed, let $y\perp x$ not be deep and let $L$ be a line through $y$ not contained in $H$. Let $\xi$ be a symp containing $x$ and $L$. Since $x$ is deep, $H\cap\xi= x^\perp\cap\xi$. Hence, no point of $L\setminus\{x\}$ is deep. We conclude that $L\subseteq H'$, proving the claim. Hence, if $H'$ is proper, then the set of deep points of $H$ is closed under taking perps. Since $\Gamma$ is connected, only $X$ is a non-empty subset of points closed under taking perps. The lemma now follows. 
\end{proof}

\begin{lem}\label{DPS2}
Let $H$ and $H'$ be two proper geometric hyperplanes of a dual polar space $\Gamma$. Then $H$ contains a point $x$ opposite some point $x'\in H'$. 
\end{lem}

\begin{proof}
Let $y\in H$ be arbitrary. Suppose no point of $H'$ is opposite $y$. Then $H'\subseteq y^{\not\equiv}$. Since no point at distance $2$ can be removed from $y^{\not\equiv}$ without losing the property of being a hyperplane --- and then also no point of $y^\perp$ can be deleted --- we see that this implies $H'=y^{\not\equiv}$,  Now let $x\in H$ be any other point of $H$ (hence $x\neq y$). \cref{Opp-notopp} yields a point $x'$ not opposite $y$ but opposite $x$. Then $x'\in H'$ and the assertion is proved.  
\end{proof}
We are now ready to deal with the condition in \cref{upanddownF4} about the connectivity of $\Gamma(S)$. 
\begin{lem}
Let $\Delta$ be a metasymplectic space in which either every line has at least four points or each plane is contained in at least $4$ symps. For each pair of points $p_0$, $p_1$ the graph $\Gamma(S)$, where $S$ is the set of points opposite both $p_0$ and $p_1$ is connected.
\end{lem}

\begin{proof}
\begin{compactenum}[(I)]
\item  Let $q_0$ and $q_1$ be two points opposite both $p_0$ and $p_1$. If $q_0$ and $q_1$ are collinear in $\Delta$, then they are connected in $\Gamma(S)$. Suppose $q_0$ and $q_1$ are symplectic in $\Delta$ and denote by $\xi(q_0, q_1)$ the symp containing both of them. Then no $p_i$, for $i\in\{0,1\}$, can be contained in $\xi(q_0, q_1)$ and can also not be close to $\xi(q_0, q_1)$, since otherwise both $q_0$ and $q_1$ would not be opposite $p_i$.  That means $p_i^{\pperp} \cap \xi(q_0, q_1)$ has to be a unique point $r_i$. We can find a path in the graph corresponding to $\xi(q_0, q_1) \setminus (r_0^\perp \cup r_1^\perp)$ between $q_0$ and $q_1$.

\item Suppose $q_0$ and $q_1$ are special in $\Delta$. Let $q$ be the unique point collinear to both $q_0$ and $q_1$.  Then we can not find a path in $\Gamma(S)$ between $q_0$ and $q_1$ immediately, if $q$ is collinear to a point of $p_0^{\perp}$ or $p_1^{\perp}$. Let $p$ be a point in $p_0^{\perp} \cup p_1^{\perp}$ collinear to $q$. We may assume $p \perp p_0$. 

We have paths $p_0 \perp p \perp q \perp q_0$ and $p_0 \perp p \perp q \perp q_1$ and know that $p_0$ is opposite both $q_0$ and $q_1$. Therefore, we know that the pairs $(p, q_0)$, $(p, q_1)$ and $(q, p_0)$ are special and if $p_1$ has distance $2$ to $q$, then the pair $(q,p_1)$ is also special. We first consider this case and define $p'$ as the unique point of $p_1^{\perp} \cap q^{\perp}$. 

In the residue of $q$ we can see the lines $pq$, $q_0q$, $q_1q$ and $p'q$ as points. Since $\Res(q)$ is a dual polar space, we can apply \cref{higherrank} and find a path between the points $q_0q$ and $q_1q$ in $\Res(q)$ all members of which are locally opposite $pq$ and $p'q$. In $\Delta$ we can see that path as a chain of planes such that consecutive planes intersect in lines $L$ locally opposite both $pq$ and $p'q$. On each such line $L$, we choose a point distinct from both $q$ and the projection of $p_1$ onto $L$. The thus obtained path in $\Delta$ connects $q_0$ with $q_1$ within $p_0^{\equiv}\cap p_1^\equiv$.

\item Suppose $q_0$ and $q_1$ are opposite in $\Delta$. Set $Z_0=q_0^\perp\cap q_1^{\Join}$. If $q\in Z_0$ is opposite both $p_0$ and $p_1$, then  we can find paths from $q_0$ to $q$ and, by Case (II), from $q$ to $q_1$ in $\Gamma(S)$ and the concatenation is a path between $q_0$ and $q_1$. Suppose such $q\in Z_o$ does not exist, in other words, $(p_0^{\not\equiv}\cup p_1^{\not\equiv})\cap Z_0=Z_0$. Now both $p_0^{\not\equiv}$ and  $p_1^{\not\equiv}$ are geometric hyperplanes of $Z_0$. By \cref{DPS1}, we may assume $Z_0\subseteq p_0^{\not\equiv}$. Set $Y_0=Z_0\cap p_1^{\not\equiv}$. Likewise, defining $Z_1=q_1^\perp\cap q_0^{\Join}$, one of $p_0^{\not\equiv}$ and $p_1^{\not\equiv}$ contains $Z_1$. The intersection of $Z_1$ with the other is denoted as $Y_1$. Let $Y_1'$ be the projection of $Y_1$ onto $Z_0$. Then \cref{DPS2} yields a point $x_0\in Y_0$ and a point $x_1'\in Y_1'$ with $q_0x_0$ locally opposite $q_0x_1'$. Projecting $x_1'$ back onto $Z_1$, we obtain a point $x_1\in Y_0$.  Now pick $q_1'\in q_1x_1\setminus\{x_1,q_1\}$ and note that, by this choice, $q'_1\in \Gamma(S)$. Then the projection $q_0'$ of $q_1'$ onto $q_0x_0$ is a point distinct from $x_0$ and hence, belongs to  $\Gamma(S)$. Since $x_0$ is special to $x_1$, Case (II) implies that $x_0$ and $x_1$ are connected in $\Gamma(S)$. Hence, $q_0\perp x_0$ is connected to $q_1\perp x_1$ in $\Gamma(S)$ and the lemma is proved. \qedhere 
\end{compactenum}
\end{proof}

Now we handle the condition in \cref{upanddownF4} about the existence of a point opposite three given points.  It follows from \cref{Tits3.30} that this condition is automatically satisfied whenever the building has no residues isomorphic to the unique projective plane of order $2$ (with $3$ points per line). It follows from the main result in \cite{Bus-Mal:25} that only triples of points that form a geometric line have no common opposite. But such a triple is determined by any pair of its elements.  
Now, it also follows from the proof of \cref{upanddownF4} in \cite[Section 8]{Bus-Sch-Mal:24} that the conclusion of \cref{upanddownF4}, possibly without the claim of $pp_3$ being opposite $p_2p_3$, holds for all projectivities that can be written as a product 
\[p_1\per p_2\per p_3\per \ldots\per p_{2\ell-1}\per p_{2\ell}\per p_1,\]
such that we can find a common opposite for the triples $\{p_{1},p_{2k-1},p_{2k+1}\}$, for all $k=2,\ldots,\ell-1$. Hence, if we were able to replace every subsequence $p_{2k-1}\per p_{2k}\per p_{2k+1}$, $k=2,\ldots,\ell-1$, for which $\{p_1,p_{2k-1},p_{2k+1}\}$ is a geometric line, with a sequence $p_{2k-1}\per p_{2k}\per q_{2k}\per p_{2k}\per p_{2k+1}$, where $q_{2k}\equiv p_{2k}$ and none of $\{p_1,p_{2k-1},q_{2k}\}$ or $\{p_1,q_{2k},p_{2k+1}\}$ is a geometric line, then the conclusion of \cref{upanddown} would still hold. But this can be achieved since, by \cref{Opp-notopp}, there exists a point $q_{2k}$ opposite $p_{2k}$ and distinct from $p_1$.  
Since the proof of \cref{upanddown} in the present paper does not use the condition of the existence of a point opposite three given points, we conclude

\begin{lem}\label{UpanddownF4}
 Let $\Delta$ be a metasymplectic space.  Let $p$ be a given point. Denote by $\Pi_4(p)$ the set of all self-projectivities $p\per p_2\per p_3\per p_4\per p$ of $p$ of length $4$ with $p\perp p_3$, $p_2\perp p_4$ and $pp_3$ opposite $p_2p_4$. Suppose that $\Pi_4(F)$ is geometric. Then $\Pi^+(F)=\<\Pi_4(F)\>$. 
 \end{lem}

\subsection{Collineations pointwise fixing a (sub)hyperpane}We will see that self-projectivities of length~$3$ in $\Pi(x)$, for $x$ a point in a polar space $\Delta$, pointwise fix a hyperplane of $\Res_\Delta(x)$. Therefore, we prove some results about collineations of polar spaces pointwise fixing a hyperplane. We call a collineation of a polar space $\Delta$ that pointwise fixes a geometric hyperplane of $\Delta$ a \emph{reflection} (In \cite{Die:73} it is called a \emph{symmetry}).
\def\cC{\mathcal{C}}
\begin{lem}\label{exrefl}
Let $\Delta$ be a separable quadric of Witt index $r\geq 1$ in $\PG(V)$, for some vector space $V$ over the field $\K$. We extend the perp-relation to all subspaces of $\PG(V)$ according to the non-degenerate polarity defined by the quadratic form defining $\Delta$.  Let $\theta$ be a central collineation of $\PG(V)$ with centre $c$ and axis the hyperplane $H:=c^\perp$. If $\theta$ maps some point $x$ of $\Delta$ to a distinct point $x^\theta$ of $\Delta$, then $\theta$ preserves $\Delta$ and consequently defines a collineation of $\Delta$.
\end{lem}

\begin{proof}
If $\ch\K\neq 2$, then this is a symmetry as defined by Dieudonn\'e in \cite[Section~10]{Die:73} and the result follows from that reference.

Now suppose $\ch\K=2$. This situation is less standard, especially when $\K$ is not perfect. We have to show that, if $y$ is an arbitrary point of $\Delta$, then $\theta$ maps $y$ to a point of $\Delta$. Considering the plane $\<x,c,y\>$, this boils down to showing, given a conic $\cC$ in a plane $\PG(2,\K)$ and a line $L$ through the nucleus $n$ of $\cC$, the existence of a non-trivial central elation with a given arbitrary centre $c$ on $L$, with $c\neq n$. We can take the conic with equation $Y^2=XZ$, and $L$ has equation $X=kZ$, for some $k\in\K$. The centre $c$ has coordinates $(\ell,1,k\ell)$. Then the non-trivial elation pointwise fixing $L$ with centre $c$ and preserving $\cC$ has matrix
\[\begin{pmatrix}
(1+ab^2)^{-1} & 0 & b^2(1+ab^2)^{-1} \\
ab(1+ab^2)^{-1} & 1 & b(1+ab^2)^{-1}\\
a^2b^2(1+ab^2)^{-1} & 0 & (1+ab^2)^{-1}
\end{pmatrix},\]
as one can check with elementary calculations. This concludes the proof.
\end{proof}

\begin{lem}\label{rank1}
Let $\Delta$ be a separable quadric of Witt index $r\geq 1$ in $\PG(V)$, for some vector space $V$ over the field $\K$. We extend the perp-relation to all subspaces of $\PG(V)$ according to the non-degenerate polarity defined by the quadratic form defining $\Delta$.  Let $\theta$ be a collineation pointwise fixing a subhyperplane $G$ of $\PG(V)$. Set $L=G^\perp$. We assume that $L$ is not a tangent, that is, that either $L\subseteq G^\perp$ with $L\cap\Delta=\varnothing$, or $L\cap G=\varnothing$, and the former only occurs in characteristic $2$.

Then $\theta$ is the product of two reflections.
\end{lem}

\begin{proof}
Suppose first $L\cap G=\varnothing$. Then $L$ is either disjoint from $\Delta$, or intersects it in exactly two points.  Let $x$ be a point of $\Delta$ not contained in $G\cup L$. We also assume that, if $L$ contains two points $x_1,x_2$ of $\Delta$, then $x$ is not contained in $x_1^\perp\cup x_2^\perp$.  As $\theta$ fixes every subspace of $G$, it stabilises every subspace through $L$ and hence, it stabilises the plane $\langle x, L \rangle$. So we can define $\{a\}=xx^\theta\cap L$. If $x\perp x^\theta$, then $a$ belongs to $\Delta$ and so is one of $x_1,x_2$, contradicting our choice of $x$. Hence, $a$ does not belong to $\Delta$. Also, $a \notin x^{\perp}$ and  $a \notin (x^{\theta})^{\perp}$, since $a$ is on a secant through $x$ and $x^\theta$. Since $L^{\perp}=G$, $a \in L$ implies $G \subseteq a^{\perp}$.

Let $\theta_{1}$ be the unique central collineation with centre $a$ that maps $x$ to $x^{\theta}$ and fixes $a^{\rho}$ pointwise. By \cref{exrefl}, $\theta_1$ defines a collineation of $\Delta$. Define $\theta_{2}$ to be $\theta \theta_{1}^{-1}$.

Let $c$ be the point $\langle x, G \rangle \cap L$. Then $c$ maps to $\langle x^{\theta}, G \rangle \cap L = \langle x^{\theta_{1}}, H \rangle \cap L$ under $\theta_1$. Consequently $c^\theta=c^{\theta_1}$ and so $c$ is fixed under $\theta\theta_1^{-1}=\theta_2$.

With that, $\theta\theta_{1}^{-1}$ fixes $G$ pointwise and also fixes $x$ and $c$. Hence, it pointwise fixes $\langle x,H\rangle$, and hence, $\theta=(\theta\theta_1^{-1})\theta_1$ is the product of two central collineations.

Secondly, suppose $L \subseteq H$. Now $L$ does not contain any point of $\Delta$ and we may pick $x$ in $\Delta$, but not in $G$, arbitrarily. As before, we define $a$ as the intersection of $xx^\theta\cap L$. We also know that $x$ is not collinear to $x^\theta$.

Define $\theta_{1}$ as the elation that fixes $a^{\perp}$ pointwise, has centre $a$ and maps $x$ to $x^{\theta}$. Then, again by \cref{exrefl}, $\theta_1$ induces a collineation in $\Delta$.

The line $xx^{\theta}$ already has two points in the quadric, so there are no more points and therefore $x^{\theta}$ maps back to $x$. That means $\theta$ is an involution (which we can of course also deduce from the fact that $\ch\K=2$).

The map $\theta \theta_{1}^{-1}$ fixes $H$ pointwise and fixes both $x$ and $x^\theta$. Hence, it pointwise fixes $\langle H,x\rangle$. Consequently, $\theta=(\theta\theta_1^{-1})\theta_1$ is the product of two central collineations.
\end{proof}

\section{Projectivity groups for upper residues in polar spaces}\label{upperPS}
\subsection{Projectivity groups of points}
\subsubsection{The generic case}
\begin{prop}\label{pointsgeneral}
Let $\Delta$ be a thick embeddable polar space of rank $n \geq 3$. Let $p_{1}$ and $p_{2}$ be two opposite points of $\Delta$, let $\Gamma$ be the polar space $p_{1}^{\perp} \cap p_{2}^{\perp}$ of rank $n-1$ inside $\Delta$ and let $H$ be a hyperplane of $\Gamma$, obtained as intersection of a projective hyperplane of some projective space with some embedding of $\Gamma$ in that projective space. Let $k$ and $k'$ be two distinct non-collinear points of $\Gamma$, both not contained in $H$. Suppose  $k^{\perp} \cap H = k'^{\perp} \cap H$. Then there exists a point $p_{3}$ in $\Delta$ that is opposite both $p_{1}$ and $p_{2}$ with $p_{3}^{\perp} \cap \Gamma = H$, such that the odd projectivity $p_{1} \barwedge p_{2} \barwedge p_{3} \barwedge p_{1}$ fixes all lines joining $p_1$ with a point of $H$, and moves $p_1k$ to $p_1k'$.
\end{prop}

\begin{proof}
The fact that $p_{1} \barwedge p_{2} \barwedge p_{3} \barwedge p_{1}$ fixes all lines through $p_1$ having a point in $H$ implies that $p_3$ is collinear to all points of $H$. Also, the fact that $p_{1} \barwedge p_{2} \barwedge p_{3} \barwedge p_{1}$ maps $p_1k$ to $p_1k'$ implies that $p_3$ is on a line $L$ that intersects both $p_1k'$ and $p_2k$. We select an arbitrary point $k''$ on the line $p_1k'$ distinct from both $p_1$ and $k'$, and we let $L$ be the unique line through $k''$ intersecting $p_2k$. Set $K=k^\perp\cap H$. If $H$ contains lines, then $K$ is a hyperplane of $H$ and so $H$ is determined by $K$ and some point $h\in H\setminus K$, meaning that, if a point is collinear to $K$ and $h$, then it is collinear to $H$. If $H$ only contains points, then the same follows from our assumption that $H$ arises as the intersection of a projective hyperplane of some projective space with some embedding of $\Gamma$ in that projective space. 

We now define $p_3$ as the unique point of $L$ collinear to $h$. Since all of $p_1,p_2,k$ and $k'$ are collinear to $K$, also $k''$ and all other points of $L$ are collinear to $K$. In particular, $p_3$ is collinear to $K$, and since it is also collinear to $h$, it is collinear to $H$. 

It follows that $p_{1} \barwedge p_{2} \barwedge p_{3} \barwedge p_{1}$  fixes all lines joining $p_1$ with a point of $H$, and moves $p_1k$ to $p_1k'$.
\end{proof}

\begin{coro}\label{corkk'}
Let $\Delta$ be a thick polar space of rank $n\geq 3$ and let $H$ be a hyperplane of $\Delta$. Let $k,k'$ be two points not contained in $H$ such that $k^\perp \cap H=k'^\perp\cap H$. Then there exists a unique collineation $\theta$ of $\Delta$ pointwise fixing $H$ and mapping $k$ to $k'$.  
\end{coro}

\begin{proof}
First we prove uniqueness, 

Let $K$ be the set $k^{\perp} \cap H = k'^{\perp} \cap H$. Let $x$ be some point of $K$. Then the line $kx$ maps to the line $k'x=k^\theta x$. Let $a$ be some point on $kx$ not equal to $k$ or $x$. Then $a^{\theta}$ has to be on the line $k^{\theta}x$. Let $y$ be some point in $(a^{\perp} \cap H) \setminus (k^{\perp} \cap H)$. Then $a$ is the unique point of $kx$ collinear to $y$, and with that, $a^{\theta}$ is uniquely determined as the unique point on $k^{\theta}x$ collinear to $y^\theta=y$. That means the images of all points in $k^{\perp}$ and $H$ are uniquely determined. Playing the same game with all points of $k^\perp\setminus H$, and then again and again shows that $\theta$ is uniquely determined as the complement of a hyperplane in a polar space is always connected.

Next we prove existence. If $\Delta$ is non-embeddable, then by \cite{Coh-Shu:90}, the only hyperplanes are of the form $p^\perp$, for some point $p$, and then $\theta$ is a central elation (see \cite[Chapter~5]{Mal:24}). If $\Delta$ is embeddable, then we can view it as the intersection of $p_1^\perp$ and $p_2^\perp$, for two opposite points $p_1,p_2$ of a polar space of rank $n+1$. Then existence follows from   \cref{pointsgeneral}.
\end{proof}

\begin{remark}\label{symplPS}
It can happen that a (thick) polar space $\Delta$ of rank at least $3$ possesses a hyperplane $H$, but that there is no nontrivial collineation  of $\Delta$ pointwise fixing it. This is not in contradiction with \cref{corkk'}, as for such hyperplanes $H$ there do not exist two distinct points $k,k'$ with $k^\perp\cap H=k'^\perp\cap H$. This situation occurs for instance in symplectic polar spaces over fields of characteristic $2$. Indeed, Let $\Delta$ be a symplectic polar space of rank $r$ at least $3$ over a field $\K$ with $\kar\K=2$. The universal embedding of such a space corresponds with a (inseparable) quadric, and hence there are geometric hyperplanes of $\Delta$ which, as subsets of points of $\PG(2r-1,\K)$ (in which $\Delta$ is naturally embedded), generate $\PG(2r-1,\K)$. Let $H$ be such a  hyperplane. Select $2r-1$ points of $H$ that generate a hyperplane $J$ of $\PG(2r-1,\K)$.   Then $J=k^\perp$ for a unique point $k$ of $\Delta$. Hence there is no second point $k'$ with $k^\perp \cap H=k'^\perp\cap H$. 
\end{remark}

If we combine our results of \cref{special=general}, \cref{triangles}, \cref{lines} and \cref{pointsgeneral}, we obtain the following theorem. 

\begin{theorem}\label{theopoints}
Let $\Delta$ be a polar space of rank at least $3$ with at least $4$ points per line, and let $p$ be a point. Then $\Pi(x)$ is the subgroup of the automorphism group of $\Res(p)$ generated by all reflections of $\Res(p)$. If $\Delta$ is not a separable orthogonal polar space, then $\Pi^+(p)=\Pi(x)$. If $\Delta$ is a separable orthogonal polar space, then $\Pi^+(p)$ is the subgroup of $\Aut\Res(p)$ consisting of products of an even number of reflections.  
\end{theorem}

We have the following immediate consequence (taking into account that a polar space, which is not a separable quadric, automatically has at least $5$ points per line, except, if it is a symplectic polar space over $\mathbb{F}_2$).  

\begin{coro}\label{1=2}
Let $\Delta$ be an embeddable, but not separable orthogonal or symplectic polar space of rank at least $2$. Then the group of collineations generated by all reflections coincides with the group of collineation, which are the product of an even number of reflections. 
\end{coro}
\subsubsection{Hermitian case and symplectic polar spaces}
Noting that, in the case that $\Delta$ is embeddable, in the proof of \cref{pointsgeneral}, the hyperplanes $H$ arise as intersections of subspaces with $\Gamma$, we gather some immediate consequences of \cref{theopoints} in the following statements.

\begin{coro}\begin{compactenum}[$(i)$]
\item If $\Delta$ is the polar space arising from a non-degenerate Hermitian form (including symmetric bilinear forms) in a finite dimensional vector space over a (commutative) field, then, for each point $p$, the group $\Pi(p)$ is the full linear group preserving the Hermitian form that defines $\Res(p)$. 
\item  If $\Delta$ is a symplectic polar space, and $p$ is a point of $\Delta$, then $\Pi(p)=\Pi^+(p)$ is the simple symplectic group corresponding to $\Res(p)$, except, if $\Res(p)$ has rank $2$ and each line exactly has $3$ points, in which case the group is isomorphic to the symmetric group on $6$ letters. 
\end{compactenum}\label{herm-sympl}
\end{coro}

\begin{proof}
The bilinear case of $(i)$ follows directly from Propositions~ 8 and~14 of \cite{Die:73}. Now let the form that defines $\Delta$ be Hermitian with nontrivial field automorphism $\sigma$. Since $\Pi^+(p)$ contains the little projective group of $\Res_\Delta(p)$, \cite[Th\'eor\`eme~5]{Die:73} implies that it suffices to show that each field element $x$ with $xx^\sigma=1$ can be obtained as a determinant of a reflection. It suffices to consider the 2-dimensional case (vector space dimension), where this is immediate: the mapping $(x,y)\mapsto (x,ay)$, with $aa^\sigma=1$, preserves the Hermitian form $xx^\sigma+kyy^{\sigma}$, has determinant $a$ and fixes the hyperplane $(0,1)$ (and also the point $(1,0)$, which is the perp of $(0,1)$). Now $(i)$ follows.  

We now show $(ii)$. By \cref{symplPS}, $\Pi^+(p)$ is generated by reflections for which the fixed hyperplanes are point perps. These are elations and hence $\Pi^+(p)$ is the simple symplectic group corresponding to $\Res(p)$, if lines have at least $4$ points. By \cref{special=general}, $\Pi(p)=\Pi^+(p)$.    The statement $(ii)$ for polar spaces with exactly three points per line follows from \cref{MR0}, as the stabiliser of a point in the simple symplectic group is the symplectic group of the residue. (This argument can in fact also be used as an alternative for the larger symplectic polar spaces.)
\end{proof}

There remain three things to be addressed in more detail: (1) If $\Delta$ is a separable orthogonal polar space, can we be more specific about when $\Pi^+(p)=\Pi(p)$ (in other words, are there generic situations in which this equality \emph{always} or perhaps \emph{never} holds true)? (2) If $\Delta$ is non-embeddable, then the description of $\Pi(p)$ above is not very transparent; can we provide a description using the bilinear form defining the dual (embeddable) generalised quadrangle? (3) The case when lines of $\Delta$ have exactly three points. We begin with (1). 

\subsubsection{Separable orthogonal polar spaces}

\begin{prop}\label{+and-}
Let $\Delta$ be a separable, orthogonal polar space of rank $n$ with at least $4$ points per line, and anisotropic codimension $r$. Let $p$ be a point. If $r$ is even, then $\Pi^+(p)$ has index $2$ in $\Pi(p)$. 
If $r$ is odd, then $\Pi^+(p)=\Pi(p)$. 
\end{prop}

\begin{proof}
If $r$ is even, then each member of $\Pi^+(p)$ preserves the system of generators of the (imaginary) hyperbolic quadric over a splitting field of the quadric, whereas each single reflection interchanges these systems.  

Now let $r$ be odd. Since $\Pi(p)$ contains reflections with arbitrary spinor norm, in particular also non-trivial elements with trivial spinor norm, we find products of two reflections (which each are products of three perspectivities) with arbitrary spinor norm. The proposition follows. 
\end{proof}

\subsubsection{Polar spaces with $3$ points per line}
Now, we address (3). For symplectic polar spaces (which are isomorphic to parabolic polar spaces), this is already contained in \cref{herm-sympl}$(ii)$. Hermitian polar spaces never contain exactly three points per line.  There remain the so-called \emph{elliptic} ones, having anisotropic dimension $2$.  For these, \cref{pointsgeneral} implies that $\Pi(p)$ contains the special (projective) orthogonal group (which is generated by reflections by\cite[Proposition~14]{Die:73}). Now, it is clear that \cref{+and-} is also valid in the present case, and therefore $\Pi^+(p)$ is the simple orthogonal group (which has index $2$ in the special orthogonal group). 
\subsubsection{Non-embeddable polar spaces}
In order to deal with the thick non-embeddable polar spaces, we first take a look at residues isomorphic to the pseudo-quadratic polar spaces of rank $2$ that are dual to separable orthogonal polar spaces of rank $2$. We will need the notion of a similitude of a (quadratic or Hermitian) form.

A \emph{similitude} of a form is a linear transformation preserving the form up to a non-zero constant, and in case of a quadratic form, it is called \emph{direct} if the form is hyperbolic over a splitting field and the similitude preserves each natural system of maximal singular subspaces. 

Let $\K$ be a field with an involution $\sigma$, and let $\FF$ be the subfield pointwise fixed by $\sigma$. Let $\Delta(\K,\FF)$ be the pseudo-quadratic polar space of rank $3$ given by the pseudo-quadratic form \[X^\sigma_{-3}X_3+X^\sigma_{-2}X_2+X^\sigma_{-1}X_1\in\FF.\] In what follows, it is also allowed that $\K$ has characteristic $2$ and is an inseparable quadratic extension of $\FF$ (and then $\sigma$ is the identity).

Let $p_i$ be the base point corresponding to the $x_i$-coordinate, $i\in\{-3,-2,-1,1,2,3\}$. Let $p$ be the point with coordinates $(1,a_{-2},0,0,a_2,k-a_{-2}^\sigma a_2)$, with $a_{-2},a_2\in\K$ and $k\in\FF$.  Also, let $p(\ell_1,\ell_2)$ be the point with coordinates $(0,\ell_1,0,0,\ell_2,0)$, with $\ell_1,\ell_2\in\FF$. We assume that $p$ is not collinear to either $p_{-3}$ or $p_3$, that is, $k-a_{-2}^\sigma a_2\neq 0$. 

An elementary calculation shows that the projectivity $p_{-3}\per p_3\per p\per p_{-3}$ maps the plane $\<p_{-3},p_{-1},p(\ell_1,\ell_2)\>$ to the plane \[\<p_{-3},p_{-1},p((a_{-2}^\sigma a_2+a_2^\sigma a_{-2}-k)\ell_1-a_{-2}^\sigma a_{-2}\ell_2,a_2a_2^\sigma\ell_1-k\ell_2)\>.\]
It follows that, if $\Gamma$ is the separable, orthogonal polar space of rank $2$ dual to $p_{-3}^\perp\cap p_3^\perp$, and $L$ is the line of $\Gamma$ corresponding to the point $p_{-1}$, then the action on $L$ by the above projectivity, is given in binary coordinates by

\[\begin{pmatrix}\ell_1\\ \ell_2\end{pmatrix}\mapsto \begin{pmatrix} a_{-2}^\sigma a_2+a_2^\sigma a_{-2}-k & -a_{-2}^\sigma a_{-2} \\ a_2^\sigma a_2 & -k\end{pmatrix}\cdot\begin{pmatrix}\ell_1\\ \ell_2\end{pmatrix}=:A\cdot\begin{pmatrix}\ell_1\\ \ell_2\end{pmatrix}.\] The determinant of the $2\times 2$-matrix $A$ in the above expression is \[(k-a_{-2}^\sigma a_2)\cdot(k-a_{-2}^\sigma a_2)^\sigma.\]

Let $\Gamma$ be embedded in $\PG(5,\FF)$. Note, that the fixed points of the projectivity above form an ovoid of $p_{-3}^\perp\cap p_3^\perp$ and hence, a spread of $\Gamma$. Taking six base points on three spread lines, amongst which we can choose $L$, we claim that the determinant of the matrix of the corresponding collineation is \[[(k-a_{-2}^\sigma a_2)\cdot(k-a_{-2}^\sigma a_2)^\sigma]^3.\]
Indeed, let $M$ be a second line with two base points. We can choose the basis such that the equation of $\Gamma\cap\<L,M\>$ is $x_0x_2+x_1x_3=0$, where $L=\<(1,0,0,0),(0,1,0,0)\>$ and $M=\<(0,0,1,0),(0,0,0,1)\>$. The projectivity now fixes every line of $\Gamma\cap\<L,M\>$ disjoint from $L$ (and also $L$).  One  calculates that the matrix of the collineation is the block matrix  \[\begin{pmatrix}A & 0 \\ 0 & A\end{pmatrix}.\] Doing this for the third line with base points again, the claim follows. 

Hence, we conclude that special and general projectivity group of the upper residue of $p_{-3}$, viewed as collineation group of $\Gamma$, consist of all similitudes of the associated quadratic form, with factor a norm of the quadratic extension $\K/\FF$ (note that each square of $\FF$ is such a norm). Writing the equation (in some other basis) of $\Gamma$ as $y_{-2}y_2+y_{-1}y_1=yy^\sigma$, where we view the underlying vector space as $\FF\times\FF\times\K\times\FF\times\FF$, with generic coordinates $(y_{-2},y_{-1},y,y_1,y_2)$, we see that the factor of each similitude is a norm. 

We can now state:

\begin{prop}
Let $\Delta$ be the thick non-embeddable polar space associated with the Cayley division algebra $\O$ over the field $\FF$ with standard involution $x\mapsto\overline{x}$. Let $p$ be any point of $\Delta$. Then $\Pi(p)=\Pi^+(p)$ is the group of all direct similitudes of the quadratic form \[q:\FF\times\FF\times\O\times\FF\times\FF\to\FF:(x_{-2},x_{-1},x,x_1,x_2)\mapsto x_{-1}x_1+x_{-2}x_2+x\overline{x}, \] and the factor of such similitude is a norm of some element of $\O$.  
\end{prop}

\begin{proof}
Let $p$  be a point of $\Delta$. The rank $2$ polar space $\Res_\Delta(p)$ is dual to a separable quadric of Witt index $2$ in $\PG(11,\FF)$, given by the quadratic form mentioned in the statement. Clearly, $\Pi(p)=\Pi^+(p)$ is contained in the group of all direct similitudes of that quadratic form. So, it sufficed to prove that each such similitude can occur. But this follows from the discussion above taking into account that the pseudo-quadratic polar space   $\Delta(\K,\FF)$, with $\K$ any $2$-dimensional subfield of $\O$ containing $\FF$, is a sub polar space of $\Delta$, with the property that all planes of $\Delta$ through a line of $\Delta(\K,\FF)$ belong to $\Delta(\K,\FF)$. 
\end{proof}

\subsection{Subspaces of dimension at least $1$.}

\begin{lem}\label{atleast1herm}
Let $\Delta$ be a polar space of rank $n \geq 3$ embedded into $\PG(V)$ that is not a separable quadric. Let $A_{1}, A_{2}, A_{3}$ and $A_{4}$ be singular subspaces of dimension $m \leq n-2$, such that $A_{1} \cap A_{3} =: B$ and $A_{2} \cap A_{4} =: C$, with $B$ and $C$ opposite and of dimension $m-1$, such that $A_{i}$ is opposite $A_{i+1}$ for $i \in \mathbb{Z} / 4 \mathbb{Z}$ and such that $A_{1}$ and $A_{3}$ are not contained in a common subspace. Then $A_{1} \barwedge A_{2} \barwedge A_{3} \barwedge A_{4} \barwedge A_{1}$ is a product of two collineations in $\Res_{\Delta} (A_{1})$ that each fix a hyperplane.
\end{lem}

\begin{proof}We may assume that $\PG(V)$ is minimal, that is, the relation $\perp$ defines a non-degenerate polarity of $\PG(V)$. 
 
(I) We will first assume that $m \leq n-3$. Then $A_{1}^{\perp} \cap A_{2}^{\perp}$ has rank $n-m-1$. Since $n-m-1 \geq 2$, $A_{1}^{\perp} \cap A_{2}^{\perp} \cap A_{3}^{\perp}$ is a hyperplane of $A_{1}^{\perp} \cap A_{2}^{\perp}$ 
that we will denote by $H$. If $A_{4}^{\perp} \cap A_{1}^{\perp} \cap A_{2}^{\perp} = H$, then $H$ is fixed pointwise. So suppose $A_{4}^{\perp} \cap A_{1}^{\perp} \cap A_{2}^{\perp} \neq H$. Define: 
\begin{align*}
x_{1} := C^{\perp} \cap A_{1}, \hspace{.4cm}&
x_{2} := B^{\perp} \cap A_{2},\\
x_{3} := C^{\perp} \cap A_{3}, \hspace{.4cm}&
x_{4} := B^{\perp} \cap A_{4}
\end{align*}
Since $B$ and $C$ are opposite, $x_1\neq x_3$ and $x_2\neq x_4$.
In $\PG(V)$, there are planes $\langle x_{1}, x_{2}, x_{3} \rangle$ and $\langle x_{1}, x_{3}, x_{4} \rangle$ which share the line $\langle x_{1}, x_{3} \rangle$ in $\PG(V)$. On the line $\langle x_{1}, x_{3} \rangle$ in $\PG(V)$, there exists some other point that is contained in $\Delta$, since $\Delta$ is not a separable quadric, and we denote it by $x_{5}$. Since every point of $C$ is collinear to $x_{1}$ and $x_{3}$, it follows that every point of $C$ is collinear to $x_{5}$. The subspace $\langle C, x_{5} \rangle$ is a subspace in $\Delta$ and we denote it by $A_{5}$. Now $A_{5}^{\perp} \cap A_{1}^{\perp} \cap A_{2}^{\perp} = H$ and with that, we can write $\varphi$ as the product of the projectivities 
\[\varphi_{1} \colon A_{1} \barwedge A_{2} \barwedge A_{3} \barwedge A_{5} \barwedge A_{1}\mbox{ 
and }\varphi_{2} \colon A_{1} \barwedge A_{5} \barwedge A_{3} \barwedge A_{4} \barwedge A_{1}.\]
Then $\varphi_{1}$  fixes $H$ pointwise and, similarly, $\varphi_{2}$ fixes a hyperplane of $\Res_{\Delta} (A_{1})$ pointwise.  Therefore, $\varphi = \varphi_{1} \varphi_{2}$ is the product of two collineations in $\Res_{\Delta} (A_{1})$ that each fix a hyperplane pointwise.

(II) Now suppose $m = n-2$. Then the intersection $A_{1}^{\perp} \cap A_{2}^{\perp} \cap A_{3}^{\perp}$ could be empty. In that case, we consider the subspace of $\PG(V)$ spanned by $A_3^{\perp}$ and intersect it with the subspace spanned by $A_{1}^{\perp} \cap A_2^{\perp}$. This intersection will be a hyperplane of the subspace spanned by $A_{1}^{\perp} \cap A_2^{\perp}$ and we take this as $H$. If the hyperplane spanned by $A_{4}^{\perp}$ intersects the space spanned by $A_{1}^{\perp} \cap A_2^{\perp}$ in $H$ as well, we are done. Otherwise, we do the same construction as before and obtain a subspace $A_{5}$, such that the hyperplane spanned by $A_{5}^{\perp}$ intersects the space spanned by $A_{1}^{\perp} \cap A_2^{\perp}$ in $H$ and $\varphi = \varphi_{1} \varphi_{2}$ as before with the same definitions for $\varphi_{1}$ and $\varphi_{2}$, and the same conclusion holds. 
\end{proof}

\begin{lem}\label{atleast1orth}
Let $\Delta$ be a polar space of rank $n \geq 3$ embedded into $\PG(V)$ that is a separable quadric. Let $A_{1}, A_{2}, A_{3}$ and $A_{4}$ be singular subspaces of dimension $m \leq n-2$, such that $A_{1} \cap A_{3} =: B$ and $A_{2} \cap A_{4} =: C$, with $B$ and $C$ opposite and of dimension $m-1$, such that $A_{i}$ is opposite $A_{i+1}$ for $i \in \mathbb{Z} / 4 \mathbb{Z}$ and such that $A_{1}$ and $A_{3}$ are contained in a common subspace and $A_{2}$ and $A_{4}$ are contained in an opposite subspace. Then $\theta \colon A_{1} \barwedge A_{2} \barwedge A_{3} \barwedge A_{4} \barwedge A_{1}$ is a product of two collineations in $\Res_{\Delta} (A_{1})$ that each fix a hyperplane in the ambient projective space.
\end{lem}

\begin{proof}
Since $\Delta$ is a separable quadric, we can view the perp relation as a non-degenerate polarity.

Let $\Gamma := A_{1}^{\perp} \cap A_{2}^{\perp}$. Then $\Gamma$ has rank $n-m-1$. Since $m \leq n-2$, the rank of $\Gamma$ is at least $1$. Let $H$ be  the hyperplane $A_{3}^{\perp} \cap \<\Gamma\>$ and $G$ be the subhyperplane $A^\perp_{4} \cap H$. If $G=H$, we are done. So suppose $G$ is a proper subset of $H$.

Define:
\begin{align*}
x_{1} := C^{\perp} \cap A_{1}, \hspace{.4cm}&
x_{2} := B^{\perp} \cap A_{2},\\
x_{3} := C^{\perp} \cap A_{3}, \hspace{.4cm}&
x_{4} := B^{\perp} \cap A_{4}
\end{align*}

Let $c_{1}$ be the point $A_{2}^{\perp} \cap \langle A_{1}, A_{3} \rangle$ and $c_{2}$ be the point $A_{1}^{\perp} \cap \langle A_{2}, A_{4} \rangle$. Note that both $c_1$ and $c_2$ belong to $\Gamma$. We have $ A_{2}^{\perp} \cap \langle A_{1}, A_{3} \rangle \subseteq C^{\perp} \cap \langle A_{1}, A_{3} \rangle = x_{1}x_{3}$. That means $c_{1}$ is on the line $x_{1}x_{3}$ and analogously $c_{2}$ is on the line $x_{2}x_{4}$. Every point of $G$ is collinear to $x_{1}$ and $x_{3}$ and with that also to $c_{1}$ and analogously to $c_{2}$. That means $G \subseteq c_{1}^{\perp} \cap c_{2}^{\perp}$, and since $G$ is a subhyperplane of $\Gamma$, we have $G=\<\Gamma\>\cap c_1^\perp\cap c_2^\perp$. We may view $\theta$ as a collineation in $\Gamma$; the fixed points correspond precisely with the points of $G$.

If $c_1\perp c_2$, then the singular subspace generated by $A_1$ and $A_3$ is not opposite the one generated by $A_2$ and $A_4$, contradicting our hypothesis. So $c_1$ and $c_2$ are not collinear. Then we can apply \cref{rank1} and the proof is complete.
\end{proof}

\begin{lem}\label{oddeven}
Let $\Delta$ be a polar space of rank $r\geq 3$, which is a separable quadric in $\PG(V)$. Let $1\leq d\leq r-2$. Let $U_1$ and $U_2$ be two opposite singular subspaces of dimension $d$.  Let $W_1$ and $W_2$ be two singular subspaces of dimension $d+1$, containing $U_1$ and $U_2$, respectively, and intersecting in a point $p$. Let $U_3$ be a singular subspace of dimension $d$ intersecting $\<U_1,U_2\>$ in a $(d-1)$-dimensional subspace $H$. Suppose $p$ is not collinear to all points of $U_3$. Then the image of $W_2$ under $U_2\per U_3\per U_1$ is $W_1$ if, and only if, $d$ is odd.  
\end{lem}

\begin{proof}
Set $u_i=\proj_{U_i}H$, $i=1,2$. Since $U_1,U_2\subseteq p^\perp$, we find that $H\subseteq p^\perp$. It follows that the line $pu_i$ is the projection of $H$ onto $W_i$, $i=1,2$. Consequently, $\proj_{W_i}U_3=:x_i\in pu_i$, $i=1,2$. Note that by assumption, $p\notin\{x_1,x_2\}$. We conclude that the image of $W_2$ under $U_2\per U_3\per U_1$ is $W_1$ if, and only if,  $x_1\perp x_2$ (the ``if''-part follows from the fact that the $(2d+1)$-dimensional subspace of $\PG(V)$ generated by $U_1$ and $U_2$ does not contain singular subspaces of dimension $d+1$ as it intersects $\Delta$ in a non-degenerate hyperbolic quadraic). This is equivalent to $u_1\perp u_2$, which, on its turn, is equivalent to $u_1\in\<H,u_2\>$. This happens if, and only if, $U_1$ and $U_2$ belong to the same natural system of maximal singular subspaces of the hyperbolic quadric obtained by restricting $\Delta$ to $\<U_1,U_2\>$. Since hyperbolic quadrics contain opposite maximal singular subspaces of the same natural system if, and only if, their rank is even (and equivalently, the dimension of the maximal singular subspaces is odd), the assertion follows.    
\end{proof}

This has now the following consequences.

\begin{theorem}\label{theosubspaces}
Let $\Delta$ be a polar space of rank $r$ at least $3$ with at least $4$ points per line, and let $U$ be a singular subspace of dimension at most $r-2$. Then $\Pi^+_\geq(U)$ is the subgroup of the automorphism group of $\Res(U)$ generated by products of two reflections of $\Res(U)$. If $\Delta$ is not a separable orthogonal polar space, then $\Pi_\geq^+(U)=\Pi_\geq(U)$ is also generated by all reflections. If $\Delta$ is a separable orthogonal polar space, then $\Pi_\geq(U)=\Pi_\geq^+(U)$, if $\dim U$ is odd, and $\Pi(U)$ is the subgroup of the automorphism group of $\Res(U)$ generated by all reflections, if $\dim U$ is even.  In the latter case, $\Pi_\geq(U)$ is not equal to $\Pi_\ge^+(U)$, as soon as the anisotropic dimension of $\Delta$ (as a quadric), is even.   
\end{theorem}

\begin{proof}
Let $G_1$ be the group of collineations of $\Res(U)$ generated by the reflections, and let $G_2\leq G_1$ be its subgroup consisting of those members, which are a product of an even number of reflections. 

By \cref{gate2} and \cref{theopoints}, $G_2\leq\Pi^+_\geq(U)$. By \cref{atleast1herm} and \cref{atleast1orth}, we have $\Pi^+_\geq(U)\leq G_2$. This proves the first part of the theorem.

 If $\Delta$ is not a separable orthogonal polar space, then $\Pi_\geq^+(U)=\Pi_\geq(U)$ by \cref{special=general}.  Then the second assertion follows directly from \cref{1=2}. 
 
Now let $\Delta$ be a separable orthogonal polar space. If $\dim U$ is odd, then \cref{special=general} implies $\Pi_\ge(U)=\Pi_\ge^+(U)$. Now let $\dim U$ be even. Then \cref{oddeven} implies that there exists a non-trivial self-projectivity of $U$ of length $3$ that fixes a hyperplane of $\Res(U)$. It follows now that $\Pi_\ge(U)$ is generated by reflections. Finally, if   the anisotropic dimension of $\Delta$ (as a quadric) is even, then this is also the case for $\Res(U)$. Then each member of $\Pi_\ge^+(U)$ preserves the system of generators of the (imaginary) hyperbolic quadric over a splitting field of the quadric, whereas each single reflection interchanges these systems.  
\end{proof}

The projectivity groups for upper residues in symplectic polar spaces follow directly from \cref{MR0}: 

\begin{prop}
If $\Delta$ is a symplectic polar space of rank $r$, and $U$ is a singular subspace of $\Delta$ of dimension $d\leq r-2$, then $\Pi_\geq(U)=\Pi_\geq^+(U)$ is the simple symplectic group corresponding to $\Res(U)$, except if $\Res(U)$ has rank $2$ and each line has exactly $3$ points, in which case the group is isomorphic to the symmetric group on $6$ letters.
\end{prop}

Similarly as for points, we can also handle the case of elliptic polar space with three points per line. 

\begin{prop}
Let $\Delta$ be a polar space of rank $r$ with three points per line and distinct from a symplectic polar space. Let $U$ be a singular subspace of dimension $d\leq r-2$. Then $\Pi_\geq^+(U)$ is the simple orthogonal group associated to $\Res(U)$.  If $d$ is odd, then $\Pi_\geq(U)=\Pi^+_\geq(U)$, whereas, if $d$ is even, $\Pi_\geq(U)$ is the special orthogonal group associated to $\Res(U)$.
\end{prop}

A special case of \cref{theosubspaces} is worth mentioning.

\begin{coro}\label{lemNEPS}
Let $\Delta$ be a $\sigma$-Hermitian polar space over a commutative field $\K$ with anisotropic dimension $0$, and let $\mathbb{F}$ be the fix field of $\sigma$. Let $U$ be a submaximal singular subspace. Then $\Pi^+_\geq(U)=\Pi_\geq(U)$ is permutation equivalent to $\PGL_2(\mathbb{F})$. 
\end{coro}

Since each non-embeddable polar space contains an ideal $\sigma$-Hermitian polar subspace (\emph{ideal} means that every plane of the large space through a line of the subspace belongs entirely to the subspace), we immediately deduce from \cref{lemNEPS} the following result,

\begin{coro}\label{linesNEPS}
Let $\Delta$ be a non-embeddable polar space with planes isomorphic to projective planes over an octonion division algebra with centre $\K$. Let $L$ be a line of $\Delta$. Then $\Pi^+_\geq(L)=\Pi_\geq(L)$ is permutation equivalent to $\PGL_2(\mathbb{\K})$. 
\end{coro}

\section{Projectivity groups for projective residues in polar spaces}\label{projres}

\subsection{Non-maximal projective residues} 

\subsubsection{Embeddable polar spaces} This case is straightforward, using \cref{gate2}. Indeed, the special projectivity group of a subspace $D$ of dimension $d$ of a Desarguesian projective space over the skew field $\L$ is the full linear group $\PGL(d,\L)$, so it also coincides with the special projectivity group inherited from the polar space.  For the general projectivity group, we have to add a duality, which is linear, if the polar space is orthogonal or symplectic, but which is semi-linear with companion skew field automorphism $\sigma$, if the polar space is associated to a $\sigma$-quadratic form. If $\dim D=1$, then this duality is just a collineation of $\PG(1,\L)$ (with $\L$ the underlying skew field) induced by $\sigma$. \label{embnm} In the orthogonal case, the duality belongs to $\PGL_2(\L)$ and so the general and special groups coincide; in the other cases the special group has index 2 in the general group. 

\subsubsection{Non-embeddable polar spaces} Here, the residues have type $1$ or $2$ in Bourbaki labelling. We look at the lower residue of a line, and at the planar line pencils. We first consider the lower residues of a line. 

\begin{prop}\label{lowerlines}
Let $\Delta$ be the non-embeddable polar space associated to the Cayley division algebra $\O$ over the field $\K$. Let $L$ be a line of $\Delta$. Then $\Pi_{\leq}^+(L)$ is permutation equivalent to the group of direct similitudes of the quadric of Witt index $1$, defined by the anisotropic form given by the norm form of $\O$ over $\K$.  Also, $\Pi_{\leq}(L)$ is isomorphic to the group generated by $\Pi^+_{\leq}(L)$ and a permutation induced by the standard involution of $\O$, when considering $L$ as a projective line over $\O$.  
\end{prop}

\begin{proof}
We begin by noting that the mentioned permutation group is equal to the projectivity group of a line in the Moufang plane $\PG(2,\O)$, as proved by Grundh\"ofer \cite[Satz]{Gru:83}.
So, in view of \cref{gate2},  it suffices to show that every (lower) even projectivity of $L$ in $\Delta$ coincides with a projectivity inside some plane of $\Delta$ containing $L$.  Thanks to \cref{upanddown}, \cref{openneer}, \cref{lines} and \cref{higherrank}, it suffices to show this for lower projectivities of the form $L\per L_2\per L_3\per L_4\per L$, where $L$ and $L_3$ intersect in a point $p_1$ and generate a plane $\pi_1$, where $L_2$ and $L_4$ intersect in a point $p_2$ and generate a plane $\pi_2$, and where $p_1$ and $p_2$ are opposite, and $\pi_1$ and $\pi_2$ are also opposite. Set $c_i=\proj_{\pi_1}L_i$, $i=2,4$.  
Then one observes that, for each $x\in L$,  $x_3:=\proj_{L_3}\proj_{L_2}x$ is contained in the line $xc_2$, and $x':=\proj_L\proj_{L_4}x_3$ lies on the line $x_3c_4$. Hence, $L\per L_2\per L_3\per L_4\per L$ coincides with the projectivity $L\per c_2\per L_3\per c_4\per L$ inside $\pi_1$. This completes the proof of the first part of the proposition.

 The second part follows from \cref{embnm} by noting that $\Delta$ contains (many) polar spaces over quaternion subfields.  
\end{proof}

\begin{prop}\label{NEpencil}
Let $\Delta$ be the non-embeddable polar space associated to the Cayley division algebra $\O$ over the field $\K$. Let $P$ be a planar line pencil of $\Delta$. Then $\Pi^+(P)$ is permutation equivalent to the group of direct similitudes of the quadratic form of Witt index $1$ defined by the norm form of $\O$ over $\K$.  Also, $\Pi(P)$ is isomorphic to the group generated by $\Pi^+(L)$ and the permutation induced by the standard involution of $\O$, when considering $P$ as a projective line over $\O$.  
\end{prop}

\begin{proof} We identify a line pencil with its point-plane pair. 
As in the previous proof, it suffices to show that, for line pencils $P_2,P_3,P_4$, with $P\sim P_3$ and $P_2\sim P_4$ (see the next paragraph for more details on the adjacency $\sim$) 
such that both $P$ and $P_3$ are opposite both $P_2$ and $P_4$, the projectivity $P\per P_2\per P_3\per P_4\per P$ is a projectivity of $P$ inside the plane it defines. 

Assume that $P=\{p,\pi\}$, with $p$ a point in the plane $\pi$, and likewise $P_i=\{p_i,\pi_i\}$, $i=2,3,4$. Then $P\sim P_3$ means that either $\pi=\pi_3$ and $p\neq p_3$, or $p=p_3$ and $\pi\cap\pi_3$ is a line. Suppose first that $\pi=\pi_3$. Let $L$ be the line $pp_3$ and set $M:=\proj_\pi p_2$. Note $M\cap L\notin\{p,p_3\}$.  Let $K$ be an arbitrary line of $P$. Then its image under $P\per P_2$ is the line $K_2$ through $p_2$ and $\proj_{\pi_2}K$. It follows that $\proj_\pi K_2=K\cap M$. But then the image $K_3$ of $K_2$ under $P_2\per P_3$ is the line joining $p_3$ with  $K\cap M$. It follows that $P\per P_2\per P_3$ coincides with the projectivity $p\per M\per p_3$ in the projective plane $\pi$. 

Now assume $p=p_3$ and $\pi\cap\pi_3$ is a line $L\in P\cap P_3$. Set $M=\proj_\pi p_2$, $M_3=\proj_{\pi_3}p_2$ and $M_2=\proj_{\pi_2}p$.  Let $K$ again be an arbitrary line of $P$ and set $x=K\cap M$, $x_2= \proj_{\pi_2}K=\proj_{M_2}x$, $K_2=p_2x_2$, $x_3=\proj_{\pi_3}K_2=\proj_{M_3}x_2$, and $K_3=px_3$. Then we find that $K_3$ is the image of $K$ under $P\per P_2\per P_3$, and $x_3$ is the image of $x$ under $M\per M_2\per M_3$ (lower perspectivities). Hence, $P\per P_2\per P_3$ can be written as the product of $p\per M$ in $\pi$, with $M\per M_2\per M_3$ in $\Delta$, and finally $M_3\per p$ in $\pi_3$. Since a similar decomposition holds for $P_3\per P_4\per P$, we conclude, using  \cref{lowerlines},  that $P\per P_2\per P_3\per P_4\per P$ is contained in the projectivity group of $p$ inside $\pi$. The first part of the proposition is proved. 

The second part follows from the second part of \cref{lowerlines} by noting that, if $Q=\{q,\alpha\}$ is an arbitrary line pencil opposite $P$, the perspectivity $P\per Q$ is the product of the perspectivity $p\per \proj_{\pi}q$ inside $\pi$, the lower perspectivity $\proj_{\pi}q\per\proj_{\alpha}p$ and the perspectivity $\proj_{\alpha}p\per q$ inside the plane $\alpha$. 
\end{proof}

\begin{remark}\label{similitudes}
The standard involution of a split octonion algebra interchanges the two systems of maximal singular subspaces, as follows immediately from the fact that it interchanges two of the eight standard basis vectors and maps the others to their additive inverse (in the representation of Zorn). Hence, the general projectivity groups in  \cref{lowerlines} and \cref{NEpencil} are the groups of all similitudes of the said quadratic forms. This also holds for the quadratic associative division algebras over the field $\K$ with non-trivial standard involution. For further use, given a quadratic alternative division algebra $\A$ over $\K$ with non-trivial standard involution, we denote the group of all permutations of the projective line $\PG(1,\A)$ corresponding to direct similitudes by $\PGL_2^+(\A)$, and the group of all permutations corresponding to all similitudes by $\PGL_2(\A)$, except, if $\A$ is commutative (and this notation would be ambiguous), then we use $\PGL^+_2(\A/\K)$ and $\PGL_2(\A/\K)$, respectively.
\end{remark}

\subsection{Maximal singular subspaces}

\subsubsection{Reduction to the composition of four perspectivities}

\begin{prop}\label{MSSred}
Let $\Delta$ be a polar space of rank $r\geq 3$. Let $M$ be a maximal singular  subspace. of $\Delta$. Suppose the set of projectivities $M\per M_2\per M_3\per M_4\per M$, with $M\equiv M_2\equiv M_3\equiv M_4\equiv M$ and $U_0:=M\cap M_2$ and $U_1:=M_1\cap M_3$ opposite submaximal singular subspaces, is geometric. Then $\Pi^+(M)$ is generated by all such projectivities. 
\end{prop}

\begin{proof}
This follows from \cref{openneer} and \cref{higherrank}.
\end{proof}

As we will see below, the projectivities $M\per M_2\per M_3\per M_4\per M$, mentioned in the previous proposition, are homologies. So proving that the said set is geometric is equivalent to showing that the set of factors of these homologies  are closed under conjugation with any element of the underlying skew field.

\subsubsection{Computation of the special projectivity group in the general (non-symplectic embeddable) case}\label{SPG}
\def\L{\mathbb{L}}
\def\q{\mathfrak{q}}
We use the standard form, see \cite[Chapter~4]{Mal:24}. Let $\Delta$ be an embeddable polar space with pseudo-quadratic description in $\PG(V)$, for some right vector space $V$ of some skew field $\L$ and let $\sigma$ be an involution of $\L$. Let $V_0$ be a subspace of $V$ of codimension $2n$ and let 
$$\{e_i\mid i\in\{-n,-n+1,\ldots,-1,1,\ldots,n\}\}$$
 be a (suitably chosen) basis of a complementary subspace, such that 
$$\{\<e_i\>\mid  i\in\{-n,-n+1,\ldots,-1,1,\ldots,n\}\}$$
is a polar frame. For general vectors $v,w\in V$ we write $v=v_0 + \sum_{i\in I} e_ix_i$ and $w=w_0+\sum_{i\in I} e_iy_i$ (where $v_0,w_0\in V_0$ and $x_i,y_i\in \L$ for all $i\in I$). Then there exists an anisotropic $\sigma$-quadratic form $\q_0$ on $V_0$, where $\q_0(v_0)=g_0(v_0,v_0)+\L_\sigma$, for the $(\sigma,1)$-linear form $g_0:V_0\times V_0\to\L$, with associated $\sigma$-Hermitian form $f_0$, such that    
\begin{align*}
\q(v)&=x_{-n}^\sigma x_n +\cdots + x_{-1}^\sigma x_1 + \q_0(v_0),\\
f(v,w)&=x_{-n}^\sigma y_n + x_{n}^\sigma y_{-n}+ \cdots +x_{-1}^\sigma y_1 + x_{1}^\sigma y_{-1} + f_0(v_0,w_0)\end{align*}
define the points and collinearity, respectively, of the polar space $\Delta$.   

Let $v_0,w_0\in V_0$, $t,u\in\L$ be such that 
\[D:=D(v_0,w_0,t,u):=(g_0(w_0,w_0)+u-u^\sigma)(g_0(v_0,v_0)+t-t^\sigma)+f_0(w_0,v_0)+1\neq 0.\] We abbreviate $g_1:=g_0(v_0,v_0)+t-t^\sigma$ and $g_2:=g_0(w_0,w_0)+u-u^\sigma$ and define the following four maximal singular subspaces.
\[\begin{cases}
M_1=\<e_{-1},e_{-2}.\ldots,e_{-n}\>,\\
M_2=\<e_1,e_2,\ldots,e_n\>,\\
M_3=\<g_1e_1+v_0+e_{-1},e_{-2},\ldots,e_{-n}\>,\\
M_4=\<g_2^\sigma e_{-1}+w_0+e_1, e_2,\ldots,e_n\>.
\end{cases}\]
Then, by our condition above, we have $M_1\equiv M_2\equiv M_3\equiv M_4\equiv M_1$ and $M_1\sim M_3$ and $M_2\sim M_4$ (cf.~\cref{defsim}). Set $\theta:=M_1\per M_2\per M_3\per M_4\per M_1$.  Then $\theta$ pointwise fixes the hyperplane $\<e_{-2},\ldots,e_{-n}\>$ of $M_1$ and also the point $\<e_{-1}\>$. Hence, $\theta$ is a homology, and we determine its factor. We can work in the subspace  $W$ generated by $e_{-2},e_{-1},e_1,e_2$ and $V_0$. We write a generic element of that subspace as $(x_{-2},x_{-1},x_0,x_1,x_2)$, where all elements belong to $\L$, except $x_0\in V_0$. 

We consider a generic point $(1,x,0,0,0)$ on $\<e_{-2},e_{-1}\>$ distinct from $e_{-2}$, so we may assume $x\neq 0$. Its image in $M_2\cap W$ under the perspectivity $M_1\per M_2$ is the point $(0,0,0,1,- x^\sigma)$.  The image of that point in $M_3\cap W$ under the perspectivity $M_2\per M_3$ is the point $(x^{-1},1,v_0,g_1,0)$. Likewise, the image of the latter in $M_4\cap W$, under the perspectivity $M_3\per M_4$, is the point $(0,g_2^\sigma,w_0,1,-x^{\sigma} D^\sigma)$. Finally, the image of that point in $M_1\cap W$, back again under the perspectivity $M_4\per M_1$, is the point $(1,Dx,0,0,0)$.

We can now formulate the main result of this subsection.  The set 
$$\{D(v_0,w_0,t,u)\mid v_0,w_0\in V_0, t,u\in\L\}$$
 will be referred to as the \emph{norm set}. 

\begin{prop}\label{normset}
The special projectivity group of a maximal singular subspace of an embeddable polar space with pseudo-quadratic description in standard form as above, is generated by homologies with factors in the norm set. 
\end{prop}

\begin{proof}
In view of the preceding computations and \cref{MSSred}, it only remains to show that the norm set is closed under conjugation with an arbitrary element $r\in\L^\times$. This follows straight from the identities $r^{-1}f_0(w_0,v_0)r=f_0(w_0r^{-\sigma},v_0r)$, $r^\sigma(t-t^\sigma)r=(r^\sigma tr)-(r^\sigma t r)^\sigma$.
\end{proof}

It strongly depends on $g_0$, what exactly the special projectivity group of a maximal singular subspace is. We present a few special cases.

For a parabolic polar space of rank $r$, $V_0=\L$, $\sigma=\id$ and $g_0(x_0,y_0)=x_0y_0$. It follows that $f_0(x_0,y_0)=2x_0y_0$, and consequently, $D(x_0,y_0,t,u)= x_0^2y_0^2+ 2x_0y_0+1=(x_0y_0+1)^2$. Hence, the special projectivity group here is $\PGL_r(\L)$, if $r$ is odd, and is the linear subgroup of $\PGL_r(\L)$, consisting of the matrices with a square determinant, if $r$ is even.

 Now, suppose $V_0=\L\oplus\L$, $\ch\L\neq 2$ and $\sigma=\id$. Let $g_0$ be given by $g_0((x,y),(x'y'))=xx'+\ell yy'$. Then $D((x,y),(x',y'),t,u)=(1+xx'+yy')^2+\ell(xy'-x'y)^2$ and so the special projectivity group is generated by homologies with factor a norm of the quadratic extension of $\L$ corresponding to $g_0$. In the finite case, every element of $\L$ is a norm, and then we simply have $\PGL_r(q)$, $q=|\L|$. 
 
 For Hermitian polar spaces over commutative fields with non-trivial involution, see the next paragraph. 
 
 \subsubsection{Special projectivity groups for symplectic polar spaces and some Hermitian ones} The previous subsection does not cover the symplectic polar spaces (it would, if we considered skew-Hermitian forms rather than Hermitian forms). However, in this case, the simple group (generated by the root groups) is  easy to describe and hence, we can directly use \cref{MR0}. The same method can be applied to Hermitian polar spaces over commutative fields with non-trivial field involution. 
 
 \begin{prop}\label{MSSsymplectic}
 Let $\Delta$ be a symplectic polar space of rank $r$, defined over the field $\K$. Let $U$ be a maximal singular subspace. Then $\Pi^+(U)\cong \PGL_r(\K)$ in its standard permutation representation. Also, $\Pi(U)\cong\PGL_r(\K)\rtimes \<\rho\>$, with $\rho$ any linear polarity. 
 \end{prop}
 
 \begin{proof}
 Let $\Delta$ be described by the alternating form \[(x_{-r}y_r+x_{-r+1}y_{r-1}+\cdots x_{-1}y_1) - (y_{-r}x_r+y_{-r+1}x_{r-1}+\cdots y_{-1}x_1).\]
 Using obvious notation, the subspace $U$, spanned by $e_{-1},\ldots, e_{-r}$, is a maximal singular subspace, and the subspace $W$ generated by $e_1,\ldots,e_r$ is an opposite maximal singular subspace. Set $X_{-}=(x_{-1},\ldots,x_{-r})$  and $X_+=(x_1,\ldots,x_r)$, which we also read as row matrices, then, for every non-singular $r\times r$ matrix $M$, the linear transformation \[X_{-}\mapsto X_{-}\cdot M^t, \hspace{1cm} X_{+}\mapsto X_{+}\cdot M^{-1}\] induces a linear collineation in $\PG(2r-1,\K)$, preserving the above alternating form. Since such a collineation of $\Delta$ belongs to the automorphism group generated by the root elations, as follows from \cite[Th\'eor\`eme~1]{Die:73}, we find, using \cref{MR0}, that $\Pi^+(U)\cong \PGL_r(\K)$. Now a single perspectivity induces a linear duality (over $\K$) from one maximal singular subspace to another, concluding the proof of the proposition. 
 \end{proof}

A similar argument can be given to obtain a slightly more concrete version of  \cref{normset} in case of Hermitian polar spaces over a commutative field with non-trivial involution. For a given field $\K$ and subfield $\FF$, let $\SL_r(\K;\FF)$ be the multiplicative group of $r\times r$ matrices with a determinant in $\FF^\times$, and let $\PSL_r(\K;\FF)$ be the corresponding projective group, that is, the quotient group with its centre. 

 \begin{prop}\label{MSShermitian}
 Let $\Delta$ be a Hermitian polar space of rank $r$, defined over the field $\K$, with non-trivial involution $\sigma$ and fixed field $\FF$. Let $U$ be a maximal singular subspace. \begin{itemize}\item[$(i)$] If the anisotropic part is trivial, then $\Pi^+(U)\cong \PSL_r(\K;\FF)$ in its standard permutation representation. Also, $\Pi(U)\cong\PSL_r(\K;\FF)\rtimes \<\rho\>$, with $\rho$ any polarity with companion involution $\sigma$.  \item[$(ii)$] If the anisotropic part is non-trivial, then $\Pi^+(U)\cong\PGL_r(\K)$ in its standard permutation representation. Also, $\Pi(U)\cong \PGL_2(\K)\rtimes\<\rho\>$, with $\rho$ any polarity with companion involution $\sigma$.\end{itemize}
 \end{prop}
 
 \begin{proof}
 First, we note that a single perspectivity induces a duality with companion field automorphism $\sigma$ (as follows from the computations preceding \cref{normset}). Hence, the special projectivity group will be a linear group. Let $\theta$ be any (linear) collineation of $\Delta$ in the little projective group $G^\dagger$ of $\Delta$ stabilising $U$. Since the automorphism group of $\Delta$ acts transitively on opposite (ordered) pairs of maximal singular subspaces,  we may assume that $\theta$ stabilises an opposite maximal singular subspace. We consider the standard Hermitian form 
 \[(x_{-1}^\sigma x_1+\cdots+x_{-r}^\sigma x_r)+(x_1^\sigma x_{-1}+\cdots+x_r^\sigma x_{-r})+f_0(x_0,x_0),\] where $f_0$ is an anisotropic Hermitian form on some  vector space $V_0$ with shorthand coordinates $x_0$. 
We take for $U$ the subspace generated by $e_{-r},\ldots,e_{-1}$. By the above, we may assume that $\theta$ also stabilises the maximal singular subspace $W$, generated by $e_1,\ldots,e_r$. Let $\theta$ be given by \[\begin{pmatrix} x_{-1} & x_{-2}&\ldots&x_{-r}\end{pmatrix}\mapsto \begin{pmatrix} x_{-1} & x_{-2}&\ldots&x_{-r}\end{pmatrix}\cdot M^t,\] with, for now, $M$ any non-singular $r\times r$ matrix over $\K$ (the transpose is taken for notational reasons) with determinant $k\in\K^\times$. It is clear that the action of $\theta$ on $W$ is described by  \[\begin{pmatrix} x_1 & x_2&\ldots&x_r\end{pmatrix}\mapsto \begin{pmatrix} x_1 & x_2&\ldots&x_r\end{pmatrix}\cdot M^{-\theta}.\] 

\begin{itemize}
\item[$(i)$] If $V_0$ is trivial, then the determinant of $\theta$ as a linear map in $\PG(2r-1,\K)$ is equal to $\ell:=kk^{-\theta}$.  Then, according to \cite[Th\'eor\`eme~5]{Die:73}, $\theta$ belongs to $G^\dagger$ if, and only if, $\ell=1$, that is, $k=k^\theta$. Now the statement of $(i)$ follows. 
\item[$(ii)$] Suppose now that $V_0$ is not trivial. Since $\PGL_r(\K)$ is the full group of linear transformations, it suffices to show the assertion in the case of $\dim V_0=1$. Then we may assume $f_0(x_0,x_0)=a_0x_0^\sigma x_0$. We define the action of $\theta$ on the coordinate $x_0$ as $x_0\mapsto \ell^{-1} x_0$. We find that $f_0(\ell x_0,\ell x_0)=a_0\ell^\sigma\ell x_0^\sigma x_0=a_0x_0^\sigma x_0$, since $\ell^\sigma\ell=1$, as one computes.  With that, the determinant of $\theta$ as a linear map in $\PG(2r,\K)$ is equal to $1$, and hence, by  \cite[Th\'eor\`eme~5]{Die:73} again, $\theta$ belongs to $G^\dagger$. This proves the assertions in $(ii)$.\qedhere
\end{itemize} 
 \end{proof}

\subsubsection{Non-embeddable polar spaces}\label{planeNE}
In this case, we could do a computation similar to the one performed in the previous subsection. However, using a result from \cite{Pas-Mal:23}, it is more efficient to use \cref{triangles}. Let $\O$ be a Cayley division algebra over the field $\K$ with standard involution $x\mapsto \overline{x}$, and let $\PG(2,\O)$ be the associated projective plane, coordinatised as in \cref{NEPS}. Then, for all $k,\ell\in\K$, the mapping with action on the points $(x,y)$, with $y\neq 0$, given by 
\[(x,y)\mapsto [k\overline{y}^{-1}\overline{x},\ell\overline{y}^{-1}]\]
induces, by \cite[Section 6.2]{Pas-Mal:23}, a polarity $\rho(k,\ell)$ of $\PG(2,\O)$, which we call a \emph{standard polarity}. 

\begin{remark}\label{remarkstandardpolarity}
Replacing $\O$ with any quadratic associative division algebra $\A$, the above expression remains a polarity in the corresponding projective plane, and we also call such a polarity a \emph{standard polarity} (with respect to the given quadratic structure, that is, relative to the field $\K$, which might not be determined by  $\A$, if $\A$ is commutative). If the standard involution of $\A$ is non-trivial, then the polarity is Hermitian. 
\end{remark}

The following result is now an immediate consequence of the computations in \cite[Section~6.2]{Pas-Mal:23}.

\begin{lem}\label{standardpolarity}
Let $\Delta$ be a thick non-embeddable polar space of rank $3$. Then every perspectivity $\pi_1\per\pi_2\per\pi_3\per\pi_1$, with $\pi_1,\pi_2,\pi_3$ pairwise opposite planes, is a standard polarity. 
\end{lem}

Noting that the class of standard polarities is geometric, we immediately conclude that $\Pi(\pi)$, for $\pi$ a plane of a non-embeddable thick polar space, is the group generated by all standard polarities.   We can be slightly more specific. Recall, that a projective collineation is a collineation, which induces a projectivity on at least one line (and hence on all lines). The set of such collineations forms a group, called the \emph{full projective group}, which we denote by $\PGL_3(\O)$, with slight abuse of notation (because there does not exist a $3$-dimensional vector space over $\O$). 

\begin{prop}
Let $\Delta$ be the thick non-embeddable polar space associated with the Cayley division algebra $\O$ over the field $\K$. Let $\pi$ be any plane of $\Delta$. Then $\Pi^+(\pi)$ is the full projective group of $\PG(2,\O)$, whereas $\Pi(\pi)$ is the full projective group extended with a standard polarity.  
\end{prop}

\begin{proof}
By \cite[Section~6.2]{Pas-Mal:23}, the action of $\rho(k,\ell)$ on lines $[m,k]$, with $k\neq 0$, is given by \[[m,k]\mapsto (k^{-1}\ell \overline{m} \overline{k}^{-1},\ell\overline{k}^{-1}).\]
We can calculate the action of $\rho(k,\ell)\rho(1,1)$ on the points $(x,y)$, $y\neq 0$ and obtain \[(x,y)\mapsto(k\ell^{-1}x,\ell^{-1}y).\]
Setting $\ell=1$, we see that we obtain all homologies with axis $[0]$ and centre $(0)$. Then the assertion follows from \cite[Satz~3]{Pic:59}.  
\end{proof}


\section{Projectivity groups of points of metasymplectic spaces}\label{pointsMSS}

For the first proposition of this section, one might benefit from having a clear visualisation of an apartment of a building of type $\mathsf{F_4}$. We provide one below. It is a graph, for which its vertices correspond to the vertices of type $1$ (with the usual Bourbaki labelling) of the said apartment, the edges are those of type $2$ , the $3$-cliques are those of type $3$, and the (skeletons of the) octahedra are the vertices of type $4$. 


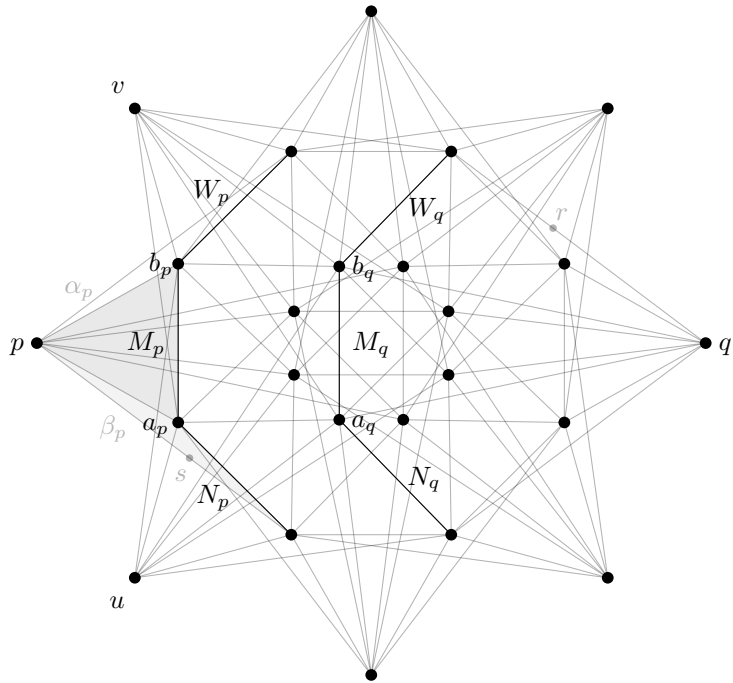
\begin{center}
\begin{tikzpicture}[scale=2.2]

    \def\LineOpacity{0.3}

    \draw[thin, transparent] (0:2) 
        -- (45:2) 
        -- (90:2) 
        -- (135:2) 
        -- (180:2) 
        -- (225:2) 
        -- (270:2) 
        -- (315:2) 
        -- cycle; 

    \foreach \angle in {0, 45, 90, 135, 180, 225, 270, 315} 
        \fill (\angle:2) circle (1pt); 
    
    \draw[thin, opacity=\LineOpacity, rotate around={22.5:(0,0)}] (0:1.25) 
        -- (45:1.25) 
        -- (90:1.25) 
        -- (135:1.25) 
        -- (180:1.25) 
        -- (225:1.25) 
        -- (270:1.25) 
        -- (315:1.25) 
        -- cycle; 

    \foreach \angle in {0, 45, 90, 135, 180, 225, 270, 315} 
        \fill[black] (\angle+22.5:1.25) circle (1pt); 
    
    \draw[thin, transparent, rotate around={22.5:(0,0)}] (0:0.5) 
        -- (45:0.5) 
        -- (90:0.5) 
        -- (135:0.5) 
        -- (180:0.5) 
        -- (225:0.5) 
        -- (270:0.5) 
        -- (315:0.5) 
        -- cycle; 

    \foreach \angle in {0, 45, 90, 135, 180, 225, 270, 315} 
        \fill[black] (\angle+22.5:0.5) circle (1pt); 

    \foreach \angle in {0, 45, 90, 135, 180, 225, 270, 315} {
        \draw[thin, black, opacity=\LineOpacity] (\angle+22.5:0.5) -- (\angle+22.5+135:0.5); 
    }

    \foreach \angle in {0, 45, 90, 135, 180, 225, 270, 315} {
        \draw[thin, black, opacity=\LineOpacity] (\angle:2) -- (\angle+22.5:1.25); 
        \draw[thin, black, opacity=\LineOpacity] (\angle:2) -- (\angle-22.5:1.25); 
    }

    \foreach \angle in {0, 45, 90, 135, 180, 225, 270, 315} {
        \draw[thin, black, opacity=\LineOpacity] (\angle:2) -- (\angle+67.5:1.25); 
        \draw[thin, black, opacity=\LineOpacity] (\angle:2) -- (\angle-67.5:1.25); 
    }

    \foreach \angle in {0, 45, 90, 135, 180, 225, 270, 315} {
        \draw[thin, black, opacity=\LineOpacity] (\angle+22.5:1.25) -- (\angle+22.5+45:0.5); 
        \draw[thin, black, opacity=\LineOpacity] (\angle+22.5:1.25) -- (\angle+22.5-45:0.5); 
    }

    \foreach \angle in {0, 45, 90, 135, 180, 225, 270, 315} {
        \draw[thin, black, opacity=\LineOpacity] (\angle:2) -- (\angle+22.5:0.5); 
        \draw[thin, black, opacity=\LineOpacity] (\angle:2) -- (\angle-22.5:0.5); 
    }

    \foreach \angle in {0, 45, 90, 135, 180, 225, 270, 315} {
        \draw[thin, black, opacity=\LineOpacity] (\angle:2) -- (\angle+22.5+90:0.5); 
        \draw[thin, black, opacity=\LineOpacity] (\angle:2) -- (\angle-22.5-90:0.5); 
    }


    \coordinate (p) at (180:2); 
    \node at (p) [behind path, fill opacity = 0, fill=white, label={[label distance=-1mm] left:{\small \( p \)}}] {};

    \coordinate (q) at (0:2); 
    \node at (q) [behind path, fill opacity = 0, fill=white, label={[label distance=-1mm] right:{\small \( q \)}}] {};
    
    \coordinate (v) at (135:2); 
    \node at (v) [behind path, fill opacity = 0, fill=white, label={[label distance=-2mm] above left:{\small \( v \)}}] {};

    \coordinate (u) at (225:2); 
    \node at (u) [behind path, fill opacity = 0, fill=white, label={[label distance=-2mm]below left:{\small \( u \)}}] {};

   \begin{scope}[rotate around={22.5:(0,0)}] 
    \coordinate (q') at (45:1.25); 
    \coordinate (r) at ($ (q) + 0.6*(q') - 0.6*(q) $); 
    \node at (r) [opacity=0.3, circle, fill, inner sep=1pt, label={[label distance=-2mm, opacity=0.3] above right: {\small \( r \)}}] {};
    
    \coordinate (p') at (225:1.25); 
    \coordinate (s) at ($ (p) + 0.6*(p') - 0.6*(p) $); 
    \node at (s) [opacity=0.3, circle, fill, inner sep=1pt, label={[label distance=-2mm, opacity=0.3] below left: {\small \( s \)}}] {};

    \coordinate (ap) at (180:1.25); 
    \node at (ap) [label={[label distance=-1.5mm] left: {\small \( a_p \)}}] {};

    \coordinate (bp) at (135:1.25); 
    \node at (bp) [label={[label distance=-2mm] left: {\small \( b_p \)}}] {};

    \coordinate (aq) at (225:0.5); 
    \node at (aq) [label={[label distance=-1.1mm] right: {\small \( a_q \)}}] {};
    
    \coordinate (bq) at (90:0.5); 
    \node at (bq) [label={[label distance=-1mm] right: {\small \( b_q \)}}] {};

    \coordinate (ap') at (225:1.25); 
    \draw[black] (ap) -- (ap');
    \coordinate (Np) at ($ (ap) + 0.5*(ap') - 0.5*(ap) $); 
    \node at (Np) [label={[label distance=-3mm] below left: {\small \( N_p \)}}] {};
    
    \coordinate (aq') at (270:1.25); 
    \draw[black] (aq) -- (aq');
    \coordinate (Nq) at ($ (aq) + 0.5*(aq') - 0.5*(aq) $); 
    \node at (Nq) [label={[label distance=-1.1mm] right: {\small \( N_q \)}}] {};    
    
    \coordinate (bp') at (90:1.25); 
    \draw[black] (bp) -- (bp');
    \coordinate (Wp) at ($ (bp) + 0.5*(bp') - 0.5*(bp) $); 
    \node at (Wp) [label={[label distance=-3mm] above left: {\small \( W_p \)}}] {};

    \coordinate (bq') at (45:1.25); 
    \draw[black] (bq) -- (bq');
    \coordinate (Wq) at ($ (bq) + 0.5*(bq') - 0.5*(bq) $); 
    \node at (Wq) [label={[label distance=-1mm] right: {\small \( W_q \)}}] {};

    \draw[black] (ap) -- (bp);
    \coordinate (Mp) at ($ (ap) + 0.5*(bp) - 0.5*(ap) $); 
    \node at (Mp) [label={[label distance=-1mm] left: {\small \( M_p \)}}] {};              
    
    \draw[black] (aq) -- (bq);
    \coordinate (Mq) at ($ (aq) + 0.5*(bq) - 0.5*(aq) $); 
    \node at (Mq) [label={[label distance=-1mm] right: {\small \( M_q \)}}] {};       
    \end{scope}
    
    \fill[black, opacity=0.09] (p) -- (bp) -- (ap) -- cycle;    
    \fill[black, opacity=0.05] (p) -- (ap) -- (ap') -- cycle;   
    \coordinate (alphap) at ($ (p) + 0.3*(bp) - 0.3*(p) $); 
    \node at (alphap) [label={[opacity=0.3, label distance=0mm] above: {\small \( \alpha_p \)}}] {};         
    \coordinate (betap) at ($ (p) + 0.3*(ap') - 0.3*(p) $); 
    \node at (betap) [label={[opacity=0.3, label distance=-1mm] below: {\small \( \beta_p \)}}] {};      
\end{tikzpicture}
\end{center}


As in \cite{Bus-Sch-Mal:24}, we call a collineation of a polar space  a \emph{homology}, if it pointwise fixes two opposite planes, which we call the \emph{axes} of the homology.

\begin{prop}\label{F4type1}
Let $p$ be a point of a metasymplectic space $\Gamma$. Let $q$ be any opposite point, and let $s\perp p$ and $r\perp q$ be further points such that $ps\equiv qr$, $p\equiv r$ and $q\equiv s$. Then the projectivity $\rho_{r,s}:p \per q \per s \per r \per p$ is a homology of the polar space $\Delta_p$, corresponding to the residue at $p$, which pointwise fixed planes $\alpha_s$ and $\alpha_r$, the planes corresponding to the lines $sp$ and the projection of $rq$ onto $p$, respectively. Moreover, if $L$ is any line of $\Delta_p$ intersecting $\alpha_r$ and $\alpha_s$ non-trivially, then $r$ and $s$ can be chosen, such that the projectivity   $\rho_{r,s}$ maps a given plane of $\Delta$ through $L$, not adjacent to either $\alpha_r$ or $\alpha_s$, to any other such plane.  
\end{prop}

\begin{proof}
The line $L$ of $\Delta_p$ corresponds to a plane $\alpha$ through $p$ in $\Gamma$. Let $\Sigma$ be an apartment in $\Gamma$ containing the opposite lines $ps$ and $qr$, containing the points $p$ and $q$, and containing the plane $\beta$ (this will only be crucial in the last part of the proof; hence, in the first part, we will not use this information). Let $\alpha_p$ and $\beta_p$ be two planes in $\Sigma$ through $p$ that intersect in a line, and let $s\in \beta_p$. Also, for the last part of the proof (Claim 3), we choose these planes such that $\alpha_p=\alpha$.  Let $M_p$ and $N_p$ be the lines in $\alpha_p$ and $\beta_p$, respectively, that are closest to $q$ (meaning the points of these lines have distance two to $q$ and are special to $q$). We denote by $a_p$ the intersection point of $M_p$ and $N_p$. Let $M_q$ and $N_q$ be the lines in $q^{\perp} \cap \Sigma$, such that each point of $M_q$ ($N_q$, respectively) is collinear to a unique point of $M_p$ ($N_p$). Let $u$ be the point of $\Sigma$ that is collinear to all points of $M_p$, $N_p$, $M_q$ and $N_q$. Let $W_p$ be the line of $\Sigma$ that intersects $\alpha_p$ in exactly one point, whose points are all collinear to $p$ and at distance $2$ to $q$ and which does not contain a point collinear to $N_p$. Let $b_p$ be the intersection point of $M_p$ and $W_p$. Again, let $W_q$ be the line in $q^{\perp} \cap \Sigma$, such that each point of $W_q$ is collinear to a unique point of $W_p$. Consider the plane $\langle M_p, u \rangle$ and let $S\subseteq\<M_p,u\>$ be the line through $a_p$ collinear to $s$ (this line exists considering the symp containing $p$ and $u$). We note that the plane $\langle S, s \rangle$ intersects $\alpha_p$ exactly in the point $a_p$, is not contained in $p^{\perp}$ and since $s$ has distance $3$ to $q$, we have $u\notin S$.

Let $\xi$ be the octahedron of $\Sigma$ containing $M_p$, $M_q$ and $u$. Let $v$ be the point in $\xi$ that is opposite $u$ and collinear to every point of $M_p$ and $M_q$. Note that $v$ is also collinear to $W_p$ and $W_q$. Hence, $v$ and $q$ are contained in a common symp. 

Let $a_q$ be the point $M_q \cap N_q$; that is the unique point of $M_q$ and $N_q$ that is collinear to $a_p$. Then $\langle a_p, a_q, v \rangle$ spans a plane of $\Sigma$.  
We denote the point $M_q \cap W_q$ by $b_q$. 

\textbf{Claim 1.} \emph{We claim that the projectivity $p \per q \per s \per r \per p$ stabilises the planes $\<p,N_p\>$ and $\<p,W_p\>$ and fixes the lines $pa_p$ and $pb_p$.}

Indeed, the plane $\langle p, N_p \rangle$ maps to $\langle q, N_q \rangle$ under the first projection and back to $\langle p, N_p \rangle$ under the second projection. It then maps to some plane through $r$, but since that plane is the image of $\langle p, N_p \rangle$, it maps back to $\langle p, N_p \rangle$  under the fourth projection. Similarly, $\langle p, W_p \rangle$ projects to $\langle q, W_q \rangle$, to some plane through $s$ and then back to $\langle q, W_q \rangle$ and $\langle p, W_p \rangle$.

The line $pa_p$ maps to $qa_q$, $sa_p$, the line between $r$ and $R \cap va_q$ and back to $pa_p$. 
The line $pb_p$ maps to the line between $q$ and $W_q \cap M_q$, maps to the line between $s$ and $S \cap ub_p$, maps to the line between $r$ and $W_q \cap M_q$ and back to $pb_p$.

\textbf{Claim 2.} \emph{We claim that the projectivity $p \per q \per s \per r \per p$ stabilises all planes through $ps$.}

This follows from Step 1 by varying $\Sigma$ through $\{p, ps\}$ and $\{q, qr\}$. 

\textbf{Claim 3.} \emph{We claim that we can always define $s$ and $r$ in a way, such that the projectivity $p \per q \per s \per r \per p$ moves a line $px$ with $x$ on $M_p$ to a line $px'$ with $x'$ on $M_p$, $x \neq x'$ and $x,x' \notin \{ a_p, b_p \}$. }

Let $x$ and $x'$ be two distinct points on $M_p$ not equal to $a_p$ or $b_p$. Let $x_{1}$ be the unique point of $M_q$ collinear to $x$. Then $px$ moves to $qx_{1}$ under the first projection. 

Considering the symp spanned by $\xi$, we can see that $x_1$ must be collinear to a unique point $x_2$ of $S$. With that, $qx_1$ has to move to $sx_2$ under the second projection. The points $x_2$ and $x'$ are both each collinear to a line of $\langle v, a_q, b_q \rangle$. These lines are distinct as $x_1$ belongs to one of them but not the other. We define $x_3$ as the intersection point of these two lines.  
Let $R$ be the connection line between $b_q$ and $x_3$. Every point of $R \setminus \{b_q\}$ is collinear to the same line of $\langle q, W_q \rangle$ through $b_q$ that we will denote by $R'$. In $\Sigma$ we see that $W_q$ and $qb_q$ form a triangle and we denote the third point of that triangle in $\Sigma$ by $c_q$. We define $r$ to be the intersection point of $R'$ and $qc_q$.  Then $sx_2$ moves to $rx_3$ under the third projection and to $px'$ under the last projection as desired.
\end{proof}

We now first handle the case of type $1$ vertices in  $\mathsf{F_4}(\K,\A)$, with $\A$ a quadratic alternative division algebra over $\K$. We may assume that $|\A|>2$ as otherwise the groups of projectivity coincide with the full group of collineations $\mathsf{Sp}_6(2)$ by \cref{MR0}, since that group is also generated by the elations. 

We denote by $\mathsf{PSp}_{2\ell}(\K)$ the group of all collineations of $\mathsf{C_{\ell,1}}(\K,\K)$ preserving the associated alternating form, $\ell\geq 2$. It is a simple group. We denote by $\overline{\mathsf{PS}}\mathsf{p}_{2\ell}(\K)$ the group generated by  $\mathsf{PSp}_{2\ell}(\K)$ and the \emph{diagonal} collineations, that is, the linear collineations of the underlying vector space mapping the associated alternating form to a nonzero scalar multiple and represented by a diagonal matrix. 

If $\kar\K=2$, and $\K'$ is an overfield of $\K$ all of whose squares are contained in $\K$, then the polar space $\mathsf{C}_{\ell,1}(\K',\K)$ is a polar subspace of $\mathsf{C}_{\ell,1}(\K',\K')$, and we denote the restriction of $\overline{\mathsf{PS}}\mathsf{p}_{2\ell}(\K')$ to $\mathsf{C}_{\ell,1}(\K',\K)$ by $\overline{\mathsf{PS}}\mathsf{p}_{2\ell}(\K',\K)$. 

Likewise, for $\A$ a separable quadratic extension of $\K$, when $\mathsf{U}_6(\A/\K)$ is the (simple) unitary group preserving a Hermitian form and whose elements correspond to matrices of determinant $1$, then we denote by $\overline{\mathsf{U}}_6(\A/\K)$ the group generated by $\mathsf{U}_6(\A/\K)$ and all diagonal automorphisms with diagonal elements in $\K$ (with respect to the standard form). 

Now let $\A$ be a quaternion division algebra over $\K$. Let $\mathsf{U}_6(\A)$ denote the (simple) collineation group of $\mathsf{C_{3,1}}(\A,\K)$ generated by the elations, then the group generated by all elations and diagonal automorphisms (with factors in $\K$ and with respect to the standard form) is denoted by $\overline{\mathsf{U}}_6(\A)$. Note that we do not need to mention $\K$ in the notation, since it is unique as the centre of $\A$. 

Finally, let $\A$ be a Cayley division algebra over $\K$. Here, there is no form of the corresponding polar space available (since it is non-embeddable) and hence we cannot consider diagonal automorphisms as in the previous paragraphs.   
However, we can either consider the group generated by all elations and homologies (which for the previous cases would have boiled down to the same groups), or all elations and the diagonal automorphisms of the universal embedding of the corresponding dual polar space, see the proof of \cref{sharply}. We denote these automorphism groups by $\mathsf{E}_{7,3}^{(28)}(\A)$ and $\overline{\mathsf{E}}_{7,3}^{(28)}(\A)$, respectively.

\begin{lem}\label{sharply}
Let $\Delta$ be a polar space of rank $3$ isomorphic to $\mathsf{C_{3,3}}(\A,\K)$, with $\A$ a quadratic alternative division algebra over $\K$,  
and let $\alpha$ and $\beta$ be two disjoint planes. Let $L$ be a line intersecting both $\alpha$ and $\beta$ non-trivially. Then there exists a unique homology with axes $\alpha$ and $\beta$ and mapping a given plane $\pi$ though $L$, not adjacent to either $\alpha$ or $\beta$, to an arbitrary other given plane like that. In particular, there is a unique set of homologies with axes $\alpha$ and $\beta$, acting transitively on the set of planes through $L$, not adjacent with either $\alpha$ or $\beta$, and hence, such a set is geometric.
\end{lem}

\begin{proof}
This is best seen through the (universal) representation of the corresponding dual polar space in a projective space over $\K$, as established uniformly in  \cite{Bru-Mal:16}. We take the slightly more symmetric algebraic description from \cite[Definition~10.1]{Sch-Sch-Mal-Vic:23}. Let $V$ be the vector space $\K^4\oplus\A^3\oplus\K^3\oplus\A^3\oplus\K$ (note that we view $\A$ as a vector space over $\K$ in the natural way, that is, coming from the algebra over $\K$). Then the projective (or Zariski) closure of the point set given by the following parameter form, with the induced line set, forms the dual polar space $\mathsf{C_{3,3}}(\A,\K)$. We denote by $x\mapsto\overline{x}$ the standard involution in $\A$. 
\[(1,\ell_1,\ell_2,\ell_3,x_1,x_2.x_3,x_1\overline{x}_1-\ell_2\ell_3,x_1\overline{x}_1-\ell_2\ell_3,x_1\overline{x}_1-\ell_2\ell_3,\] \[\ell_1\overline{x}_1-x_2x_3,\ell_1\overline{x}_1-x_2x_3,\ell_1\overline{x}_1-x_2x_3,\ell_1x\overline{x}_1+\ell_2x_2\overline{x}_2+\ell_3x_3\overline{x}_3-x_1(x_2x_3)-\overline{x}_3(\overline{x}_2\overline{x}_1)-\ell_1\ell_2\ell_3),\]

where $\ell_1,\ell_2,\ell_3\in\K$ and $x_1,x_2,x_3\in\A$. We let $\alpha$ and $\beta$ correspond to the points $a=(1,0,0,\ldots,0)$ and $b=(0,0,\ldots,0,1)$, respectively, We write coordinates of points of $\PG(V)$ as $14$-tuples with entries in $\K\cup\A$ in the natural way, according to the definition of $V$ above. Then one sees that the set of points $a^\perp\cap b^{\not\equiv}$ spans the subspace $(0,*,*,*,*,*,*,0,\ldots,0)$ and $a^{\not\equiv}\cap b^\perp$ spans the subspace $(0,\ldots,0,*,*,*,*,*,*,0)$. Each homology with axes $\alpha$ and $\beta$ has to fix all points of these subspaces, and one can see that such a homology stems from the following linear map, where $t\in\K^\times$:
\begin{align*}
&V \to V \\
&(k,\ell_1,\ell_2,\ell_3,x_1,x_2,x_3,m_1,m_2,m_3,y_1,y_2,y_3,k') \\ &\mapsto  (k,t\ell_1,t\ell_2,t\ell_3,tx_1,tx_2,tx_3,t^2m_1,t^2m_2,t^2m_3,t^2y_1,t^2y_2,t^2y_3,t^3k')
\end{align*}
According to \cite{Bru-Mal:16}, the line generated by $a'=(0,0,1,0,0,\ldots,0)$ and $b'=(0,\ldots,0,1,0,0,0,0,0,0)$ corresponds to a line joining a point of $\alpha$ with one of $\beta$ (in the coordinates of \cite{Bru-Mal:16} these points are $(0,0,0,0,0)$ and $(0,0)$, and the line has type \textbf{(I)}). The point $a'$ is collinear to $a$ and $b'$ to $b$. The above linear mappings, for varying $t\in\K^\times$, act (sharply) transitively on the points of $a'b'\setminus\{a',b'\}$. 

This proves the lemma.
\end{proof}

\begin{coro}\label{type1F4}
Let $v$ be a vertex of type $1$ in $\mathsf{F_4}(\K,\A)$, with $\A$ a quadratic alternative division algebra over $\K$, $|\A|>2$. Then $\Pi^+(\{p\})$ is generated by all homologies of the corresponding polar space. 
\end{coro}

\begin{proof}
This follows from \cref{UpanddownF4}, \cref{sharply} and \cref{F4type1}. 
\end{proof}

In view of \cref{special=generalMSS}, this proves the case (C3) of \cref{tableF4}

We now consider vertices of type $4$ in $\mathsf{F_4}(\K,\A)$. Since $\mathsf{F_{4,1}}(\K,\K')\cong \mathsf{F_{4,4}}(\K'^2,\K)$, if $\kar\K=2$ and $\K'$ is an inseparable extension of $\K$, we may restrict ourselves to the case, where $\A$ is separable; that is, either equal to $\K$ with $\kar\K\neq 2$, or a separable quadratic extension of $\K$ (in any characteristic), or a quaternion or octonion division algebra over $\K$. 

First we consider the case $\A=\K$, with $\kar\K\neq 2$. Let $p$ be a point of $\mathsf{F_{4,4}}(\K,\K)$. We have the following analogue of \cref{sharply}.

\begin{lem}\label{sharply2}
Let $\Delta$ be a parabolic polar space of rank $3$ over the field $\K$ with $\kar\K\neq 2$, 
and let $\alpha$ and $\beta$ be two disjoint planes. Let $L$ be a line intersecting both $\alpha$ and $\beta$ non-trivially. Then there exists a unique homology with axes $\alpha$ and $\beta$ mapping a given plane $\pi$ though $L$, not adjacent to either $\alpha$ or $\beta$, to an arbitrary other given plane like that. In particular, there is a unique set of homologies with axes $\alpha$ and $\beta$ acting transitively on the set of planes through $L$, not adjacent with either $\alpha$ or $\beta$, and hence such set is geometric. Also, there is a unique homology with axes $\alpha$ and $\beta$ that is a reflection. 
\end{lem}

\begin{proof}
We consider the standard equation of a parabolic quadric in $\PG(6,\K)$, that is,
\[X_{-3}X_3+X_{-2}X_2+X_{-1}X_1=X_0^2.\]
We write the coordinates as $(X_{-3},X_{-2},X_{-1},X_0,X_1,X_2,X_3)$. We can take, with self-explaining notation, \[\begin{cases} \alpha=(*,*,*,0,0,0,0), \\ \beta=(0,0,0,0,*,*,*),\\ L=\<(1,0,0,0,0,0,0),(0,0,0,0,0,1,0)\>.\end{cases}\] Then  $\pi_k=\<L,(0,0,k^{-1},1,k,0,0)\>$, $k\in\K^\times$, is a generic plane through $L$ in $\Delta$ not adjacent to either $\alpha$ or $\beta$. With that, the linear map
\[(x_{-3},x_{-2},x_{-1},x_0,x_1,x_2,x_3)\mapsto (x_{-3},x_{-2},x_{-1},kx_0,k^2x_1,k^2x_2,k^2x_3)\] maps $\pi_1$ to $\pi_k$. This is a reflection if, and only if, it fixes each point with $x_0=0$, that is, if, and only if, $k^2=1$, or $k\in\{1,-1\}$. 

Now suppose a homology $\varphi$ with axes $\alpha$ and $\beta$ also stabilises the plane $\pi_1$. Then $\varphi$ stabilises each line of $\pi_1$ through $\pi_1\cap\alpha$ (as such a line is the projection onto $\pi_1$ of some fixed point of $\alpha$), and also each line of $\pi_1$ through $\pi_1\cap\beta$ (for the analogous reason).  Hence, $\pi_1$ is fixed pointwise, and this readily implies that $\PG(6,\K)$ is fixed pointwise. 

All assertions are proved. 
\end{proof}

We denote by $\mathsf{O}_7(\K)$ the group of linear collineations of a $7$-dimensional vector space $V$ over $\K$, preserving the bilinear form associated to the standard equation of a parabolic quadric in $\PG(6,\K)$, and by $\mathsf{PO}_7(\K)$ the quotient with its centre. Note that the former coincides with  $\mathsf{O^+_7}(\K)$ with the notation of \cite{Die:73}, as $\dim V=7$ is odd.  We then have:

\begin{lem}\label{type4F4}
Let $p$ be a point of $\mathsf{F_{4,4}}(\K,\K)$. Then $\Pi^+(p)=\Pi(p)\cong\mathsf{PO}_7(\K)$ in its standard action. 
\end{lem}

\begin{proof} 
By  \cref{UpanddownF4}, \cref{F4type1} and \cref{sharply2}, $\Pi^+(p)$ is generated by all homologies, hence, it is contained in $\mathsf{PO}_7(\K)$. 
Now let $q$ be opposite $p$ and consider the extended equator geometry $\widehat{E}(p,q)$. If we now restrict to perspectivities $x\per y$ with $x$ and $y$ points of $\widehat{E}(p,q)$, then we see that the special projectivity group of $p$ in $\mathsf{B_{4,1}}(\K,\K)\cong \widehat{E}(p,q)$ is a subgroup of the special projectivity group $\Pi^+(p)$ inside $\Delta$. Hence \cref{herm-sympl} yields $\mathsf{PO}_7(\K)\leq \Pi^+(p)$. Also, using  \cref{pointsgeneral}, we find that there exists a self-projectivity of $\Delta$ of length $3$, exhibiting a reflection of the residue in $p$. Hence, $\Pi^+(p)=\Pi(p)=\mathsf{PO}_7(\K)$ and the proof is complete.
\end{proof}

This takes care of the first line of (B3) of \cref{tableF4} (and note that for $\kar\K=2$, the groups $\mathsf{PO}_7(\K)$ and $\overline{\mathsf{PS}}\mathsf{p}_6(\K)$ coincide, so we do not have to require that $\kar\K\neq 2$ on that line in the table).

Now let $\A$ be a separable quadratic alternative division algebra over the field $\K$, with $d:=\dim_{\K}\A\in\{2,4,8\}$. Let $\Delta$ be the polar space of rank $3$, obtained from the quadric in $\PG(5+d,\K)$ with standard equation
\[X_{-3}X_3+X_{-2}X_2+X_{-1}X_1=\Norm(X_0),\]
where the coordinate $X_0$ belongs to $\A$, viewed as vector space over $\K$ in the natural way, and where $\Norm(X_0)=X_0\overline{X}_0$, with $X_0\mapsto\overline{X}_0$ the standard involution in $\A$ with respect to the quadratic algebra structure. Let $V$ be the vector space $K^3\oplus\A\oplus\K^3$. Then we denote by $\mathsf{O}_{d+6}(\K,\A)$ or $\overline{\mathsf{O}}_{d+6}(\K,\A)$ the group of all linear transformations of  $V$ preserving the quadratic form associated with the above equation, or mapping it to a scalar multiple (the similitudes; the corresponding scalar is called the \emph{factor} of the similitude), respectively,  and we denote with $\mathsf{PO}_{d+6}(\K,\A)$ and $\overline{\mathsf{PO}}_{d+6}(\K,\A)$ the respective quotients with the centre. Let $\FF$ be a splitting field of $\A$, that is, a quadratic extension of $\K$ over which $\A$ splits as an algebra. Then there are two natural systems of maximal singular subspaces (corresponding to those of the associated hyperbolic quadric) of the corresponding quadratic form over $\FF$. The subgroups of $\mathsf{O}_{d+6}(\K,\A)$ and  $\overline{\mathsf{O}}_{d+6}(\K,\A)$, preserving each of these systems, will be denoted by $\mathsf{O}^+_{d+6}(\K,\A)$ and $\overline{\mathsf{O}}^+_{d+6}(\K,\A)$, respectively, partially following Dieudonn\'e \cite{Die:73} (instead of the bilinear form, we included the algebra in the notation). The corresponding projective groups are then $\mathsf{PO}^+_{d+6}(\K,\A)$ and $\overline{\mathsf{PO}}^+_{d+6}(\K,\A)$.

Furthermore, let $\alpha$ and $\beta$ be two opposite planes of $\Delta$. Then we say that the homology group with axes $\alpha$ and $\beta$ \emph{acts transitively}, or \emph{is a transitive homology group} if for some line $L$ intersecting both $\alpha$ and $\beta$ non-trivially (and then for each such line), and each pair of planes $\pi,\pi'$ through $L$ but not intersecting either $\alpha$ or $\beta$ in a line, there exists a homology with axes $\alpha$ and $\beta$ mapping $\pi$ to $\pi'$. 

Note that the factors of the similitudes are precisely the non-zero norms of $\A$. 

We now prepare for the determination of the projectivity groups of a point in $\mathsf{F_{4,4}}(\K,\A)$, with $\A$ separable and of dimension $2,4$ or $8$ over $\K$.

\begin{lem}\label{sep248}
Let $G$ be a group of collineations of $\mathsf{B_{3,1}}(\A,\K)$, with  $\A$ a separable quadratic alternative division algebra over the field $\K$, with $d=\dim_{\K}\A\in\{2,4,8\}$, containing  $\mathsf{PO}^+_{d+6}(\K,\A)$ and containing a transitive homology group. Then $G$ contains $\overline{\mathsf{PO}}^+_{d+6}(\K,\A)$. 
\end{lem}

\begin{proof}
It suffices to show that, for each norm $r\in\K$, there exists a member of $G$ mapping a given quadratic form describing $\mathsf{B_{3,1}}(\A,\K)$ to an $r$-multiple. We consider, with previous notation, the quadratic form  \[X_{-3}X_3+X_{-2}X_2+X_{-1}X_1-\Norm(X_0).\] As in the proof of \cref{sharply2}, we set  \[\begin{cases} \alpha=(*,*,*,0,0,0,0), \\ \beta=(0,0,0,0,*,*,*),\\ L=\<(1,0,0,0,0,0,0),(0,0,0,0,0,1,0)\>.\end{cases}\] Since $r$ is a norm, there exists some $a\in\A$ with $\Norm(a)=r$. Consider the plane $\pi_a$ spanned by $L$ and the point $p_a=(0,0,1,a,r,0,0)$.  Let $\pi_1$ be the plane spanned by $L$ and $p_1=(0,0,1,1,1,0,0)$. Then there exists some homology $h$ with axes $\alpha$ and $\beta$ mapping $\pi_1$ to $\pi_a$. Now $p_1$ and $p_a$ are the unique points of $\pi_1$ and $\pi_a$, respectively, collinear to the fixed points $(0,1,0,0,0,0,0)$ and $(0,0,0,0,0,0,1)$. Hence, the point $p_1$ is mapped by $h$ to $p_a$. Since all points of $\alpha\cup\beta$ are fixed, the matrix of $h$ is diagonal, and hence, this diagonal is, up to a scalar multiple, equal to $(1,1,1,M,r,r,r)$, where $M$ is a $d\times d$ matrix with as first column $a$ (in its coordinates over $\K$; the action is on the right). It follows that $h$ is a similitude with factor $r$ and the lemma is proved.  
\end{proof}

\begin{lem}\label{type4F4bis}
Let $p$ be a point of $\mathsf{F_{4,4}}(\K,\A)$, with $\A$ separable and of dimension $2,4$ or $8$ over $\K$. Then $\Pi^+(p)\cong\overline{\mathsf{PO}}^+_{d+6}(\K,\A)$ and $\Pi(p)\cong \overline{\mathsf{PO}}_{d+6}(\K,\A)$. 
\end{lem}

\begin{proof}
This time, the set of homologies obtained from \cref{F4type1} is not (necessarily) geometric. However, \cref{UpanddownF4} still shows that $\Pi^+(p)$ is generated by a set of homologies, which contains a transitive group of homologies. It already shows that $\Pi^+(p)\leq \overline{\mathsf{PO}}_{d+6}(\K,\A)$, such that $\Pi^+(p)$ contains a similitude with an arbitrary norm as factor. 

However, since splitting the forms over $\FF$ (see above) produces the exceptional buildings of types $\mathsf{E_6,E_7,E_8}$ for $d=2,4,8$, respectively, and in these buildings, an odd projectivity \emph{always} switches the natural classes of maximal singular subspaces of the hyperbolic quadratic associated to the splitting of the symps  (as can be read off of Table~2 of \cite{Bus-Sch-Mal:24} in the lines labelled (D4), (D5) and (D7)), we deduce that $\Pi^+(p)\leq \overline{\mathsf{PO}}^+_{d+6}(\K,\A)$.

Now we select a point $q$ opposite $p$ and restrict the projectivity group $\Pi(p)$ to self-projectivities using only  perspectivities between points of the extended equator geometry $\widehat{E}(p,q)$. Then \cref{pointsgeneral} implies that every reflection is contained in $\Pi(p)$, seen as a polar space. Hence, $\mathsf{PO}_{d+6}(\K,\A)\leq\Pi(p)$ and the previous paragraphs and \cref{sep248} imply now that $\Pi(p)$ contains $\overline{\mathsf{PO}}_{d+6}(\K,\A)$. But since $\Pi^+(p)\leq \overline{\mathsf{PO}}^+_{d+6}(\K,\A)$, because $\Pi^+(p)$ has index at most $2$ in $\Pi(p)$ and $\overline{\mathsf{PO}}^+_{d+6}(\K,\A)$ has exactly index $2$ in $\overline{\mathsf{PO}}_{d+6}(\K,\A)$, all projectivity groups now follow.   
\end{proof}


\section{Projectivity groups of non-maximal residues in metasymplectic spaces}\label{nonmax}

We start with the residues isomorphic to generalised quadrangles, that is, buildings of type $\mathsf{B_2}$ (or $\mathsf{C_2}$). The groups $\mathsf{PO}_5(\K)$ and $\overline{\mathsf{PO}}^+_{d+4}(\K,\A)$, with $\A$ a separable quadratic alternative division algebra over $\K$ of dimension $2,4$ or $8$, are similarly defined as their higher dimensional analogues above. We then have the following lemma.

\begin{lem}\label{14F4}
Let $F$ be a simplex of type $\{1,4\}$ of a building $\Delta:=\mathsf{F_4}(\K,\A)$. Then $\Pi^+(F)=\Pi(F)$. If $\A=\K'$ is an inseparable extension of $\K$, when $\kar\K=2$, then $\Pi(F)\cong \overline{\mathsf{PS}}\mathsf{p}_{4}(\K',\K)$; this includes $\K=\K'$ for $\kar\K=2$. If $\K=\A$, then $\Pi(F)\cong \mathsf{PO_5}(\K)\cong\overline{\mathsf{PS}}\mathsf{p}(\K)$. If $\A$ is separable with $\dim_{\K}\A\in\{2,4,8\}$, then $\Pi(F)\cong \overline{\mathsf{PO}}^+_{d+4}(\K,\A)$. 
\end{lem}

\begin{proof}
The fact that $\Pi^+(F)=\Pi(F)$ follows from \cite[Lemma~5.2]{Bus-Sch-Mal:24}. More concretely, let $F_1=\{x_1,\xi_1\}$ be an incident point-symp pair of $\mathsf{F_{4,1}}(\K,\A)$, and let $F_2=\{x_2,\xi_2\}$ be an opposite point-symp pair. Let $x_3\notin\{x_1,x_2\}$ be a point on the imaginary line defined by $x_1$ and $x_2$, that is, $x_3\in \{x_1,x_2\}^{\pperp\pperp}\setminus\{x_1,x_2\}$. Finally, let $\xi_3$ be a symp through $x_3$, intersecting $E(x_1,x_2)$ in a point of the hyperbolic line of $E(x_1,x_2)$, defined by $E(x_1,x_2)\cap\xi_i$, $i=1,2$. Since $E(x_1,x_2)$ is either a symplectic, a mixed, a unitary, or a thick non-embeddable polar space, we can choose $\xi_3$ disjoint from $\xi_1\cup\xi_2$, and it follows that $\xi_3$ is opposite both $\xi_1$ and $\xi_2$. Now one can check  that $F_1\per F_2\per F_3\per F_1$ is the identity. 

Set $F=\{x,\xi\}$, with $x$ of type $1$ and $\xi$ of type $4$. By \cref{MR0}, $\Pi^*(F)$ is the group induced on $\Res_\Delta(F)$ by the little projective group $G^\dagger$ of $\Delta$.   Then also $\Res_\Delta(\xi)$ is stabilised and it follows that  $\Pi^+(F)$ is the restriction of the stabiliser of $x$ in $\Pi^+(\xi)$ to $\Res_\Delta(F)$. All assertions now follow.    
\end{proof}

We note the following isomorphisms (for appropriate definitions of the unitary groups, similarly to their higher dimensional analogues): 
\[\mathsf{PO}_{5}(\K)\cong\overline{\mathsf{PS}}\mathsf{p}_4(\K); \hspace{.5cm} 
 \overline{\mathsf{PO}}^+_{6}(\K,\L)\cong \overline{\mathsf{U}}_4(\L/\K); \hspace{.5cm}
 \overline{\mathsf{PO}}^+_{8}(\K,\H)\cong \overline{\mathsf{U}}_4(\H),\]
 where $\L$ is a separable quadratic extension of $\K$, and $\H$ is a quaternion division algebra over $\K$.
 
We now turn to the planes and the rank $1$ residues. Recall the notion of a standard polarity in planes $\PG(2,\A)$, with $\A$ a quadratic alternative division algebra over $\K$ (see \cref{remarkstandardpolarity}). 
 
 \begin{lem}\label{planesofF4}
 If $F$ is a simplex of a building $\Delta:=\mathsf{F_4}(\K,\A)$, whose residue is isomorphic to $\PG(2,\B)$, with $\B\in\{\K,\A\}$, then $\Pi^+(F)\cong\PGL_3(\B)$ and $\Pi(F)\cong \PGL_3(\B)\rtimes 2$, where the extension is with a standard polarity if $\B=\A$, and with an ordinary orthogonal polarity if $\B=\K$. If $\Res_\Delta(F)$ is the rank~$1$ building over $\B$, then $\Pi^+(F)=\Pi(F)\cong \PGL_2(\B)$, if $\B=\K$ or, if the standard involution of $\B$ is trivial; otherwise $ \Pi^+(F)=\PGL^+_2(\B)$ and $\Pi(F)\cong \PGL_2(\B)$.  
 \end{lem}
 
 \begin{proof}
 Let $\alpha$ be a plane of $\mathsf{F_{4,1}}(\K,\A)$, let $\xi$ be a symplecton containing $\alpha$ and set $F=\{\alpha,\xi\}$ (then $F$ is a simplex of type $\{3,4\}$ of the corresponding building). 
 By \cref{MR0}, $\Pi^+(F)$ is induced by the little projective group $G^\dagger$ of the building. Since also $\Pi^+(\{\xi\})$ is induced by $G^\dagger$, we find that $\Pi^+(F)$ contains the stabiliser of $\alpha$ in the group $\Pi^+(\{\xi\})$ . The latter is determined in \cref{type1F4} (for the inseparable case), \cref{type4F4} and \cref{type4F4bis} (for the other cases). 
 These stabilisers are the full linear groups, hence, $\Pi^+(F)=\PGL_3(\K)$.  Similarly, if $F$ is a simplex of type $\{1,2\}$, then $\Pi^+(F)=\PGL_3(\A)$. Now we determine the general groups. So we exhibit a self-projectivity of length $3$. Let $F=\{x,L\}$, with $x$ a point incident with some line of $\mathsf{F}_{4,i}(\K,\A)$, with $i\in\{1,4\}$, and let $F'=\{x',L'\}$ be an opposite simplex. The residue of $F$ induces a plane $\pi$ in $E(x,x')$, viewed as polar space (the set of points symplectic with both $x$ and $x'$ via a symp through $L$). Similarly, the residue of $F'$ also induces a plane $\pi'$ in $E(x,x')$. Now, let $\zeta$ be a symp through $L$, intersecting $E(x,x')$ in a point $w\in\pi$. The plane $\beta'$ through $L'$ nearest to $\zeta$ is characterised by the property that every symp $\zeta'$ through $\beta'$ is special to $\zeta$.  It follows that such symps $\zeta'$ intersect $E(x,x')$ in a point symplectic to $E(x,x')\cap\zeta$. This implies that these intersections form the projection of $w$ to $\pi'$. Hence, we can consider a self-projectivity of length $3$ in a symp of $\mathsf{F}_{4,5-i}(\K,\A)$. If $i=4$, then it suffices to consider a self-projectivity in a parabolic polar space, and a direct computation shows that this is always a linear duality. If $i=1$, then,     according to \cref{standardpolarity} (which is valid for arbitrary $\A$ and not only for $\A=\O$), this is a standard polarity.
 
 Now suppose that $F$ has size $3$. Then the residue of $F$ is also a residue in a residue that is a projective plane, and hence, $\Pi^+(F)$ is the same as the special projectivity group in that projective plane. Completely similarly as in the previous paragraph, a self-projectivity of length $3$ can be equivalently seen as a self-projectivity of length $3$ of either (the lower residue of) a line or a planar line pencil. This reduces the determination of $\Pi(F)$ to the polar spaces $\mathsf{C_{3,1}}(\A,\K)$ and $\mathsf{B_{3,1}}(\K,\A)$. For the former, the results follow from   \cref{lowerlines}, \cref{NEpencil} and \cref{similitudes}; for the latter, it follows from \cref{embnm}. 
 \end{proof}
 
 We note the following isomorphisms: 
\begin{align*}
\mathsf{PO}_{3}(\K)\cong\PGL_2(\K); \hspace{.2cm} 
\overline{\mathsf{PO}}^+_{4}(\K,\L)\cong \PGL^+_2(\L/\K); \hspace{.2cm} \\
\overline{\mathsf{PO}}^+_{6}(\K,\H)\cong \PGL_2^+(\H),
\overline{\mathsf{PO}}^+_{10}(\K,\O)\cong \PGL_2^+(\O),
\end{align*}
where $\L$ is a separable quadratic extension of $\K$, $\H$ is a quaternion division algebra over $\K$, and $\O$ an octonion division algebra over $\K$.

\section{Conclusion for buildings of type $\mathsf{F_4}$}\label{conclusionF4}

We summarise the results obtained in the previous two sections for metasymplectic spaces in the following concluding theorem, including a tabular form.
\begin{theorem}

Let $ \mathsf{F_{4}}(\K,\A)$ be a building of type $\mathsf{F_4}$, with $\K$ a field, and $\A$ a quadratic alternative division algebra over $\K$. Let $F$ be a simplex with irreducible residue. Then $\Pi^+(F)$ and $\Pi(F)$ are as given in \emph{\cref{tableF4}}, where \begin{compactenum}[$*$] \item $\L$ denotes a separable quadratic extension of $\K$, \item $\H$ denotes a quaternion division algebra over $\K$, \item $\O$ denotes a Cayley (octonion) division algebra over $\K$,  \item $\K'$ denotes an inseparable extension of $\K$, with $\kar\K=2$, \item $\A'$ denotes a separable quadratic alternative division algebra over $\K$ with $\dim_{\K}\A'=d\in\{2,4,8\}$. \end{compactenum} 
\end{theorem}

\begin{table}[htb]
\renewcommand{\arraystretch}{1.4}
\begin{tabular}{|c|c|c|c|c|c|c|}\hline
Reference&$\Delta$& $\Res_\Delta(F)$ &  $\cotyp(F)$ & $\Pi^+(F)$&$\Pi(F)$& \\ \hline
\multirow{3}{*}{(A1)} & $\mathsf{F_4}(\K,\A)$ & $\mathsf{A_1}(\K)$ & $\{1\},\{2\}$ & $\PGL_2(\K)$ &$\PGL_2(\K)$&$\checkmark$\\
 & $\mathsf{F_4}(\K,\K)$ & $\mathsf{A_1}(\K)$ & $\{3\},\{4\}$ & $\PGL_2(\K)$ &$\PGL_2(\K)$&\\
&  $\mathsf{F_4}(\K,\A')$ & $\mathsf{A_1}(\A')$ & $\{3\},\{4\}$ &$\overline{\mathsf{PO}}^+_{d+2}(\K,\A')$ &$\overline{\mathsf{PO}}_{d+2}(\K,\A')$&\\ \hline
 \multirow{3}{*}{(A2)} &$\mathsf{F_4}(\K,\A)$ & $\mathsf{A_2}(\K)$& $\{1,2\}$  & $\PGL_3(\K)$ & $\PGL_3(\K)\rtimes 2$&\\ 
 & $\mathsf{F_4}(\K,\A)$ & $\mathsf{A_2}(\A)$ & $\{3,4\}$ & $\PGL_3(\A)$ &  $\PGL_3(\A)\rtimes2$&\\ 
 & $\mathsf{F_4}(\K,\K')$ & $\mathsf{A_2}(\K)$ & $\{1,2\}$& $\PGL_3(\K)$ & $\PGL_3(\K)\rtimes2$ &\\ \hline
\multirow{2}{*}{(B2)} & $\mathsf{F_4}(\K,\K)$ & $\mathsf{B_2}(\K,\K)$ & \multirow{2}{*}{$\{2,3\}$} & $\mathsf{PO}_5(\K)$ & $\mathsf{PO}_5(\K)$ &$\checkmark$ \\
& $\mathsf{F_4}(\K,\A')$ & $\mathsf{B_2}(\K,\A')$ & & $\overline{\mathsf{PO}}^+_{d+4}(\K,\A')$ & $\overline{\mathsf{PO}}^+_{d+4}(\K,\A')$ &$\checkmark$ \\ \hline
(C2)&$\mathsf{F_4}(\K,\K')$ & $\mathsf{B_2}(\K,\K')$ & $\{2,3\}$ & $\overline{\mathsf{PS}}\mathsf{p}_4(\K',\K)$ & $\overline{\mathsf{PS}}\mathsf{p}_4(\K',\K)$ & $\checkmark$ \\  \hline
\multirow{2}{*}{(B3)}&$\mathsf{F_4}(\K,\K)$ & $\mathsf{B_3}(\K,\K)$ & \multirow{2}{*}{$\{1,2,3\}$}& $\mathsf{PO}_7(\K)$ & $\mathsf{PO}_7(\K)$ &\\
& $\mathsf{F_4}(\K,\A')$ & $\mathsf{B_3}(\K,\A')$ & & $\overline{\mathsf{PO}}^+_{d+6}(\K,\A')$ & $\overline{\mathsf{PO}}_{d+6}(\K,\A')$ &  \\ \hline
\multirow{5}{*}{(C3)}& $\mathsf{F_4}(\K,\K)$ & $\mathsf{C_3}(\K,\K)$ & \multirow{5}{*}{$\{2,3,4\}$} &   $\overline{\mathsf{PS}}\mathsf{p}_6(\K)$& $\overline{\mathsf{PS}}\mathsf{p}_6(\K)$ & $\checkmark$ \\ 
&$\mathsf{F_4}(\K,\L)$ & $\mathsf{C_3}(\L,\K)$ &  &  $\overline{\mathsf{U}}_6(\L/\K)$ &  $\overline{\mathsf{U}}_6(\L/\K)$ & $\checkmark$ \\ 
&$\mathsf{F_4}(\K,\H)$ & $\mathsf{C_3}(\H,\K)$ &  & $\overline{\mathsf{U}}_6(\H)$ & $\overline{\mathsf{U}}_6(\H)$ & $\checkmark$ \\ 
&$\mathsf{F_4}(\K,\O)$ & $\mathsf{C_3}(\O,\K)$ &  &  $\overline{\mathsf{E}}_{7,3}^{28}(\O)$ &  $\overline{\mathsf{E}}_{7,3}^{28}(\O)$ &$\checkmark$ \\ & $\mathsf{F_4}(\K,\K')$ & $\mathsf{C_3}(\K',\K)$ &  &   $\overline{\mathsf{PS}}\mathsf{p}_6(\K',\K)$& $\overline{\mathsf{PS}}\mathsf{p}_6(\K',\K)$ & $\checkmark$ \\  \hline
\end{tabular}\vspace{.3cm}
\caption{Projectivity groups in the exceptional case $\mathsf{F_4}$\label{tableF4}}
\end{table}

The last column of \cref{tableF4} contains a checkmark ($\checkmark$), if \cite[Lemma~5.2]{Bus-Sch-Mal:24} automatically yields $\Pi^+(F)=\Pi(F)$; see also \cref{specialgeneralMSS}.  Note that there are other cases, for which $\Pi^+(F)=\Pi(F)$ generically holds, in contrast to the simply laced case, see \cite{Bus-Sch-Mal:24}.

\textbf{Acknowledgment.} The authors are grateful to Peter Abramenko, Theo Grundh\"ofer and Bernhard M\"uhlherr for some interesting discussions about some lemmata proved in this paper.


\end{document}